\RequirePackage{fix-cm}
\documentclass{svjour3}                     % onecolumn (standard format)
\smartqed  % flush right qed marks, e.g. at end of proof
\usepackage{graphicx}
\usepackage[pagebackref=true,breaklinks=true,colorlinks,citecolor=blue,bookmarks=false]{hyperref}
\usepackage{adjustbox}
\usepackage{multirow}
\usepackage{bbm}
\usepackage{comment}
\usepackage{siunitx}
\usepackage{relsize}
\usepackage{ifthen}

\usepackage[caption=false]{subfig}

\usepackage{graphics} % for pdf, bitmapped graphics files
\usepackage{rotating}
\usepackage{color}
\usepackage{enumerate}
\usepackage[T1]{fontenc}
\usepackage{psfrag}
\usepackage{epsfig} % for postscript graphics files
\usepackage{booktabs}
\usepackage{graphicx,url}
\usepackage{multirow}
\usepackage{array}
\usepackage{latexsym}
\usepackage{amsfonts}
\usepackage{amsmath}
\usepackage{amssymb}
\usepackage{xstring}
\usepackage{algorithm}
\usepackage{algpseudocode}
\usepackage{multirow}
\usepackage{xcolor}
\usepackage{prettyref}
\usepackage{flexisym}
\usepackage{bigdelim}
\usepackage{breqn} % load this last
\usepackage{listings}

\usepackage{enumitem}
\usepackage{xspace}
\usepackage{bm}
\usepackage{tikz}
\usetikzlibrary{matrix,calc}

\usepackage{mdwlist}

\makecompactlist{itemize}{stditemize}

\newcommand{\cf}{\emph{cf.}\xspace}

\newcommand{\bdmath}{\begin{dmath}}
\newcommand{\edmath}{\end{dmath}}
\newcommand{\beq}{\begin{equation}}
\newcommand{\eeq}{\end{equation}}
\newcommand{\bdm}{\begin{displaymath}}
\newcommand{\edm}{\end{displaymath}}
\newcommand{\bea}{\begin{eqnarray}}
\newcommand{\eea}{\end{eqnarray}}
\newcommand{\beal}{\beq \begin{array}{ll}}
\newcommand{\eeal}{\end{array} \eeq}
\newcommand{\beas}{\begin{eqnarray*}}
\newcommand{\eeas}{\end{eqnarray*}}
\newcommand{\ba}{\begin{array}}
\newcommand{\ea}{\end{array}}
\newcommand{\bit}{\begin{itemize}}
\newcommand{\eit}{\end{itemize}}
\newcommand{\ben}{\begin{enumerate}}
\newcommand{\een}{\end{enumerate}}

\newcommand{\calA}{{\cal A}}

\newcommand{\calC}{{\cal C}}
\newcommand{\calD}{{\cal D}}

\newcommand{\calF}{{\cal F}}

\newcommand{\calH}{{\cal H}}
\newcommand{\calI}{{\cal I}}

\newcommand{\calL}{{\cal L}}

\newcommand{\calN}{{\cal N}}

\newcommand{\calQ}{{\cal Q}}

\newcommand{\calS}{{\cal S}}
\newcommand{\calT}{{\cal T}}
\newcommand{\calU}{{\cal U}}
\newcommand{\calV}{{\cal V}}

\newcommand{\etal}{\emph{et~al.}\xspace}

\newcommand{\eg}{\emph{e.g.,}\xspace}
\newcommand{\ie}{\emph{i.e.,}\xspace}

\newcommand{\hide}[1]{}

\newcommand{\hiddenText}{{\color{gray} hidden text.}}
\newcommand{\hideWithText}[1]{\hiddenText}

\DeclareMathOperator*{\argmin}{arg\,min}

\newcommand{\norm}[1]{\left\| #1 \right\|}

\newcommand{\tran}{^{\mathsf{T}}}

\newcommand{\trace}[1]{\mathrm{tr}\left(#1\right)}

\newcommand{\rank}[1]{\mathrm{rank}\left(#1\right)}

\newcommand{\inv}{^{-1}}

\newcommand{\Real}[1]{ { {\mathbb R}^{#1} } }

\newcommand{\at}[1]{^{(#1)}}

\newcommand{\SOthree}{\ensuremath{\mathrm{SO}(3)}\xspace}

\newcommand{\scenario}[1]{{\smaller \sf#1}\xspace}

\newcommand{\blue}[1]{{\color{blue}#1}}

\newcommand{\red}[1]{{\color{red}#1}}

\newcommand{\linkToPdf}[1]{\href{#1}{\blue{(pdf)}}}
\newcommand{\linkToPpt}[1]{\href{#1}{\blue{(ppt)}}}
\newcommand{\linkToCode}[1]{\href{#1}{\blue{(code)}}}
\newcommand{\linkToWeb}[1]{\href{#1}{\blue{(web)}}}
\newcommand{\linkToVideo}[1]{\href{#1}{\blue{(video)}}}
\newcommand{\linkToMedia}[1]{\href{#1}{\blue{(media)}}}
\newcommand{\award}[1]{\xspace} % {{\red{#1}}} % omit awards

\renewcommand{\blue}[1]{{#1}}
\newcommand{\R}{\mathbb{R}}
\newcommand{\cA}{\mathcal{A}}

\newcommand{\cL}{\mathcal{L}}
\renewcommand{\norm}[1]{\left\lVert #1 \right\rVert}
\newcommand{\inprod}[2]{\left\langle #1, #2 \right\rangle}

\newcommand{\sym}[1]{\mathbb{S}^{#1}}
\newcommand{\calAadj}{\calA^{*}}

\newcommand{\bmat}{\left[ \begin{array}}
\newcommand{\emat}{\end{array}\right]}
\newcommand{\psd}[1]{\sym{#1}_{+}}
\newcommand{\pd}[1]{\sym{#1}_{++}}
\newcommand{\parentheses}[1]{\left(#1\right)}

\newcommand{\hatq}{\hat{q}}
\newcommand{\hattheta}{\hat{\theta}}

\newcommand{\abs}[1]{\left|#1\right|}

\newcommand{\hatx}{\hat{x}}
\newcommand{\hatX}{\widehat{X}}
\newcommand{\barX}{\overline{X}}

\newcommand{\name}{\scenario{STRIDE}}
\newcommand{\gnc}{\scenario{GNC}}
\newcommand{\namelong}{SpecTrahedRon Inexact projected gradient Descent along vErtices}

\newcommand{\cbrace}[1]{\left\{ #1\right\}}

\newcommand{\lbfgs}{L-BFGS}

\newcommand{\sdpnal}{\scenario{SDPNAL+}}
\newcommand{\sdpt}{\scenario{SDPT3}}
\newcommand{\mosek}{\scenario{MOSEK}}
\newcommand{\cdcs}{\scenario{CDCS}}

\newcommand{\sketchy}{\scenario{SketchyCGAL}}

\newcommand{\nchoosek}[2]{\left(\substack{#1 \\ #2} \right)}
\newcommand{\bard}{\bar{d}}
\newcommand{\mmom}{m_{\mathrm{mom}}}
\newcommand{\mloc}{m_{\mathrm{loc}}}
\newcommand{\pfeas}{\eta_p}
\newcommand{\dfeas}{\eta_d}
\newcommand{\gap}{\eta_g}
\newcommand{\subopt}{\eta_s}
\newcommand{\ee}[1]{\mathrm{e}{#1}}
\newcommand{\fmincon}{\texttt{fmincon}}
\newcommand{\rounding}{\texttt{rounding}}

\newcommand{\usphere}[1]{\calS^{#1}}
\newcommand{\manopt}{\scenario{Manopt}}
\newcommand{\tz}{\tilde{z}}
\newcommand{\tw}{\tilde{w}}

\newcommand{\fpop}{p}
\newcommand{\fpopstar}{\fpop^{\star}}
\newcommand{\fpsdp}{f_{\mathrm{P}}}
\newcommand{\fpsdpstar}{\fpsdp^{\star}}
\newcommand{\fdsdp}{f_{\mathrm{D}}}
\newcommand{\fdsdpstar}{\fdsdp^{\star}}
\newcommand{\xstar}{x^{\star}}
\newcommand{\Xstar}{X^{\star}}
\newcommand{\Wstar}{W^{\star}}
\newcommand{\xistar}{\xi^{\star}}
\newcommand{\zstar}{z^{\star}}

\newcommand{\ipgm}{iPGM}
\newcommand{\barz}{\bar{z}}
\newcommand{\barx}{\bar{x}}
\newcommand{\ystar}{y^{\star}}

\newcommand{\setpop}{\calF_{\mathrm{POP}}}
\newcommand{\setsdpp}{\calF_{\mathrm{P}}}
\newcommand{\nlp}{\texttt{nlp}}

\newcommand{\sgspgm}{sGS-APG}
\newcommand{\Sstar}{S^{\star}}
\newcommand{\calUadj}{\calU^*}
\newcommand{\etaproj}{\eta_{\mathrm{proj}}}
\newcommand{\tol}{\texttt{tol}}
\newcommand{\trinum}{\mathfrak{t}}

\newcommand{\sign}{\mathrm{sgn}}

\newcommand{\spadmm}{sPADMM}
\newcommand{\sos}{\texttt{sos}}

\newcommand\crule[3][black]{\textcolor{#1}{\rule{#2}{#3}}}
\newcommand{\redsquare}{\crule[red]{0.15cm}{0.15cm}}
\newcommand{\redcirc}{$\red{\bullet}$}
\newcommand{\redstar}{{\small $\red{\bigstar}$}}
\begin{document}

\title{An Inexact Projected Gradient Method with Rounding and Lifting by Nonlinear Programming for Solving Rank-One Semidefinite Relaxation of Polynomial Optimization
}
% \subtitle{Do you have a subtitle?\\ If so, write it here}

\titlerunning{{\name} for Solving Rank-One Semidefinite Relaxations}        % if too long for running head

\author{Heng Yang \and
        Ling Liang  \\
        Luca Carlone \and
        Kim-Chuan Toh
}


\institute{Heng Yang \at
              Laboratory for Information and Decision Systems \\
              Massachusetts Institute of Technology \\
              \email{hankyang@mit.edu}           %  \\
           \and
           Ling Liang \at
              Department of Mathematics \\
              National University of Singapore \\
              \email{liang.ling@u.nus.edu}
           \and
           Luca Carlone \at
              Laboratory for Information and Decision Systems \\
              Massachusetts Institute of Technology \\
              \email{lcarlone@mit.edu}  
          \and
          Kim-Chuan Toh \at
              Department of Mathematics and Institute of Operations Research and Analytics\\
              National University of Singapore \\
              \email{mattohkc@nus.edu.sg} 
          \and
          The first two authors contributed equally. \\
          Code available: \url{https://github.com/MIT-SPARK/STRIDE}
}

\date{\today}

\maketitle

%!TEX root = main.tex

\begin{abstract}
	We consider solving high-order semidefinite programming (SDP) relaxations of nonconvex polynomial optimization problems (POPs) that often admit degenerate rank-one optimal solutions. 
	Instead of solving the SDP alone, we propose a new algorithmic framework that blends \emph{local search} using the nonconvex POP into \emph{global descent} using the convex SDP. In particular, we first design a globally convergent \emph{inexact} projected gradient method (\ipgm) for solving the SDP that serves as the backbone of our framework. We then accelerate \ipgm\xspace by taking long, but \emph{safeguarded}, rank-one steps generated by fast nonlinear programming algorithms. We prove that the new framework is still globally convergent for solving the SDP. To solve the  \ipgm\xspace subproblem of projecting a given point onto the feasible set of the SDP, we design a two-phase algorithm with phase one using a symmetric Gauss-Seidel based accelerated proximal gradient method (sGS-APG) to generate a good initial point, and phase two using a modified limited-memory BFGS (\lbfgs) method to obtain an accurate solution. We analyze the convergence for both phases and establish a novel global convergence result for the modified \lbfgs\xspace that does not require the objective function to be twice continuously differentiable. We conduct numerical experiments for solving second-order SDP relaxations arising from a diverse set of POPs. Our framework demonstrates state-of-the-art efficiency, scalability, and robustness in solving degenerate rank-one SDPs to high accuracy, even in the presence of millions of equality constraints.
\end{abstract}

\keywords{
	Semidefinite programming \and 
	Polynomial optimization \and 
	Inexact projected gradient method \and 
	Rank-one solutions \and 
	Nonlinear programming \and 
	Degeneracy
}

\subclass{
	90C06 \and 
	90C22 \and 
	90C23 \and 
	90C55
}

% \tableofcontents

%!TEX root = main.tex
\section{Introduction}
\label{sec:introduction}

Let $p, h_1,\dots, h_l \in \Real{}[x]$ be real-valued multivariate polynomials in $x \in \Real{d}$, we consider the following equality constrained polynomial optimization problem\footnote{Although our framework can be generalized to problems with inequality constraints, here we focus on equality constraints to simplify the presentation.}
\begin{equation} \label{eq:pop}
\min_{x \in \Real{d}} \cbrace{p(x)\ |\ h_i(x)=0, i=1,\dots,l}. \tag{POP}
\end{equation}
Assuming problem \eqref{eq:pop} is feasible and bounded below, we denote $-\infty < \fpopstar < \infty$ as the global minimum and $\xstar$ as a global minimizer. Finding $(\fpopstar,\xstar)$ has applications in various fields of engineering and science~\cite{Lasserre09book-moments}. It is, however, well recognized that problem~\eqref{eq:pop} is NP-hard {in general (\eg it can model the binary constraint $x \in \{+1,-1\}$ via $x^2 = 1$)}.
% unless in some special cases (\eg~$\min_{x\tran x = 1} x\tran Q x $ for a symmetric matrix $Q$, has $(p^\star,\xstar)$ as the minimum eigenvalue and associated eigenvector of $Q$). 
As a result, solving convex relaxations of \eqref{eq:pop}, especially semidefinite programming (SDP) relaxations, has been the predominant method for finding $(\fpopstar,\xstar)$.
 % Following the same methodology, 
In this paper, we are interested in designing efficient, robust, and scalable algorithms for computing $(p^\star,\xstar)$ through SDP relaxations.

{\bf SDP Relaxations for Polynomial Optimization Problems}. Towards this goal, let us first give an overview of the celebrated Lasserre's moment/sums-of-squares (SOS) semidefinite relaxation hierarchy~\cite{Lasserre01siopt-LasserreHierarchy,Parrilo03MP-SOS}. Let $\kappa \geq 1$ be an integer such that $2\kappa$ is equal or greater than the maximum degree of $p,h_i,i=1,\dots,l$, and denote 
$v := [x]_\kappa$ as the standard vector of monomials in $x$ of degree up to $\kappa$. We build the \emph{moment matrix} $X := [x]_\kappa [x]_\kappa\tran$ that contains the standard monomials in $x$ of degree up to $2\kappa$. 
As a result,
$p$ and $h_i$'s can be written as linear functions in $X$.
For example, choosing $x = (x_1,x_2) \in \Real{2}$ and $\kappa=2$ leads to the following moment matrix that contains all monomials of $x$ up to degree $4$:
\beq\label{eq:examplemomentmat}
X = 
\bmat{cccccc}
	1 &  x_1  &  x_2  & x_1^2 & x_1 x_2 & x_2^2 \\
	x_1 &  x_1^2  &  x_1 x_2  & x_1^3 & x_1^2 x_2 & x_1 x_2^2 \\
	x_2 &  x_1 x_2  &  x_2^2  & x_1^2 x_2 & x_1 x_2^2 & x_2^3 \\
	x_1^2 &  x_1^3  &  x_1^2 x_2  & x_1^4 & x_1^3 x_2 & x_1^2 x_2^2 \\
	 x_1 x_2 &   x_1^2 x_2  &   x_1 x_2^2  & x_1^3 x_2 & x_1^2 x_2^2 & x_1 x_2^3 \\
	 x_2^2 &  x_1 x_2^2 &  x_2^3  & x_1^2 x_2^2 & x_1 x_2^3 & x_2^4 
\emat.
\eeq
We then consider two types of necessary linear constraints that can be imposed on $X$: (i) the entries of $X$ are not linearly independent, in fact the same monomial could appear in multiple locations of $X$ (\eg the monomial $x_1 x_2$ in \eqref{eq:examplemomentmat}); (ii) each constraint $h_i$ in the original \eqref{eq:pop} generates a set of equalities on $X$ of the form $h_i(x) [x]_{2\kappa - \deg(h_i)} = 0$, where $\deg(h_i)$ is the degree of $h_i$ (\ie if $h_i=0$, then $h_i \cdot x_1 = 0$, $h_i \cdot (x_1 x_2) = 0$ and so on).
% , as seen by the fact that if $h_i(x) = 0$ and $h_i$ has degree less than $2\kappa$, then $x_j \cdot h_i(x) = 0, \forall j = 1,\dots,d$. 
% Fortunately, Lasserre's hierarchy provides a systematic way of adding all linearly independent constraints on $X$~\cite{Lasserre01siopt-LasserreHierarchy}. 
Since $X$ is positive semidefinite by construction, one can therefore write the resulting semidefinite relaxation in standard \emph{primal} form as
\begin{equation}\label{eq:primalSDP}
\min_{X \in \sym{n}} \cbrace{\inprod{C}{X}\ |\ \calA(X) = b, \ X \succeq 0},\tag{P}
\end{equation}
where $n := \bard_\kappa \triangleq \nchoosek{d+\kappa}{\kappa}$ is the dimension of $[x]_\kappa$, $b \in \Real{m}$, $\calA: \sym{n} \rightarrow \Real{m}$ is a linear map $\calA(X) \triangleq (\inprod{A_i}{X})_{i=1}^m$ with $A_i \in \sym{n}$, that is onto and collects all independent linear constraints generated by the relaxation \blue{of $X =[x]_\kappa [x]_\kappa\tran$ to $X \succeq 0$}. Let $\calAadj: \Real{m} \rightarrow \sym{n}$ be the adjoint of $\calA$ defined as $\calAadj y \triangleq \sum_{i=1}^m y_i A_i$, then the Lagrangian dual of problem~\eqref{eq:primalSDP} reads
\begin{equation}\label{eq:dualSDP}
\max_{y \in \Real{m}, S \in \sym{n}} \cbrace{\inprod{b}{y}\ |\ \calAadj y + S = C,\ S \succeq 0},\tag{D}
\end{equation}
and admits an interpretation using sums-of-squares polynomials~\cite{Parrilo03MP-SOS}.
SDP relaxations~\eqref{eq:primalSDP}-\eqref{eq:dualSDP} correspond to the so-called \emph{dense} hierarchy, 
% of which a special case is using $\kappa=1$ for quadratically constrained quadratical programming. Although relaxations using $\kappa=1$ have been shown to be tight in a few application~\cite{Candes15SIReview-phaseretrieval,Abbe15TIT-exact,Briales17CVPR-3Dregistration,Shi21arXiv-PACE,Rosen19IJRR-sesync}, in general one needs a higher relaxation order to attain tightness. 
while many \emph{sparse} variants have been proposed to exploit the sparsity structure of $p$ and $h_i$'s to generate SDP relaxations with multiple blocks and smaller sizes, see~\cite{Waki2006,Wang21SIOPT-chordalTSSOS,Wang21SIOPT-TSSOS} and references therein. We remark that the algorithms developed in this paper can be extended to multiple blocks in a straightforward way to solve sparse relaxations. 

Throughout this paper, we assume that strong duality holds for \eqref{eq:primalSDP} and \eqref{eq:dualSDP}, and both \eqref{eq:primalSDP} and \eqref{eq:dualSDP} admit at least one solution.\footnote{Particularly, if a redundant ball constraint is included in the relaxation, then dual Slater and strong duality hold for Lasserre's moment relaxations of all orders \cite[Theorem 1]{Josz16OL-strongdualityLasserre}.}  We denote $\Xstar$ and $(\ystar,\Sstar)$ as an optimal solution for \eqref{eq:primalSDP} and \eqref{eq:dualSDP}, respectively. Then, the KKT condition for \eqref{eq:primalSDP} and \eqref{eq:dualSDP} given as
	\begin{equation}
	\label{eq-kkt-sdp}
		\cA(X) - b = 0,\quad \cA^*y + S - C = 0,\quad \inprod{X}{S} = 0,\quad X, S\succeq 0,	
	\end{equation}
	has $(\Xstar,\ystar,\Sstar)$ as one of its solutions. Moreover, it follows that $\fpopstar \geq \fpsdpstar = \fdsdpstar$, where $\fpsdpstar$ and $\fdsdpstar$ are the optimal values of~\eqref{eq:primalSDP} and~\eqref{eq:dualSDP}, respectively. The SDP relaxation is said to be \emph{tight} if $\fpopstar = \fpsdpstar$. In such cases, for any global minimizer $\xstar$ of \eqref{eq:pop}, the rank one lifting $X = [\xstar]_\kappa [\xstar]_\kappa\tran$ is a minimizer of~\eqref{eq:primalSDP}, 
	{and any rank one optimal solution  $\Xstar$ of \eqref{eq:primalSDP} corresponds to a global minimizer of \eqref{eq:pop}.}
	A special case of the relaxation hierarchy is Shor's relaxation ($\kappa=1$) for quadratically constrained quadratical programming problems~\cite{Shor90AOR-shorrelax,Luo10SPM-SDPQCQP}, of which applications and analyses have been extensively studied~\cite{Candes15SIReview-phaseretrieval,Abbe15TIT-exact,Briales17CVPR-3Dregistration,Shi21arXiv-PACE,Rosen19IJRR-sesync,Cifuentes17arxiv-localstability,Wang21MP-QCQPtightness,Burer19MP-exact}. Although $\kappa =1$ is certainly an interesting case, it cannot handle $p$ and $h_i$'s with degrees above $2$, and it is known to be not tight for many problems such as MAXCUT \cite{Goemans95JACM-maxcut}. {Moreover, it typically leads to SDPs with small $m$ that can often be handled well by existing SDP solvers \cite{Toh99OMS-sdpt3}.} Therefore, in this paper, we mainly focus on solving \emph{high-order} ($\kappa \geq 2$) relaxations that are more powerful in attaining tightness, as we describe below.
% However, the fact problem~\eqref{eq:primalSDP} admits rank-one solutions does not imply that algorithms for SDPs, such as interior point methods, can find them. For example, suppose~\eqref{eq:pop} has $K$ linearly independent global minimizers $\xstar_1,\dots,\xstar_K$, and let $X_1,\dots,X_K$ be their rank-one liftings $X_i = [\xstar_i][\xstar_i]\tran, i=1,\dots,K$. Then any $X \in \conv{X_1,\dots,X_K}$ is optimal for problem~\eqref{eq:primalSDP}. However, by the nature of interior point methods, a solution of maximum rank $K$ will be obtained~\cite{Deklerk06book-SDP}. A procedure for recovering multiple global minimizers is described in~\cite{Henrion05collection-optimality}.

% Adding inequalities to~\eqref{eq:pop} corresponds to adding more PSD constraints in~\eqref{eq:primalSDP} and can also be handled. For example, when the number os inequalities is small, a preferable strategy is to write an inequality constraint $h(x) \geq 0$ as $h(x)-s^2 = 0$ by introducing a slack variable $s$. 

% From a theoretical perspective, the moment/SOS hierarchy is extremely powerful. 
In the seminal work \cite{Lasserre01siopt-LasserreHierarchy}, Lasserre proved that $\fpsdpstar$ \emph{asymptotically} approaches $p^\star$ as $\kappa$ increases to infinity~\cite[Theorem 4.2]{Lasserre01siopt-LasserreHierarchy}. Later on, several authors showed that tightness indeed happens for a finite $\kappa$~\cite{Laurent07MP-momenttight,Nie13SIOPT-polynomial,Nie14MP-finiteconvergence}. Notably, Nie proved that, under the archimedean condition,
% (loosely speaking, problem~\eqref{eq:pop} admits compact feasible set), 
if constraint qualification, strict complementarity and second-order sufficient condition hold at every global minimizer of \eqref{eq:pop}, then the hierarchy converges at a finite $\kappa$~\cite[Theorem 1.1]{Nie14MP-finiteconvergence}. What is even more encouraging is that, empirically, for many important \eqref{eq:pop} instances, tightness {can be attained} at a very low relaxation order such as $\kappa = 2,3,4$. For instance, Lasserre~\cite{Lasserre01ICIPCO-SDPbinary} showed that $\kappa=2$ attains tightness for a set of 50 randomly generated MAXCUT problems; In the experiments of Henrion and Lasserre~\cite{Henrion03ACMTOMS-gloptipoly}, tightness holds at a small $\kappa$ for most of the~\eqref{eq:pop} problems in the literature; Doherty, Parrilo and Spedalieri~\cite{Doherty04PRA-quantumseparability} demonstrated that the second-order SOS relaxation is sufficient to decide quantum separability in all bound entangled states found in the literature of dimensions up to 6 by 6; Cifuentes~\cite{Cifuentes19arXiv-stls} designed and proved that a sparse variant of the second-order moment relaxation is guaranteed to be tight for the structured total least squares problem under a low noise assumption; Yang and Carlone applied a sparse second-order moment relaxation to globally solve a broad family of nonconvex computer vision problems~\cite{Yang19ICCV-wahba,Yang20neurips-onering,Yang20CVPR-perfectshape}
%,Yang21arxiv-certifiable}.

{\bf Computational Challenges in Solving High-order Relaxations}. The surging applications make it imperative to design computational tools for solving high-order (tight) SDP relaxations. However, high-order relaxations pose notorious challenges to existing SDP solvers. For this reason, in all applications above, the experiments were performed on solving~\eqref{eq:pop} problems of very small size, \eg the original~\eqref{eq:pop} has $d $ up to $20$ for dense relaxations. The computational challenges mainly come from two aspects. In the first place, in stark contrast to Shor's semidefinite relaxation, high order relaxations lead to SDPs with a large \blue{number} of linear constraints, even if the dimension of the original~\eqref{eq:pop} is small. Taking the binary quadratic programming as an example (\ie~$x_i \in \{+1,-1\}$ and $p$ is quadratic in \eqref{eq:pop}), when $d=60$, the number of equality constraints for $\kappa=2$ is $m=1,266,971$. In the second place, the resulting SDP is not only large, but also highly degenerate. Specifically, when the relaxation \eqref{eq:primalSDP} is tight and admits rank-one optimal solutions, it necessarily fails the primal nondegeneracy condition~\cite{Alizadeh97MP-degeneracy}, which is crucial for ensuring numerical stability and fast convergence of existing algorithms~\cite{Alizadeh98SIOPT-ipm,Zhao2010}. Due to these two challenges, our experiments show that \emph{no existing solver} can consistently solve semidefinite relaxations with rank-one solutions to a desired accuracy when $m$ is larger than $500,000$. In particular, interior point methods (IPM), such as~\sdpt~\cite{Toh99OMS-sdpt3} and~\mosek~\cite{mosek}, can only handle up to $m=50,000$ before {they} run out of memory on an ordinary workstation.\footnote{When chordal sparsity exists in the SDP such that a large positive semidefinite constraint can be decomposed into multiple smaller ones, IPMs can be more scalable \cite{Fukuda01siopt-sparsity,Zhang21MP-dualcliquetree}. However, for problems considered in this paper, there is no chordal sparsity.} First-order methods based on ADMM and conditional gradient method, such as~\cdcs~\cite{Zheng20MP-CDCS} and~\sketchy~\cite{Yurtsever21SIMDS-scalableSDP}, converge very slowly and cannot attain even modest accuracy for challenging instances. The best performing existing solver is~\sdpnal~\cite{Yang15MPC-sdpnalplus}, which employs an augmented Lagrangian framework where the inner problem is solved inexactly by a semi-smooth Newton method with a conjugate gradient method (SSNCG), and hence has very good scalability. The problem with {\sdpnal} is that, in the presence of large $m$ and degeneracy, solving the inexact SSN step involves solving a large and singular linear system of size $m \times m$, and CG becomes incapable of computing the search direction with sufficient accuracy. A tempting alternative approach, since~\eqref{eq:primalSDP} has rank-one optimal solutions, is to apply the low-rank factorization method of Burer and Monteiro~\cite{Burer03MP-BMlowrank},~\ie~solving the nonlinear optimization $\min_{V \in \Real{n \times r}} \{ \inprod{C}{VV\tran} \ |\ \calA(VV\tran) = b \}$ for some small $r$. We note that since $\nabla_V (\inprod{A_i}{VV\tran}) = 2A_i V$, if $\rank{V^\star} = 1$ at the optimal solution $V^\star$, then \blue{$\rank{2[ A_1 V^\star, \ldots,A_m V^\star]} \leq n \ll m$} no matter how large $r$ is. Therefore, $V^\star$ fails the linear independence constraint qualification (LICQ) and it may not be a KKT point, making it pessimistic for nonlinear programming solvers to find. In fact, successful applications of B-M factorization require $m \approx n$ so that the dual multipliers can be obtained in closed-form from nonlinear programming solutions~\cite{Rosen19IJRR-sesync,Rosen20wafr-scalableSDP,Boumal16neurips-BM}. Nevertheless, exploiting low rankness and nonlinear programming tools is a great source of inspiration, and as we will show, our solver {also exploits similar insights but} in a fashion that is not subject to the large $m$ and degeneracy of problem~\eqref{eq:primalSDP}.

{\bf Our Contribution.} In this paper, we contribute the first solver that can solve high-order SDP relaxations of \eqref{eq:pop} to high accuracy in spite of large $m$ and degeneracy.  {Our solver employs a globally convergent \emph{inexact projected gradient method} (iPGM) as the backbone for solving the convex SDP \eqref{eq:primalSDP}. Despite having favorable global convergence, the iPGM steps are typically \emph{short} and lead to slow convergence, particularly in the presence of degeneracy. Therefore, we blend short iPGM steps with \emph{long} rank-one steps generated by fast nonlinear programming (NLP) algorithms. The intuition is that, once we have achieved a sufficient amount of descent via iPGM, we can switch to NLP to perform local refinement to ``jump'' to a nearby rank-one point for rapid convergence. Due to this reason, we name our framework \emph{\namelong} (\name), which essentially strides along the rank-one vertices of the spectrahedron until the global minimum of the SDP is reached.}

Before presenting the mathematical details of {\name} in a general setup, let us
{first} introduce an informal version of {\name} by running an illustrative example.

%!TEX root = main.tex
% \Section{A Simple Univariate Example}
% \label{sec:simpleexample}

% In this section, we develop the high-level intuition and algorithmic framework of~\name~using a simple univariate example. 
\subsection{A Univariate Example}
\begin{figure}[ht!]
\centering
\includegraphics[width=\textwidth]{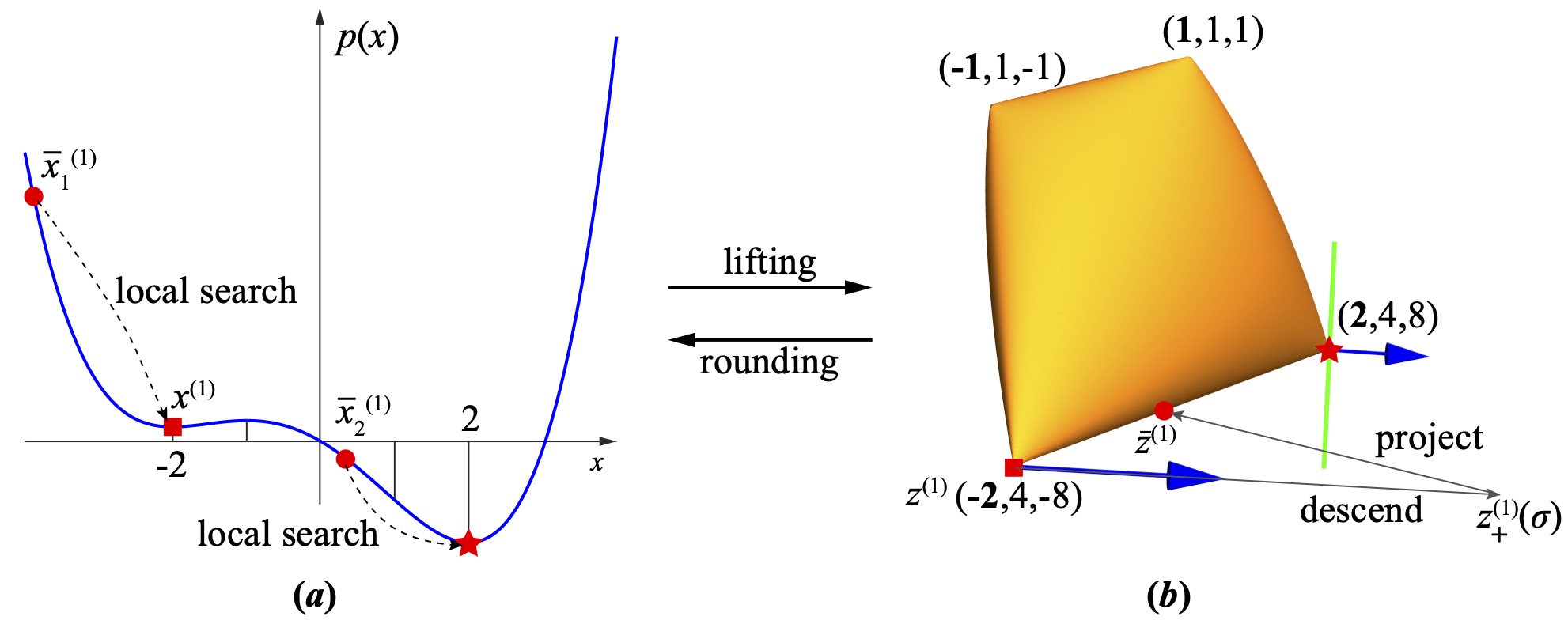}
\caption{Overview of our algorithm using a univariate example. (a) Plot of the univariate polynomial in~\eqref{eq:univariate}. (b) Visualization of the feasible set of the SDP~\eqref{eq:univariaterelax} arising from a second-order moment relaxation of~\eqref{eq:univariate}. ``\redsquare'': suboptimal local minimum of the POP and its corresponding (lifted) vertex in the SDP. ``\redstar'': globally optimal solutions of POP and SDP. ``\redcirc'': A descent point in SDP and corresponding rounded hypotheses in POP. Blue arrows in (b): direction along the negative gradient of the SDP objective function. Green line in (b): a hyperplane of constant objective going through the optimal solution.}
\label{fig:simpleexample}
\end{figure}

% {\bf A Univariate Example}. 
Consider minimizing a quartic polynomial subject to a single quartic equality constraint
\bea\label{eq:univariate}
\min_{x \in \Real{}} \ \ \cbrace{p(x):= x^4 + \frac{2}{3}x^3 - 8 x^2 - 8 x \ \middle|\  h(x) := (x^2-4)(x^2-1) = 0 },
\eea
where it is obvious that $x$ can only take four real values $\{+2, -2, +1, -1\}$ due to the constraint $h(x)$, among which $\xstar = 2$ attains the global minimum $\fpopstar = - \frac{80}{3} \approx -26.67$. A plot of $p(x)$ is shown in Fig.~\ref{fig:simpleexample}(a). However, let us discard our ``global'' view of $p(x)$ and design an algorithm to numerically obtain $(\xstar,p^\star)$ through its SDP relaxation. Towards this, according to our brief introduction of Lasserre's hierarchy, let us denote $v(x):=[1,x,x^2]\tran$ to be the vector of monomials in $x$ of degree up to $2$, and build the moment matrix $X:= v v\tran$. 
% \bea
% X := vv\tran = \bmat{ccc}
% 1 & x & x^2 \\
% x & x^2 & x^3 \\
% x^2 & x^3 & x^4
% \emat.
% \eea
$X$ has monomials of $x$ with degree up to 4, which we denote as $z_i := x^{i}, i = 1,\dots,4$. With $h(x)=0$, we further have $x^4 = 5x^2 - 4$,~\ie~$z_4 = 5 z_2 - 4$. Now we can write the cost function $p(x)$ as a linear function of $z$: $f(z) = z_4 + \frac{2}{3}z_3 - 8 z_2 - 8 z_1 = \frac{2}{3} z_3 - 3 z_2 - 8 z_1 - 4$. With this construction, the second-order relaxation of~\eqref{eq:univariate} reads
\bea\label{eq:univariaterelax}
\min_{z \in \Real{3}}\ \ \cbrace{f(z) := \frac{2}{3} z_3 - 3 z_2 - 8z_1 - 4\ \middle|\ X(z):=\bmat{ccc}
1 & z_1 & z_2 \\
z_1 & z_2 & z_3 \\
z_2 & z_3 & 5z_2 - 4\emat \succeq 0
}.
\eea
Problem~\eqref{eq:univariaterelax} allows us to explicitly visualize its feasible set, which is in general called a \emph{spectrahedron}. Using the fact that $X(z)\succeq 0$ if and only if the coefficients of its characteristic polynomial weakly alternate in sign~\cite[Proposition A.1]{Blekherman12book-sdpgeometry},
% we can describe the spectrahedron of~\eqref{eq:univariaterelax} with scalar inequalities
% \begin{equation}
% \begin{aligned}
% 3-6 z_2 \leq 0 \\
% - z_1^2 + 4 z_2^2 + 2 z_2 - z_3^2 - 4 \geq 0 \\
% 5 z_1^2 z_2 - 4 z_1^2 - 2 z_1 z_2 z_3 + z_2^3 - 5 z_2^2 + 4 z_2 + z_3^2 \leq 0
% \end{aligned}.
% \end{equation}
we plot the spectrahedron in Fig. \ref{fig:simpleexample}(b) using Mathematica \cite{Mathematica}. The four rank-one vertices are annotated with associated coordinates $(z_1,z_2,z_3)$, where $z_1 = x \in \{+2,-2,+1,-1\}$ corresponds to the four values that $x$ can take. The negative gradient direction $-\nabla f(z) = [8,3,-\frac{2}{3}]\tran$ is annotated by blue solid arrows. The hyperplane defined by $f(z)$ taking a constant value is shown as a green solid line in Fig.~\ref{fig:simpleexample}(b). Clearly, the rank-one vertex $\zstar = (2,4,8)$ (``\redstar''~in (b)) attains the minimum of the SDP \eqref{eq:univariaterelax}, which corresponds to $\xstar = 2$ (``\redstar''~in (a)) in POP~\eqref{eq:univariate}.

%!TEX root = main.tex

We now state our numerical algorithm that finds both $\xstar$ and $\zstar$.

{\bf Initialization}. We use a standard nonlinear programming (NLP) solver (\eg an interior point method interfaced via \fmincon~in Matlab) to solve problem~\eqref{eq:univariate}. We initialize at $x^{(0)} = -10$, and NLP converges to $x^{(1)} = -2$, with a cost $p(x^{(1)}) \approx -5.33$ (``\redsquare''~in~Fig.~\ref{fig:simpleexample}(a)). Note that $x^{(1)} = -2$ is a strict local minimum because choosing a Lagrangian multiplier of ``$0$'' for the constraint $h$ leads to $\nabla p(-2)=0$, $\nabla^2 p(-2) = 24$, \blue{thereby} satisfying first-order optimality and second-order sufficiency~\cite{Jorge2006}.

{\bf Lift and Descend}. The initialization step converges to a suboptimal solution. Unfortunately, it is very challenging, if not impossible, for NLP to be aware of this suboptimality, not to mention about escaping it. The SDP relaxation now comes into play. By lifting $x^{(1)}$ to $X = v(x^{(1)}) v(x^{(1)})\tran$ 
% $[1,x^{(1)},(x^{(1)})^2][1,x^{(1)},(x^{(1)})^2]\tran$ 
according to the SDP relaxation, we obtain the suboptimal rank-one vertex ``\redsquare''~in Fig.~\ref{fig:simpleexample}(b) with $z^{(1)} = [-2,4,-8]\tran$.
% \bea \nonumber
% z^{(1)} = \bmat{c} -2 \\ 4 \\ -8 \emat,\quad 
% X(z^{(1)}) = \bmat{ccc}
% 1 & -2 & 4 \\
% -2 & 4 & -8 \\
% 4 & -8 & 16
% \emat.
% \eea
Although it is challenging to descend from $x^{(1)}$ in POP~\eqref{eq:univariate}, descending from $z^{(1)}$ in the convex SDP~\eqref{eq:univariaterelax} is relatively easy. We first take a step along the negative direction of the gradient to arrive at $z^{(1)}_{+}(\sigma) := z^{(1)} - \sigma \nabla f(z^{(1)}) = [-2+8\sigma,4+3\sigma,-8-\frac{2}{3}\sigma]\tran$ for a given step size $\sigma > 0$, and then project $z^{(1)}_{+}$ back to the spectrahedron in Fig.~\ref{fig:simpleexample}(b), denoted as $\barz^{(1)}$. Note that the projection step solves the problem $\min_{z \in \Real{3}} \{ \| X(z) - X(z^{(1)}_{+}) \|^2 \ |\ X(z) \succeq 0 \}$, which is itself a quadratic SDP. For this small-scale example, we can compute the projection exactly (using \eg IPMs). Taking $\sigma = 5$ yields $z^{(1)}_+ = [38,19,-\frac{34}{3}]\tran$ and $\barz^{(1)} = [-0.35,3.99,-1.67]\tran$,
% \bea \nonumber
% &
% z^{(1)}_+ = \bmat{c}
% 38 \\
% 19 \\
% -\frac{34}{3}
% \emat\!\!, \quad
% X(z^{(1)}_+) = \bmat{ccc}
% 1 & 38 & 19 \\
% 38 & 19 & -\frac{34}{3} \\
% 19 & -\frac{34}{3} & 91
% \emat\!\!,\\
% &
% \barz^{(1)} = \bmat{c}
% -0.35 \\
% 3.99 \\
% -1.67
% \emat\!\!, \quad
% X(\barz^{(1)}) = \bmat{ccc}
% 1 & -0.35 & 3.99 \\
% -0.35  & 3.99 & -1.67 \\
% 3.99  & -1.67  & 15.97
% \emat\!\!,\nonumber
% \eea
and the projection $\barz^{(1)}$ is plotted in Fig.~\ref{fig:simpleexample}(b) as ``\redcirc''.

{\bf Rounding and Local Search}. The step we just performed is called \emph{projected gradient descent}, and it is well known that iteratively performing this step guarantees \blue{convergence} to the SDP optimal solution \cite{Bertsekas99book-nlp}. However, the drawback is that it typically requires many iterations for the convergence to happen. Now we will show that, by leveraging local search back in POP \eqref{eq:univariate}, we can quickly arrive at the rank-one optimal point. The intuition is, from Fig.~\ref{fig:simpleexample}(b), we can observe that $\barz^{(1)}$ has moved ``closer'' to $\zstar$, and hence maybe $\barz^{(1)}$ already contains useful local information about $\zstar$. Particularly, let us perform a spectral decomposition of the top-left $2\times 2$ block of $X(\barz^{(1)})$:
\bea\nonumber 
&
\bmat{cc}
1 & -0.35 \\
-0.35  & 3.99 
\emat = \lambda_1 v_1 v_1\tran + \lambda_2 v_2 v_2\tran, \quad \text{with} \\
&
v_1 = 0.12\bmat{c}
1 \\
-8.59
\emat,
\ v_2 = 0.99\bmat{c}
1 \\
0.12
\emat,
\ \lambda_{1,2} = (4.04, 0.96),\nonumber
\eea
where we have written $v_{1,2}$ by normalizing their leading entries as $1$, and ordered them such that $\lambda_1 \geq \lambda_2$. We then generate two \emph{hypotheses} $\barx^{(1)}_{1,2} = (-8.59,0.12)$ (``\redcirc''~in Fig.~\ref{fig:simpleexample}(a)) and use NLP to perform local search starting from the hypotheses $\barx^{(1)}_{1,2}$, respectively. Surprisingly, although $\barx^{(1)}_1$ still converges to $x^{(1)} = -2$, $\barx^{(1)}_2$ converges to $\xstar = 2$. We call this step \emph{rounding and local search}.

{\bf Lift and Certify}. Now that NLP has visited both $x^{(1)}$ and $\xstar$, with $\xstar$ attaining a lower cost, we are certain that $\xstar$ is at least a better local minimum. To \emph{certify} the global optimality of $\xstar$, we need the SDP relaxation again. We lift $\xstar$ to $\zstar = [2,4,8]\tran$, and try to descend from $\zstar$ by performing another step of projected gradient descent.
% first going along the negative gradient direction and then projecting the trial point back to the spectrahedron. 
However, because $\zstar$ is optimal, this step ends up back to $\zstar$ again, \ie~no more descent is possible \cite{Bertsekas99book-nlp}. Therefore, we conclude that $\zstar$ is optimal for SDP~\eqref{eq:univariaterelax}. Because $X(\zstar)$ is rank one, and $p(\xstar) = f(\zstar)$, we can say that $\xstar$ is indeed globally optimal for the nonconvex POP~\eqref{eq:univariate} (\ie a numerical \emph{certificate} of global optimality is obtained for $\xstar$).

Despite the simplicity of this example, it clearly demonstrates the stark contrast between our algorithm {\name} and existing approaches. While all existing methods solve the SDP relaxation \emph{independently} from the original POP (\ie they relax the POP into an SDP and solve the SDP without \emph{computationally} revisiting the POP), our method alternates between \emph{local search} in POP and \emph{global descent} in SDP, and the connection between the POP and the SDP is established by \emph{lifting} and \emph{rounding}. The advantages for doing so are threefold. (i) \emph{Scalability}: Although the SDP may have millions of variables ($n$ and $m$ can grow rapidly), the original POP has only hundreds or thousands of variables (recall that $d$ is moderate) and an NLP solver can perform local search quickly (typically in a negligible amount of time compared to solving the SDP); (ii) \emph{Degeneracy}: When the optimal SDP solution is rank one and degenerate, convergence of existing scalable SDP solvers (such as first-order methods) is typically very slow or even unachievable (\cf results in Section \ref{sec:numexps}). However, the original POP is much less sensitive to the degeneracy. Once the SDP iterate gets close to being optimal, NLP can quickly arrive at the global optimal solution. (iii) \emph{Warm-starting}: Unlike IPMs, for which designing good initialization is relatively challenging, our algorithm easily benefits from warm-starting. As we will show in Section~\ref{sec:numexps}, in many practical engineering applications, there exist powerful heuristics that find the globally optimal POP solution with high probability of success. In this case, our algorithm only needs to perform the last step of lifting and certification. 

However, several questions need to be answered in order for our algorithm to be general and practical for large-scale problems.
\begin{question}
	\label{qn-1}
	Since the projection onto a spectrahedron is itself an SDP and typically cannot be done exactly when $n,m$ are large, how can one design a \emph{globally convergent} {solver for the original SDP \eqref{eq:primalSDP}} that can tolerate \emph{inexactness} of the projection subproblem?
\end{question}
\begin{question}
	\label{qn-2}
	With rounding and local search in the loop, is it still possible for the algorithm to be globally convergent?	
\end{question}
\begin{question}
	\label{qn-3}
	{How can one design a \emph{scalable and efficient} numerical algorithm to handle millions of constraints (\ie $m$) when computing the (inexact) projection onto the spectrahedron?}
\end{question}
\begin{question}
	\label{qn-4}
	{While our univariate example shows that this framework} works on the simple example above, will it work on more complicated problems arising in real world applications?
\end{question}
% \begin{question}
% 	\label{qn-5}
% 	Under what conditions can we guarantee that rounding and local search can find a better local minimum?	
% \end{question}

Section~\ref{section-iPGM} answers Question \ref{qn-1}, where we design an \emph{inexact} projected gradient method (\ipgm) for solving a generic SDP pair~\eqref{eq:primalSDP}-\eqref{eq:dualSDP}, prove its global convergence, and provide complexity analysis. Although our results are inspired by~\cite{Jiang2012}, several nontrivial extensions are made. In answering Question \ref{qn-2} (Section \ref{sec:nlpacceleration}), we incorporate rounding and local search into \ipgm~and formally present {\name}. In a nutshell, \name~follows the globally convergent trajectory driven by \ipgm, but simultaneously probes long, but safeguarded, rank-one vertices of the spectrahedron \blue{generated} from solutions of NLP, to seek rapid descent and convergence. Notably, we prove that, even with rounding and local search in the loop, \name~\blue{is guaranteed} to converge to the optimal solution of~\eqref{eq:primalSDP}-\eqref{eq:dualSDP}. In Section~\ref{sec:apgsubproblem}, we focus on solving the critical subproblem of projecting a given symmetric matrix onto the feasible set of~\eqref{eq:primalSDP} and provide an efficient answer to Question \ref{qn-3}. We propose a two-phase algorithm where phase one uses a \emph{symmetric Gauss-Seidel based accelerated proximal gradient method} (\sgspgm) to generate a good initial point, and phase two applies a modified limited-memory BFGS (L-BFGS) method to compute an accurate solution. This two-phase algorithm is simple, robust, easy to implement, and can scale to very large SDP problems with millions of constraints. Additionally, for the modified L-BFGS algorithm, we establish a novel convergence result \blue{where the objective function does not need} to be at least twice continuously differentiable. Finally, we answer Question \ref{qn-4} by providing extensive numerical results in Section~\ref{sec:numexps}. We apply~\name~to solve (dense and sparse) second-order moment relaxations arising from a diverse set of POP problems and demonstrate its superior performance. We observe that \name~is the only solver that can consistently solve rank-one semidefinite relaxations to high accuracy (\eg~KKT residuals below $1\ee{-9}$) and it is up to {$1$-$2$ orders of magnitude faster than existing SDP solvers}. 
% On SDP relaxations with high-rank solutions such as MAXCUT, \name~also compares favorably to existing solvers. 
% Thus, the high efficiency, scalability and robustness of \name provides a strong affirmative answer to. 
% Question \ref{qn-5} is interesting but theoretically challenging, and we hope more future research can be done to explain the effectiveness of \name.

{\bf Notation}.
We use $\mathcal{X}$ and $\mathcal{Y}$ to denote finite dimensional Euclidean spaces. Let $\sym{n}$ be the space of real symmetric $n\times n$ matrices, and $\psd{n}$ (resp. $\pd{n}$) be the set of positive semidefinite (resp. definite) matrices. We also write $X \succeq 0$ (resp. $X \succ 0$) to indicate that $X$ is positive semidefinite (resp. definite) when the dimension of $X$ is clear. We use $\usphere{d-1}:=\{ x \in \Real{d}\mid \norm{x} = 1 \}$ to denote the $d$-dimensional unit sphere. For $x \in \Real{n}$, $\norm{x} = \sqrt{x\tran x}$ is the standard $\ell_2$ norm. For $X \in \Real{m \times n}$, $\norm{X} = \sqrt{\trace{X\tran X}}$ denotes the Frobenius norm. $\norm{x}_{\calH} := \sqrt{\inprod{x}{\calH x}}$ denotes the $\calH$-weighted norm for a positive semidefinite mapping $\calH$.
{We use $\mathcal{I}$ to denote the identity map from $\mathbb{S}^n\rightarrow \mathbb{S}^n$.}
% Let $\mathcal{I}$ denote the identity map from $\mathbb{S}^n\rightarrow \mathbb{S}^n$. 
We use $\delta_{\calC}(\cdot)$ to denote the indicator function of set $\calC$. For a positive integer $d$, we use $\trinum(n):= \frac{d(d+1)}{2}$ to denote the $d$-th triangle number, and we use $\bard_{\kappa} := \nchoosek{n+\kappa}{\kappa}$ for some positive integer $\kappa$, which is particularly helpful when describing the number of monomials in $x \in \Real{d}$ of degree up to $\kappa$.
% (the $\kappa$-th dense moment relaxation of a POP with $d$ variables leads to an SDP with matrix $\bard_\kappa$). 
For $x \in \Real{d}$, we use $[x]_{\kappa} \in \Real{\bard_\kappa}$ to denote the full set of standard monomials of degree up to $\kappa$ (as already used in the introductory paragraphs).

% \bea
% \min_{\vxx_i \in \calX} \sum_{k=1}^N \rho (r(\vxx_i , \vp_{i,k}^a, \vq_{i,k}^b (\vtheta) ))
% \eea

% \bea
% \calL(\theta) = \sum_{k=1}^N \rho (r(\vxx_i(\vtheta) , \vp_{i,k}^a, \vq_{i,k}^b (\vtheta) ))
% \eea
% {\bf Our contributions}
% \begin{itemize}
% 	\item Inspired by \sdpnal~but make significant customizations for POPs.
% 	\item Inexact (A)PG framework, with global convergence and iteration complexity analysis
% 	\item In the subproblem, we use sGS-APG + L-BFGS, simple, easy to implement, scalable, and empirically robust.
% 	\item For L-BFGS, we establish a convergence result when the objective is not twice continuously differentiable.
% 	\item We are the first one to algorithmically connect the SDP and POP by using NLP methods inside the APG framework to accelerate convergence. Surprisingly, even with NLP in the loop, our algorithm is still guaranteed to be globally convergent.
% 	\item Highlight the extensive computation experiments.
% \end{itemize}

%!TEX root = main.tex

\section{An inexact Projected Gradient Method for SDPs}
\label{section-iPGM}
In this section, we describe an inexact Projected Gradient Method (iPGM) for solving the SDP \eqref{eq:primalSDP} and analyze its global convergence properties. Accelerating iPGM 
\blue{via rounding and lifting by} nonlinear programming will be presented in Section \ref{sec:nlpacceleration}.

Algorithm \ref{alg-iPGM} presents the pseudocode for iPGM. Denoting the feasible set of the primal SDP \eqref{eq:primalSDP} as
\begin{equation*}
	\mathcal{F}_{\mathrm{P}} := \left\{ X\in \mathbb{S}^n\;\middle |\; \mathcal{A}(X) = b,\; X\in \mathbb{S}_+^n\right\},
\end{equation*}
the $k$-th iteration of iPGM first moves along the direction of the negative gradient (\ie $-C$) with step size $\sigma_k > 0$, and then performs an inexact projection of the trial point $X^{k-1} - \sigma_k C$ onto the primal feasible set $\setsdpp$ (\cf \eqref{eq-iPGM-proj}), with inexactness conditions given in \eqref{eq-iPGM-inexact-conditions}. In Section \ref{sec:apgsubproblem}, we will present {efficient algorithms} that can fulfill the inexactness conditions in \eqref{eq-iPGM-inexact-conditions}.

\begin{algorithm}
	\caption{An inexact projected gradient method (\ipgm) for SDP \eqref{eq:primalSDP}.}
	\label{alg-iPGM}
	\begin{algorithmic}
		\item[{\bf Input}:] Initial points $(X^0, y^0, S^0)\in \mathbb{S}^n_+ \times \mathbb{R}^m \times \mathbb{S}^n_+$, a nondecreasing positive sequence $\{\sigma_k\}$ and a nonnegative sequence $\{\varepsilon_k\}$ such that $\{k\epsilon_k\}$ is summable.
		\item[{\bf Repeat} For $k \geq 1$:] 
		\item[\quad] Compute 
		\begin{equation}
		\label{eq-iPGM-proj}	
		X^{k} \approx \Pi_{\mathcal{F}_{\mathrm{P}}}\left(X^{k-1} - \sigma_k C\right),\quad X^{k} \in \mathbb{S}_+^n
		\end{equation}
		with associated dual pair $(y^{k}, S^{k})\in \mathbb{R}^m\times \mathbb{S}_+^n$ such that the following conditions hold:
		\begin{equation}
		\label{eq-iPGM-inexact-conditions}
		\begin{aligned}
			\norm{\mathcal{A}(X^{k}) - b} \leq &\; \varepsilon_k, \\
				\mathcal{A}^*y^k + S^k - C - \frac{1}{\sigma_k}\left(X^{k} - X^{k-1}\right)= &\; 0, \\
				\inprod{X^{k}}{S^{k}} = &\;0.
			\end{aligned}
		\end{equation}
		\item[{\bf Until}:] Termination conditions are met (\cf Section \ref{sec:numexps}).
		\item[{\bf Output}:] $(X^k,y^k,S^k)$.
	\end{algorithmic}
\end{algorithm}

%\begin{figure}[ht!]
%\centering
%\centerline{\fbox{\parbox{\textwidth}{
%	{\bf Algorithm iPGM: An inexact proximal gradient method for (P).} \\[5pt]
%	Given $(X^0, y^0, S^0)\in \mathbb{S}^n_+ \times \mathbb{R}^m \times \mathbb{S}^n_+$, a nondecreasing positive sequence $\{\sigma_k\}$, and a nonnegative sequence $\{\varepsilon_k\}$ such that $\{k\epsilon_k\}$ is summable. 	
%	\begin{description}
%		\item[{\bf Repeat} for $k \geq 1$.] \quad \\
%		Compute 
%		\begin{equation}
%		\label{eq-iPGM-proj}	
%		X^{k} \approx \Pi_{\mathcal{F}_{\mathrm{P}}}\left(X^{k-1} - \sigma_k C\right),\quad X^{k} \in \mathbb{S}_+^n
%		\end{equation}
%		with associated dual pair $(y^{k}, S^{k})\in \mathbb{R}^m\times \mathbb{S}_+^n$ such that the following conditions hold:
%		\begin{equation}
%			\label{eq-iPGM-inexact-conditions}
%			\begin{aligned}
%				\norm{\mathcal{A}(X^{k}) - b} \leq &\; \varepsilon_k, \\
%				\mathcal{A}^*y^k + S^k - C - \frac{1}{\sigma_k}\left(X^{k} - X^{k-1}\right)= &\; 0, \\
%				\inprod{X^{k}}{S^{k}} = &\;0.
%			\end{aligned}
%		\end{equation}
%		\item[{\bf Until} termination conditions are met.]
%	\end{description}
%}}}	
%\caption{An inexact proximal gradient method. \label{alg-iPGM}}
%\end{figure}

The following theorem states the convergence properties of Algorithm \ref{alg-iPGM} (iPGM).
\begin{theorem}
	\label{thm-iPGM}
	Let $\{(X^k, y^k, S^k)\}$ be any sequence generated by Algorithm \ref{alg-iPGM}, and suppose that there exists a constant $M>0$ such that $\norm{y^k}\leq M$ for all $k\geq 0$. Then, for all $k\geq 1$, it holds that 
	\begin{equation*}
		\left\{
			\begin{array}{l}
				\displaystyle -\frac{\norm{y^*}k\varepsilon_k}{k} \leq \inprod{C}{X^k - \Xstar} \leq \frac{1}{k}\left( \frac{1}{2\sigma_0}\norm{X^0 - \Xstar}^2 + 2M\sum_{i = 1}^k i\varepsilon_i\right), \\
				\displaystyle \norm{\mathcal{A}(X^k) - b} \leq \varepsilon_k \leq O\left(\frac{1}{k}\right), \\ 	
				\displaystyle \norm{\mathcal{A}^*y^k + S^k - C} \leq O\left( \frac{1}{\sqrt{k\sigma_k}}\right),
			\end{array}
		\right.
	\end{equation*}
	where $\Xstar$ is an optimal solution of the primal problem \eqref{eq:primalSDP} and $(\ystar, \Sstar)$ is an optimal solution of the dual problem \eqref{eq:dualSDP}. In particular, if one chooses $\sigma_k \geq O(\sqrt{k}) > 0$ for all $k\geq 1$, then we have 
	\begin{equation*}
		\max\left\{ \left\lvert \inprod{C}{X^k - \Xstar}\right\rvert, \norm{\mathcal{A}(X^k) - b}, \norm{\mathcal{A}^*y^k + S^k - C}\right\} \leq O\left(\frac{1}{k}\right).
	\end{equation*}
\end{theorem}

The proof of Theorem \ref{thm-iPGM} is given in Appendix \ref{appendix-proof-iPGM}. Although an accelerated version of Algorithm \ref{alg-iPGM} and its convergence analysis have been studied in \cite{Jiang2012}, Theorem \ref{thm-iPGM} and its proof are new, {to the best of our knowledge}. By the presented results, if one chooses $\varepsilon_k$ to be sufficiently small and $\sigma_k$ to be sufficiently large, then the convergence of primal and dual infeasibilities can be {as fast as $O(1/k)$}. However, a small $\varepsilon_k$ or a large $\sigma_k$ will make the projection problem in \eqref{eq-iPGM-proj} \blue{more} difficult to solve. Hence, for better overall efficiency, one may choose $\varepsilon_k$ and $\sigma_k$ dynamically to balance the convergence speed of Algorithm \ref{alg-iPGM} and the efficiency of the (inexact) projection onto $\calF_{\mathrm{P}}$.

%!TEX root = main.tex
\section{{\name}: Accelerating \ipgm~by Nonlinear Programming}
\label{sec:nlpacceleration}
Despite being globally convergent, it may take many iterations for Algorithm \ref{alg-iPGM} to converge to a solution of high accuracy, especially when the SDP \eqref{eq:primalSDP} has large scale and the optimal solution $\Xstar$ is low-rank and degenerate. Therefore, in this section, we propose to accelerate Algorithm \ref{alg-iPGM} by rounding and local search, just as what we have shown in the simple numerical example in Fig.~\ref{fig:simpleexample} of Section~\ref{sec:introduction}. The intuition here is simple: when the SDP relaxation \eqref{eq:primalSDP} is exact, there exist rank-one optimal solutions and they may be computed much more efficiently by performing local search in the low-dimensional \eqref{eq:pop} with proper initialization.

%Recall that the primal SDP problem \eqref{eq:primalSDP} can be rewritten as the following constrained minimization problem:
%\begin{equation}
%	\label{prob:constrained-SDP}
%	\min_{X \in \sym{n}}  f(X) + g(X) ,
%\end{equation}
%where $f:\mathbb{S}^{n}\rightarrow \mathbb{R}$ defined as $f(X) = \left\langle C, X \right \rangle$ is a linear function with $ L_f $-Lipschitz gradient and $g(X) = \delta_{\setsdpp}(X)$ with $\setsdpp = \left\{ X\in \mathbb{S}^n\mid \mathcal{A}(X) = b, \; X\succeq 0\right\}$. Note that $ L_f $ could be any positive real number. Thus, for any positive and \emph{nondecreasing} sequence $ \{\sigma_k\} $ if we apply Algorithm \ref{alg-iPGM} with $\calH_k= \frac{1}{\sigma_k}\mathcal{I}$ to \eqref{prob:constrained-SDP} directly, Theorem \ref{thm-A} and Corollary \ref{corollary-PGM} will guarantee the convergence under certain mild conditions. However, the convergence could be slow even if we adjust the step size $ \sigma_k $ carefully, especially when the optimal solution is degenerate.

With this intuition, we now develop the details of the acceleration scheme. At each iteration with current iterate $ X^{k-1}\in \mathbb{S}^n_+ $, we first follow \eqref{eq-iPGM-proj} in Algorithm \ref{alg-iPGM} and compute the projection
\[
	\barX^{k}\approx \Pi_{\setsdpp}\left( X^{k-1} - \sigma_k C \right)
\]
under the inexactness conditions described in \eqref{eq-iPGM-inexact-conditions}, where $ \sigma_k>0 $ is a given step size. Then by the analysis in Theorem \ref{thm-iPGM}, $ \overline{X}^{k}$ is guaranteed to make certain progress towards optimality. However, the point $ \overline{X}^{k} $ may not be a promising candidate for the next iteration since it may not attain rapid convergence.
% and the rank of $\barX^k$ may not be one or not even small. 
Motivated by the success of rounding and local search in the simple example in Fig.~\ref{fig:simpleexample}, we propose to compute a potentially better candidate based on $ \overline{X}^{k} $ via the following steps:
\begin{enumerate}
	\item\label{item:rounding} {\bf (Rounding)}. Let $\barX^{k} = \sum_{i=1}^n \lambda_i v_i v_i\tran$ be the spectral decomposition of $\barX^{k}$ with $\lambda_1 \geq \dots \geq \lambda_n$ in nonincreasing order. Compute $r \geq 1$ hypotheses for the \eqref{eq:pop} from the leading $r$ eigenvectors $v_1,\dots, v_r$   
	\bea
		\bar{x}^{k}_{i} = \rounding(v_i), \quad i = 1,\dots,r,
	\eea
	where the function~\rounding~can be problem-dependent and we provide examples in the numerical experiments in Section~\ref{sec:numexps}. Generally, if the SDP~\eqref{eq:primalSDP} comes from a dense relaxation, and the feasible set of the original~\eqref{eq:pop}, denoted as $\setpop$, is simple to project, then we can design~\rounding~to be
	\bea
	v_i \leftarrow \frac{v_i}{v_i[1]}, \quad \bar{x}^{k}_{i} \leftarrow \Pi_{\setpop}(v_i[x]),
	\eea
	where one first normalizes $v_i$ such that its leading entry $v_i[1]$ is equal to $1$, and then projects its entries corresponding to order-one monomials,~\ie~$v_i[x]$, to the feasible set of~\eqref{eq:pop}. Note that the rationale for designing this rounding method is that practical POP applications often involve simple constraints such as binary~\cite{Goemans95JACM-maxcut}, unit sphere~\cite{Yang19ICCV-wahba}, and orthogonality~\cite{Yang20neurips-onering,Briales17CVPR-3Dregistration}, whose projection maps are simple (\eg if $x$ is binary, then $\Pi_{\setpop}(v_i[x]) = \sign\parentheses{v_i[x]}$ just takes the sign of each entry of $v_i[x]$). If $\setpop$ is not easy to project, then we omit the projection, in which case the local search will start from an infeasible initialization.
	\item {\bf (Local Search)}. Apply a local search method for the POP (as an NLP) with the initial point chosen as $ \bar{x}^{k}_i $ for each hypothesis $i=1,\dots,r$. Denote the solution of each local search as $ \hat{x}^{k}_i $, with associated objective value $p(\hatx^{k}_i)$, choose the best local solution with minimum objective value. Formally, we have
	\bea
	\hatx_i^{k} = \nlp(\barx^{k}_i),\quad i=1,\dots,r, \quad \hatx^{k} = \argmin_{\hatx^{k}_i, i=1\dots,r} p(\hatx^{k}_i).
	\eea

	\item\label{item:lifting} {\bf (Lifting)}. Perform a rank-one lifting of the best local solution according to the SDP relaxation scheme
	\bea
	\hatX^{k} = v\parentheses{\hatx^{k}} v\parentheses{\hatx^{k}}\tran,
	\eea
	where $v: \Real{d} \rightarrow \Real{n}$ is a dense or sparse monomial lifting (\eg $v(x) = [x]_\kappa$ is the full set of monomials up to degree $\kappa$ in the case of dense Lasserre's hierarchy at order $\kappa$).
\end{enumerate}

Now, we are given two candidates for the next iteration, namely $ \barX^{k} $ (generated by computing the projection of $X^{k-1} - \sigma_kC$ onto the feasible set $\setsdpp$ inexactly) and $ \hatX^{k} $ (obtained by rounding, local search and lifting described just now), the follow-up question is: which one should we choose to be the next iterate $X^{k}$ such that the entire sequence $\{X^k\}$ is globally convergent for the SDP \eqref{eq:primalSDP}?

{The answer to this question is quite natural --we accept $\hatX^k$ if and only if it attains a strictly lower cost than $\barX^{k}$-- and guarantees that the algorithm visits a sequence of rank-one vertices (local minima via NLP) with descending costs.
With this insight, we now introduce our algorithm~\name~in Algorithm~\ref{alg-iPGMnlp}. \name~is a combination of Algorithm \ref{alg-iPGM} and the acceleration techniques in items \ref{item:rounding}-\ref{item:lifting}, but with a judicious safeguarding policy~\eqref{eq:accept-reject} that ensures (i) the rank-one iterate $\hatX^k$ is feasible (\ie $\hatX^k \in \setsdpp$), and (ii) the rank-one iterate $\hatX^k$ attains a strictly lower cost than $\barX^{k}$ and previously accepted rank-one iterates. Notice that requiring $\hatX^{k}$ to attain a lower cost than previous rank-one iterates can prevent the algorithm from revisiting the same vertex.}

%!TEX root = main.tex

\begin{algorithm}
	\caption{An inexact projected gradient method accelerated \blue{via rounding and lifting by} NLP for solving tight SDP relaxations of POPs (\name).}
	\label{alg-iPGMnlp}
	\begin{algorithmic}
		\item[{\bf Input}:] Initial points $ (X^{0}, y^0, S^0)\in \mathbb{S}_+^n\times \mathbb{R}^{m}\times\mathbb{S}_+^n $, a nondecreasing positive sequences $ \{\sigma_k\}$, a nonnegative sequence $\{\varepsilon_k\}$ such that $\{k\varepsilon_k\}$ is summable, a stopping criterion $\eta:\mathbb{S}^n_+\times \mathbb{R}^m\times \mathbb{S}_+^n\rightarrow \mathbb{R}_+$ with a tolerance $ {\tt Tol} > 0 $. Initialize $\calV = \{ X^0 \}$ if $X^0\in \setsdpp$ and $\calV = \emptyset$ otherwise. Choose a positive integer $r \in [1,n]$, and a positive constant $\epsilon>0$. 
		\item[{\bf Iterate} the following steps for $ k = 1,\dots $:]
		\item[\quad {\bf Step 1 (Projection)}.] Compute $(\overline{X}^k, y^k, S^{k})$ satisfying \eqref{eq-iPGM-proj} and \eqref{eq-iPGM-inexact-conditions} (see Section \ref{sec:apgsubproblem}).
		\item[\quad {\bf Step 2 (Certification)}.] If $ \eta(\barX^{k}, y^{k}, S^{k}) < {\tt Tol} $, {\bf output} $(\barX^{k}, y^{k}, S^{k})$ and {\bf stop}.
		\item[\quad {\bf Step 3 (Acceleration)}.] Compute $ \widehat{X}^{k} $ as in items~\ref{item:rounding}-\ref{item:lifting}:
		\bea
			& \barX^{k} = \sum_{i=1}^n \lambda_i v_i v_i\tran, \quad \bar{x}^{k}_{i} = \rounding(v_i), \quad i = 1,\dots,r, \nonumber \\
			& \hatx_i^{k} = \nlp(\barx^{k}_i), \quad i=1,\dots,r, \quad \hatx^{k} = \argmin_{\hatx^{k}_i, i=1\dots,r} p(\hatx^{k}_i), \nonumber \\
			& \hatX^{k} = v\parentheses{\hatx^{k}} v\parentheses{\hatx^{k}}\tran. \nonumber 
		\eea
		\item[\quad {\bf Step 4 (Policy)}.] Update $X^{k}$ according to 
		\begin{equation}
			\label{eq:accept-reject}
			X^{k} = 
				\begin{cases}
					\hatX^{k} & {\rm if }\; \inprod{C}{\hatX^k} < 
					\blue{\min \left(\inprod{C}{\barX^k}, \min_{X\in \calV} \left\{\inprod{C}{X}\;\middle|\; X\in \calV\right\}\right)} - \epsilon,\; \widehat{X}^{k} \in \setsdpp  \\
					\barX^{k} & {\rm otherwise }
				\end{cases}.
		\end{equation}
		If $X^{k} = \hatX^{k}$, set $\calV \leftarrow \calV \cup \{ \hatX^k \}$.
	\end{algorithmic}
\end{algorithm}

We now state the convergence result for~\name.
\begin{theorem}
	\label{thm-ipgmnlp}
	Let $\{(\barX^k, y^k, S^k)\}$ be any sequence generated by Algorithm \ref{alg-iPGMnlp}, and suppose that there exists a constant $M>0$ such that $\norm{y^k}\leq M$ for all $k\geq 0$. 
	% \red{Assuming Slater's conditions hold for SDP \eqref{eq:primalSDP}-\eqref{eq:dualSDP}} \blue{(LIANG: I think we do not need Slater's conditions here for global convergence, but only need them for solving subproblem. See also the red part in Page 3.)}, 
	If one chooses $\sigma_k \geq O(\sqrt{k})$ for all $k\geq 1$, then, $\left\{\inprod{C}{\barX^k}\right\}$ converges to the optimal value of the SDP problem \eqref{eq:primalSDP}. 
\end{theorem}
\begin{proof}
Let us only consider the case with $X^0 \in \setsdpp$ such that $\calV$ is initialized as $\{ X^0\}$ (the case with $\calV$ initialized as $\emptyset$ can be argued in a similar manner). By the construction of the set $\calV$, it contains a sequence of feasible points of~\eqref{eq:primalSDP}, and the objective value along this sequence is \emph{strictly} decreasing according to the acceptance policy~\eqref{eq:accept-reject}. Hence, we have the following relation:
\begin{equation*}
	\inf\left\{\inprod{C}{X}\;\middle|\; X\in \calV\right\} \leq \inprod{C}{X^0} - (|\calV|-1) \epsilon,
\end{equation*}
If $|\calV|=\infty$, then there exists an infinite sequence of feasible points to problem \eqref{eq:primalSDP} whose objective function value converges to $-\infty$. In this case, problem \eqref{eq:primalSDP} is unbounded, a contradiction (recall that we assumed strong duality throughout this paper in Section \ref{sec:introduction}). Hence, $|\calV|< \infty$. If $|\calV| = 1$, then none of the $\hatX^{k}$ generated by {\bf Step 3} has been accepted, in which case Algorithm \ref{alg-iPGMnlp} reduces to Algorithm \ref{alg-iPGM} iPGM. Therefore, by Theorem~\ref{thm-iPGM}, Algorithm \ref{alg-iPGMnlp} converges as $O(1/k)$ if one chooses $\sigma_k \geq O(\sqrt{k})$ for all $k\geq 1$. If $|\calV| > 1$, then some of the $\hatX^{k}$ generated by {\bf Step 3} have been accepted. Denote the last element of $\calV$ as $\hatX_{\calV}$ with associated dual variables $(y_\calV,S_\calV)$ (from {\bf Step 1}), then Algorithm \ref{alg-iPGMnlp} essentially reduces to \ipgm~with a new initial point at $(\hatX_{\calV}, y_\calV, S_\calV)$. Again, by Theorem~\ref{thm-iPGM}, Algorithm \ref{alg-iPGMnlp} converges as $O(1/k)$. This completes the proof.\qed
\end{proof}

We now make a few remarks about Algorithm \ref{alg-iPGMnlp} (\name).

The first remark is on the local search algorithm. In general, one can use any nonlinear programming method to perform local search on the original~\eqref{eq:pop} starting from an initial guess. A good choice is to use a solver that is based on a primal-dual interior point method, for which the global and local convergence properties are well studied~\cite[Chapter 19]{Jorge2006}. Notice that even if the NLP solver can fail to converge to a feasible point of \eqref{eq:pop}, the safeguarding policy \eqref{eq:accept-reject} ensures that all the rank-one points in $\calV$ are feasible
\blue{for the SDP (P)} (\ie a failed NLP solution will not be accepted). In many POPs arising from practical engineering applications, the constraint set typically defines a smooth manifold (see examples in Section~\ref{sec:numexps}). In such cases, we prefer to use unconstrained optimization algorithms on the manifold as the local search method, for example the Riemannian trust region method that admits favorable global and local convergence properties~\cite[Chapter 7]{Absil09book-manifoldopt}. A general interior point method for NLP is available through~\fmincon~in Matlab (see~\cite[Section 19.9]{Jorge2006} for other available software), while the Riemannian trust region method is available through~\manopt~\cite{manopt}. We emphasize here that the idea of using local search algorithms for accelerating iPGM is heuristic and only supported by numerical experiments in Section \ref{sec:numexps}. {In other words, while we 
guarantee global convergence of {\name} for solving the SDP \eqref{eq:primalSDP} (under mild technical assumptions), the amount of acceleration
gained by nonlinear programming may be 
problem dependent and is mostly observed
empirically.}
% We do not claim any conclusion on that whether such algorithms are able to find a global optimal solution to the original \eqref{eq:pop}. 
In the present paper, we are not able to establish conditions under which the local search algorithms can provide provably better rank-one candidates than the iterates generated by \ipgm. We think this aspect deserves deeper future research. Nevertheless, the nice property of {\name} is that, 
% assuming the SDP relaxation is tight, when the iPGM iterates $\{\barX^k\}$ approach the optimal solution set of the SDP, it is likely that $\barX^k$ can provide good initialization for the local search algorithms to find the optimal rank-one solution. Moreover, 
even if we reject all the candidates provided by local search, global convergence is still guaranteed by taking the safeguarded \ipgm~steps.

Next, we comment on the number of eigenvectors to round. In \name, the hyperparameter $r$ decides how many eigenvectors to round and how many hypotheses to generate at each iteration. Since the NLP solver performs local search very quickly, it is affordable to choose a large $r$. However, we empirically observed that choosing $r$ between $2$ and $5$ is sufficient for finding the global optimal solution. 
% \red{Based on our experience, in most situations, the largest eigenvalue corresponds to the global optimal solution at the terminating state of \name.}

Lastly, we provide a possible extension to \name. The careful reader may have noticed that Algorithm \ref{alg-iPGMnlp} uses the general Algorithm iPGM \emph{without} acceleration. We emphasize here that the accelerated version studied in \cite{Jiang2012} can also be used as the backbone of~\name~with the same global convergence property. The reason we only implement the unaccelerated version is because empirically, with a proper warmstart (described in Section \ref{subsec:initialpoints}), we observe that Algorithm \ref{alg-iPGMnlp} converges in $1$ or $2$ iterations, just like what we have shown in the simple univariate example in Section~\ref{sec:introduction}. Therefore, the acceleration scheme would not have major improvement on the efficiency of Algorithm \ref{alg-iPGMnlp}. 
% Such a fast convergence is enabled by a good initialization scheme that we will soon describe in the next subsection.

%!TEX root = main.tex
\subsection{Initialization}
\label{subsec:initialpoints}
\name~requires an initial guess $(X^0,y^0, S^0) \in \psd{n} \times \Real{m}\times \psd{n}$ as described in Algorithm~\ref{alg-iPGMnlp}. Although one can choose an arbitrary initial point and the algorithm is guaranteed to converge, a good initialization can significantly promote fast convergence.
% In our numerical tests, we found that the choice of initial point affects the overall performance of~\name. 
We first mention that, for many POPs arising from engineering applications, practitioners often have designed powerful heuristics for solving the POPs based on their specific domain knowledge. By heuristics we mean algorithms that can find the \emph{globally optimal} POP solutions with very high probability of success, but \emph{cannot} provide a certificate of global optimality (or suboptimality). Therefore, when such heuristics exist, we {could} use them to generate $x^0$ for the~\eqref{eq:pop} and perform a rank-one lifting of $x^0$ to form $X^0$, \ie $X^0 = v(x^0)v(x^0)\tran$, where $v$ is the lifting monomials. In Section~\ref{sec:wahba}, we demonstrate the effectiveness of one such heuristic in speeding up {\name} on a computer vision application.

Here we describe a general initialization scheme based on a convergent semi-proximal alternating direction method of multipliers (\spadmm)~\cite{Sun2015}. \spadmm~can be used to generate the dual initialization $(y^0,S^0)$ when POP heuristic exists, or to generate both $X^0$ and $(y^0,S^0)$ when no POP heuristics exist. We begin by reformulating problem~\eqref{eq:dualSDP} as follows:
\begin{equation}
	\label{prob:DSDP-ADMM}
	\min_{y \in \Real{m}, S \in \sym{n}} \left\{ \delta_{\mathbb{S}_+^n}^*(-S) - \left\langle b, y \right \rangle\mid \mathcal{A}^{*}y + S = C\right\}.
\end{equation}
Given $ \sigma > 0 $, the augmented Lagrangian associated with~\eqref{prob:DSDP-ADMM} is given as 
\[
	\mathcal{L}_\sigma(y, S; X) = \delta_{\mathbb{S}_+^n}^*(-S) - \left\langle b, y \right \rangle + \left\langle X, \mathcal{A}^{*}y + S - C \right \rangle + \frac{\sigma}{2}\left\lVert \mathcal{A}^{*}y + S - C \right\rVert^{2}, 
\]
for $(X,y,S)\in \sym{n}\times \Real{m}\times \sym{n}$. The \spadmm~method for solving problem \eqref{prob:DSDP-ADMM} is described in Algorithm~\ref{alg:ADMM}, and its convergence is well studied in \cite{Sun2015}. 
An implementation of Algorithm~\ref{alg:ADMM} is included in the software package \sdpnal \cite{Yang15MPC-sdpnalplus}, namely the {\tt admmplus} subroutine. Let the output of {\spadmm} be $(X,y,S)$, we use $(\Pi_{\psd{n}}(X),y,S )$ to be the initial point for Algorithm~\ref{alg-iPGMnlp} (recall that {\name} requires $X^0 \in \psd{n}$, so we project $X$ to be positive semidefinite).

\begin{algorithm}
\caption{A semi-proximal ADMM for solving \eqref{prob:DSDP-ADMM}.}
\label{alg:ADMM}
\begin{algorithmic}
	\item[{\bf Input}:] Initial points $ X^{0} = S^0= 0 \in \mathbb{S}^n $ and $ \gamma \in (0, 2) $.
	\item[{\bf Iterate} the following steps for $ k = 1,\dots $:]
	\item[\quad {\bf Step 1}.] Compute 
		\[
		\begin{split}
			\hat{y}^{k+1} & = \argmin_y \mathcal{L}_\sigma(y, S^{k}; X^{k}) =  \big(\mathcal{A}\mathcal{A}^*\big)^{-1}\left( \frac{1}{\sigma}b - \mathcal{A} \left(\frac{1}{\sigma}X^{k} + S^{k} - C\right)\right).
		\end{split}
		\]
	\item[\quad {\bf Step 2}.] Compute 
		\[
		\begin{split}
			S^{k+1} & = \argmin_S \mathcal{L}_\sigma(\hat{y}^{k+1}, S; X^{k}) \\
			        & = \frac{1}{\sigma}\left( \Pi_{\mathbb{S}_+^n}\left(X^{k} + \sigma(\mathcal{A}^*\hat{y}^{k+1} - C)\right) - (X^{k} + \sigma(\mathcal{A}^*\hat{y}^{k+1} - C)) \right).
		\end{split}
		\]
	\item[\quad {\bf Step 3}.] Compute 
		\[
		\begin{split}
			y^{k+1} & = \argmin_y \mathcal{L}_\sigma(y, S^{k+1}; X^{k}) = \big(\mathcal{A}\mathcal{A}^*\big)^{-1}\left( \frac{1}{\sigma}b - \mathcal{A} \left(\frac{1}{\sigma}X^{k} + S^{k+1} - C\right)\right).
		\end{split}
		\]
	\item[\quad {\bf Step 4}.] Compute 
		\[
			X^{k+1} = X^{k} + \gamma \sigma (S^{k+1} + \mathcal{A}^*y^{k+1} - C).
		\]
	\item[{\bf Until} termination conditions are met.]
	\item[{\bf Output}:] $(X^{k+1}, y^{k+1}, S^{k+1})$.
\end{algorithmic}
\end{algorithm}

\blue{We note that the steplength $\gamma$ in Algorithm~\ref{alg:ADMM} can take values
in the interval $(0,2)$ instead of the usual interval of $(0,(1+\sqrt{5})/2)$.
In Step 1 and 3 of the algorithm, $\hat{y}^{k+1}$ and $y^{k+1}$ are computed
by using the sparse Cholesky factorization of $\cA\cA^*$, which is computed
once at the beginning of the algorithm. One can also compute 
$\hat{y}^{k+1}$ and $y^{k+1}$ inexactly by using a precondtioned conjugate gradient
method without affecting the convergence of the overall algorithm as long as the 
residual norms of the inexact solutions are bounded by a summable error sequence;
we refer the reader to \cite{ChenSunToh2017} for the details.
We further add that for the computation of $S^{k+1}$ in Step 2, one can make use
of the expected low-rank property of $\Pi_{\mathbb{S}_+^n}\left(X^{k} + \sigma(\mathcal{A}^*\hat{y}^{k+1} - C)\right)$. In our implementation we use the subroutine 
{\tt dsyevx} in LAPACK to compute only the positive eigen-pairs
of $X^{k} + \sigma(\mathcal{A}^*\hat{y}^{k+1} - C)$.
}
%\begin{figure}[ht!]
%	\centerline{\fbox{\parbox{\textwidth}{
%	{\bf Algorithm sPADMM: A semi-proximal ADMM for \eqref{prob:DSDP-ADMM}.} \\[5pt]
%	Given $ X^{0} = S^0= 0 \in \mathbb{S}^n $ and $ \gamma \in (0, 2) $. Iterate the following steps.
%	\begin{description}
%		\item[\bf Step 1. ] Compute 
%		\[
%		\begin{split}
%			\hat{y}^{k+1} & = \argmin_y \mathcal{L}_\sigma(y, S^{k}; X^{k}) =  \big(\mathcal{A}\mathcal{A}^*\big)^{-1}\left( \frac{1}{\sigma}b - \mathcal{A} \left(\frac{1}{\sigma}X^{k} + S^{k} - C\right)\right).
%		\end{split}
%		\]
%		\item[\bf Step 2. ] Compute 
%		\[
%		\begin{split}
%			S^{k+1} & = \argmin_S \mathcal{L}_\sigma(\hat{y}^{k+1}, S; X^{k}) \\
%			        & = \frac{1}{\sigma}\left( \Pi_{\mathbb{S}_+^n}\left(X^{k} + \sigma(\mathcal{A}^*\hat{y}^{k+1} - C)\right) - (X^{k} + \sigma(\mathcal{A}^*\hat{y}^{k+1} - C)) \right).
%		\end{split}
%		\]
%		\item[\bf Step 3. ] Compute 
%		\[
%		\begin{split}
%			y^{k+1} & = \argmin_y \mathcal{L}_\sigma(y, S^{k+1}; X^{k}) = \big(\mathcal{A}\mathcal{A}^*\big)^{-1}\left( \frac{1}{\sigma}b - \mathcal{A} \left(\frac{1}{\sigma}X^{k} + S^{k+1} - C\right)\right).
%		\end{split}
%		\]
%		\item[\bf Step 4. ] Compute 
%		\[
%			X^{k+1} = X^{k} + \gamma \sigma (S^{k+1} + \mathcal{A}^*y^{k+1} - C).
%		\]
%	\end{description}
%	}}}
%	\caption{A semi-proximal ADMM. \label{alg:ADMM}}
%\end{figure}

%!TEX root = main.tex
\section{Solving the Projection Subproblem}
\label{sec:apgsubproblem}

Recall that the feasible set of~\eqref{eq:primalSDP} is $\setsdpp = \left\{ X\in \sym{n}\mid \calA(X) = b,\; X\succeq 0 \right\} $. In this section, we aim at computing the projection onto $ \setsdpp $ efficiently since it is required at each iteration in \name (Algorithm \ref{alg-iPGMnlp}). Note that since there are $m$ linear equality constraints defining $\setsdpp$ with $m$ possibly as large as a few millions, the projection algorithm has to be scalable. To this end, we propose a two-phase algorithm. In phase one (Section \ref{subsec:sGSAPG}), we use an {sGS-based} accelerated proximal gradient method to generate a reasonably good initial point for the purpose of warm starting. {Then in} phase two (Section \ref{subsec:lbfgs}), we apply a modified version of the classical limited-memory BFGS (\lbfgs) method to compute a highly accurate solution.

Formally, given a point $Z \in \sym{n}$, the projection problem is described as finding the closest point in $\setsdpp$ with respect to $Z$
\begin{equation}
	\label{prob:projection}
	\min_{X \in \sym{n}}\left\{ \frac{1}{2}\left\lVert X - Z \right\rVert^2 \mid X\in \mathcal{F}_P \right\}.
\end{equation}
Notice that the feasible set $ \setsdpp$ is a spectrahedron defined by the intersection of two convex sets, namely the hyperplane $ \left\{ X\in \sym{n} : \mathcal{A}(X) = b \right\} $ and the PSD cone $ \psd{n} $. Therefore, a natural idea is to apply Dykstra's projection (see~\eg~\cite{Patrick2011}) {for generating} an approximate solution to \eqref{prob:projection} by alternating the projection onto the hyperplane and the projection onto the {PSD} cone, both of which are easy to compute. However, Dykstra's projection is known to have slow convergence rate and it may take too many iterations until a satisfactory approximate projection is found. In this paper, instead of solving \eqref{prob:projection} directly, we consider its dual problem for which we can handle more efficiently.

It is not difficult to write down the Lagrangian dual problem (ignoring the constant term $ -\frac{1}{2}\left\lVert Z \right\rVert^2 $ and converting ``$\max$'' to ``$\min$'') of~\eqref{prob:projection} as:
\begin{equation}
	\label{prob:projection-dual}
	\min_{W, \xi}\left\{ \frac{1}{2}\left\lVert W + \mathcal{A}^*\xi  + Z \right\rVert^2 - \left\langle b, \xi \right \rangle \mid \xi \in \mathbb{R}^{m}, \; W\in \psd{n} \right\}.
\end{equation}
For the rest of this section, we assume that the KKT system for \eqref{prob:projection}-\eqref{prob:projection-dual} given as 
\begin{equation}
	\label{eq:KKT-proj}
	\mathcal{A}(X) = b,\quad \mathcal{A}^*\xi + W = X - Z,\quad X, W\succeq 0,\quad \left\langle X, W \right\rangle = 0,
\end{equation}
admits at least one solution, which is satisfied if problem \eqref{prob:projection} (or \eqref{eq:primalSDP}) satisfies the Slater's condition, \ie there exists an $X^\dagger\in \mathbb{S}^n$ such that
\[
	\cA(X^\dagger) - b = 0,\quad X^\dagger \succ 0.
\]

Before presenting our two-phase algorithm, let us briefly discuss about stopping conditions for solving the pair of projection SDPs~\eqref{prob:projection} and~\eqref{prob:projection-dual}. To this end, let $(W, \xi)\in \psd{n} \times \mathbb{R}^{m} $ be the output of any algorithm that solves the dual problem~\eqref{prob:projection-dual}. From the KKT conditions~\eqref{eq:KKT-proj}, we deem 
\bea\label{eq:XfromSandy}
	X(W, \xi):=\mathcal{A}^*\xi + W + Z\in \mathbb{S}^{n}
\eea
as an approximate solution for the primal problem~\eqref{prob:projection}. Given a tolerance parameter \tol, we accept $ X(W, \xi)$ as an approximate projection point if the relative KKT \blue{residue} is below $\tol$:
\begin{equation}
	\label{eq:relative-KKT-proj}
	\etaproj(W, \xi) := \max\left\{ \frac{\left\lVert \mathcal{A}(X(W, \xi)) - b \right\rVert}{1 + \left\lVert b \right\rVert}, \frac{\left\lVert X(W, \xi) - \Pi_{\mathbb{S}_+^n}(\calAadj \xi + Z) \right\rVert}{1 + \left\lVert  X(W, \xi) \right\rVert + \left\lVert W \right\rVert }\right\} \leq {\tt tol}.
\end{equation}
The first-order optimality conditions for problem \eqref{prob:projection-dual},~\ie
\[
	0\in \begin{pmatrix}
		W + \mathcal{A}^{*}\xi + Z + \partial \delta_{\mathbb{S}_+^n}(W) \\[5pt]
		\mathcal{A}(W + \mathcal{A}^{*}\xi + Z) - b
	\end{pmatrix} = 
	\begin{pmatrix}
		X(W, \xi) +	\partial \delta_{\mathbb{S}_+^n}(W)\\[5pt]
		\mathcal{A}(X(W, \xi)) - b
	\end{pmatrix},
\]
ensures that $\etaproj(W, \xi)$ is small if $ (W, \xi) $ is sufficiently \blue{accurate}, since $ 0\in 	X(W, \xi) + \partial \delta_{\mathbb{S}_+^n}(W) $ implies \blue{that} $ X(W, \xi) - \Pi_{\mathbb{S}_+^n}(X(W, \xi) - W) = X(W, \xi) - \Pi_{\psd{n}}(\calAadj \xi + Z) = 0 $. 
% For the rest of this section, we shall focus on solving the dual problem \eqref{prob:projection-dual} efficiently. 

Now let us present our two-phase algorithm for solving the dual~\eqref{prob:projection-dual} for the remaining part of this section.

%%%%%%%%%%%%%%%%%%%%%%%%%%%%%%%%%%%%%%%%%%%%%%%%%%%%%%%%%%%%%%%%%%%%%%%%%%%%%%%%%%%%%%%%%%%%%%%%%%%%%%%%%%%%%%
%%%%%%%%%%%%%%%%%%%%%%%%%%%%%%%%%%%%%%%%%%%%%%%%%%%%%%%%%%%%%%%%%%%%%%%%%%%%%%%%%%%%%%%%%%%%%%%%%%%%%%%%%%%%%%
%%%% Phase One: sGS-PGM
%%%%%%%%%%%%%%%%%%%%%%%%%%%%%%%%%%%%%%%%%%%%%%%%%%%%%%%%%%%%%%%%%%%%%%%%%%%%%%%%%%%%%%%%%%%%%%%%%%%%%%%%%%%%%%
%%%%%%%%%%%%%%%%%%%%%%%%%%%%%%%%%%%%%%%%%%%%%%%%%%%%%%%%%%%%%%%%%%%%%%%%%%%%%%%%%%%%%%%%%%%%%%%%%%%%%%%%%%%%%%

\subsection{Phase One: An sGS-based Accelerated Proximal Gradient Method}
\label{subsec:sGSAPG}

Problem~\eqref{prob:projection-dual} is in the form of a convex composite minimization problem 
\begin{equation}
	\label{prob:composite}
	\min_{W\in \sym{n}, \xi\in\Real{m}} \left\{ F(W, \xi):= f(W, \xi) + g(W, \xi) \right\},
\end{equation}
where $ g(W, \xi):= \delta_{\mathbb{S}_+^n}(W) $ denotes the indicator function for $ \mathbb{S}_+^n $ which is convex but non-smooth, and $f(W, \xi)$ is the following convex quadratic function
\[
	f(W, \xi):= \frac{1}{2}\left\langle \begin{pmatrix}W \\[5pt] \xi
 \end{pmatrix}, \mathcal{Q}\begin{pmatrix} W \\[5pt] \xi	
 \end{pmatrix} \right \rangle + \left\langle \begin{pmatrix} Z \\[5pt] \mathcal{A}(Z) - b	
 \end{pmatrix}, \begin{pmatrix} W \\[5pt] \xi	
 \end{pmatrix} \right \rangle 
\]
with $\calQ$ written as
\[
	\mathcal{Q}:= 
	\begin{pmatrix}
		\mathcal{I} & \mathcal{A}^* \\[5pt] \mathcal{A} & \mathcal{A} \mathcal{A}^*
	\end{pmatrix}.
\]
It is desirable to apply a proximal-type method for problem \eqref{prob:composite} whose objective is the sum of a smooth and a non-smooth functions. The key idea in a proximal-type method is to approximate the smooth part $f(\cdot)$ by a wisely chosen convex quadratic function, namely $q_k(\cdot)$, such that (i) $q_k + g$ is a good approximation to $f + g$ near the current iterate, and (ii) minimizing $q_k + g$ can be solved efficiently. In particular, one needs to choose an appropriate positive definite mapping $\calH_k$ ($k \geq 1$) for approximating the function $f(\cdot)$ by $q_k(\cdot)$ that is defined as follows:
\bea\label{eq:pgmsubproblemprojection}
	q_k(W, \xi) := f(\widetilde{W}^k, \tilde{\xi}^k) + \left\langle \nabla f(\widetilde{W}^k, \tilde{\xi}^k), \begin{pmatrix}
		W - \widetilde{W}^k \\[5pt] \xi - \tilde{\xi}^k
	\end{pmatrix} \right\rangle + \frac{1}{2} \left\lVert \begin{pmatrix}
		W - \widetilde{W}^k \\[5pt] \xi - \tilde{\xi}^k
	\end{pmatrix}\right\rVert_{ \mathcal{H}_k}^2
\eea
for a given $(\widetilde{W}^k, \tilde{\xi}^k) \in \psd{n}\times \Real{m}$. There are many possible choices for $ \mathcal{H}_k$. For example, $ \mathcal{H}_k $ can be chosen as $\mathcal{Q}$ since the function $f(\cdot)$ is itself a quadratic function. However, since the variables $W$ and $\xi$ are \blue{coupled}, minimizing the approximate function $q_k(W,\xi) + g(W, \xi)$, for $\calH_k = \calQ$, is as difficult as solving the original projection problem~\eqref{prob:composite}. Another possibility is to choose $\calH_k = L_k \mathrm{diag}(\mathcal{I}, \calI_m) -  \calQ$ where $L_k$ is a given positive scalar. However, the value of $L_k$ could be difficult to choose. Indeed, choosing $L_k$ too small or too large can both harm the efficiency and convergence. Even though one could apply certain back-tracking techniques (see \eg \cite[Section 4]{Beck2009}) for estimating $L_k$ at each iteration, the computational cost may be expensive. To decouple the variables $W$ and $\xi$ and to achieve fast convergence, we will choose the operator $\calH_k$ by applying the symmetric block Gauss-Seidel (sGS) decomposition technique in~\cite{Li2019}.

To present the sGS decomposition method, we first need to introduce some useful notation. Let $ \mathcal{Q} $ have the following decomposition
\begin{equation}
	\label{eq:Q-decomp}
	\mathcal{Q} = \mathcal{U} + \mathcal{D} + \mathcal{U}^{*},
\end{equation}
where 
\[
	\mathcal{U} = 
	\begin{pmatrix}
		0 & \mathcal{A}^* \\[5pt] 0 & 0
	\end{pmatrix},\quad \mathcal{D} =
	\begin{pmatrix}
		\mathcal{I} & 0 \\[5pt] 0 & \mathcal{A}\mathcal{A}^*
	\end{pmatrix},
\]
and $\calD \succ 0$ because $\calA$ is onto.
Define the linear operator  $ \mathcal{T}_{\mathcal{Q}} := \mathcal{U}\mathcal{D}^{-1}\mathcal{U}^{*} $. Clearly, $ \mathcal{T}_{\mathcal{Q}} $ is positive semidefinite. In sGS decomposition, we choose 
\bea\label{eq:sgsHk}
\calH_k:= \calQ + \calT_{\calQ}, \ \  \forall k.
\eea 
Let us first verify that $\calH_k$ is positive definite. To show this, we write
\bea
\calQ + \calT_{\calQ} = (\calD + \calU + \calUadj )+ \calU \calD\inv \calUadj = (\calD+\calU)\calD\inv(\calD + \calUadj),
\eea
and easily see that $\calH_k \succ 0$ due to $\calD \succ 0$ and $\calD + \calU$ being nonsingular. Inserting this choice of $\calH_k=\calQ + \calT_\calQ$ into the expression $q_k$ in~\eqref{eq:pgmsubproblemprojection}, we have
\bea
q_k(W, \xi) = f(W, \xi) + \frac{1}{2}\left\langle \begin{pmatrix}
		W- \widetilde{W}^k \\[5pt] \xi - \tilde{\xi}^k
	\end{pmatrix}, \mathcal{T}_{\mathcal{Q}}\begin{pmatrix}
		W- \widetilde{W}^k \\[5pt] \xi - \tilde{\xi}^k
	\end{pmatrix} \right \rangle,
\eea 
and the subproblem at the $k$-th iteration to be solved is 
\begin{equation}
	\label{prob:sgs-decomp}
	\min_{W \in \sym{n}, \xi \in \Real{m}}\left\{  q_k(W, \xi) + g(W,\xi) \right\}.
\end{equation}
The nice property of choosing $\calH_k$ as in~\eqref{eq:sgsHk} and arriving at the subproblem~\eqref{prob:sgs-decomp} is that an optimal solution for \eqref{prob:sgs-decomp} can be computed efficiently by solving three simpler problems, as stated in the following theorem.
\begin{theorem}
	\label{thm:sGS}
	Assume $\mathcal{A}\mathcal{A}^*$ is nonsingular, then the optimal solution $ (W^+, \xi^+)\in \psd{n} \times \Real{m} $ for \eqref{prob:sgs-decomp} can be computed exactly via the following steps:
	\begin{equation}
		\label{eq:thm-SGS}
		\left\{
			\begin{aligned}
				% \xi^\dagger:=&\argmin_{\xi \in \mathbb{R}^{m}}\left\{ \delta_{\mathbb{S}_+^n}(\widetilde{W}^k)+ f(\widetilde{W}^k, \xi)\right\}, \\[5pt]
				\xi^\dagger:=&\argmin_{\xi \in \mathbb{R}^{m}}\left\{ f(\widetilde{W}^k, \xi)\right\}, \\[5pt]
				W^+:=& \argmin_{W \in \mathbb{S}^{n}}\left\{ \delta_{\mathbb{S}_+^n}(W) + f(W, \xi^\dagger) \right\}, \\[5pt]
				y^+:=&\argmin_{\xi\in \mathbb{R}^{m}}\left\{ f(W^+, \xi)\right\}.
				% y^+:=&\argmin_{\xi\in \mathbb{R}^{m}}\left\{ \delta_{\mathbb{S}_+^n}(W^+)+ f(W^+, \xi)\right\}.
			\end{aligned}
		\right.
	\end{equation}
\end{theorem}
The proof of Theorem~\ref{thm:sGS} is given in~\cite[Theorem 1]{Li2019} and is omitted here. The reason that procedure~\eqref{eq:thm-SGS} is simple is because (i) with $W$ fixed as $\widetilde{W}^k$ or $W^{+}$ (first and third line in~\eqref{eq:thm-SGS}), optimizing $\xi$ boils down to solving a linear system and can be solved exactly due to the sparsity of $\calA$; (ii) with $\xi$ fixed as $\xi^\dagger$ (second line in~\eqref{eq:thm-SGS}), optimizing $W$ resorts to performing a projection onto $\psd{n}$ and can be done in closed form via computing eigenvalue decompositions.

With these preparations, we now formally present the sGS-based accelerated proximal gradient method (\sgspgm) for solving~\eqref{prob:projection-dual} in Algorithm~\ref{alg:sgs-apg}, where we also give the closed-form solution of the procedure~\eqref{eq:thm-SGS} in~\eqref{eq:sgs-apg-1}. Note that in Algorithm \ref{alg:sgs-apg}, the computation of $\xi^k$ in {\bf Step 1} is only needed for checking terminations.

%!TEX root = main.tex
%\begin{figure}[ht!]
%	\centerline{\fbox{\parbox{\textwidth}{
%	{\bf Algorithm sGS-PGM: An sGS-based proximal gradient method for \eqref{prob:projection-dual}.} \\[5pt]
%	Given tolerance $ {\tt tol} > 0 $, $ \widetilde{W}^1 = W^{0}\in \mathbb{S}_+^n $ and $ t_1 = 1 $. For $ k = 1, 2, \dots $, perform
%	\begin{description}
%	\item[{\bf Step 1.}] (sGS update). Compute 
%		\begin{equation}
%			\label{eq:sgs-apg-1}
%			\begin{aligned}
%				\tilde{\xi}^{k} = (\calA\calAadj)\inv \left( b- \mathcal{A}(Z) - \mathcal{A}(\widetilde{W}^{k}) \right), \\[5pt]
%				W^{k} = \Pi_{\mathbb{S}_+^n}\parentheses{-\mathcal{A}^*\tilde{\xi}^k - Z}, \\[5pt]
%				\xi^{k} = (\calA\calAadj)\inv\left( b- \mathcal{A}(Z) - \mathcal{A}(W^{k}) \right).
%			\end{aligned}
%		\end{equation}
%	\item[{\bf Step 2.}] (Check termination). Output $ (W^k, \xi^k) $ if $ \etaproj(W^k, \xi^k) \leq \tol $.
%	\item[{\bf Step 3.}] Compute $ t_{k+1} = \frac{1+\sqrt{1+4t_k^2}}{2} $, and $ \widetilde{W}^{k+1} = W^{k} + \frac{t_k - 1}{t_{k+1}}(W^{k} - W^{k-1}) $.
%	\end{description}
%	}}}
%	\caption{An sGS-based proximal gradient method for \eqref{prob:projection-dual}. \label{alg:sgs-apg}}
%\end{figure}

\begin{algorithm}
\caption{An sGS-based accelerated proximal gradient method for \eqref{prob:projection-dual}.}
\label{alg:sgs-apg}
\begin{algorithmic}
	\item[{\bf Input}:] Initial points $ \widetilde{W}^1 = W^{0}\in \mathbb{S}_+^n $, termination tolerance $ {\tt tol} > 0 $, and $ t_1 = 1 $.
	\item[{\bf Iterate} the following steps for $ k = 1,\dots $:]
	\item[\quad {\bf Step 1 (sGS update).}] Compute 
		\begin{equation}
			\label{eq:sgs-apg-1}
			\begin{aligned}
				\tilde{\xi}^{k} = (\calA\calAadj)\inv \left( b- \mathcal{A}(Z) - \mathcal{A}(\widetilde{W}^{k}) \right), \\[5pt]
				W^{k} = \Pi_{\mathbb{S}_+^n}\parentheses{-\mathcal{A}^*\tilde{\xi}^k - Z}
				= \blue{\Pi_{\mathbb{S}_+^n}\parentheses{\mathcal{A}^*\tilde{\xi}^k + Z} 
				- \parentheses{\mathcal{A}^*\tilde{\xi}^k + Z}}, 
				\\[5pt]
				\xi^{k} = (\calA\calAadj)\inv\left( b- \mathcal{A}(Z) - \mathcal{A}(W^{k}) \right).
			\end{aligned}
		\end{equation}
		\item[\quad {\bf Step 2.}] Compute $ t_{k+1} = \frac{1+\sqrt{1+4t_k^2}}{2} $, and $ \widetilde{W}^{k+1} = W^{k} + \frac{t_k - 1}{t_{k+1}}(W^{k} - W^{k-1}) $.
	\item[{\bf Until}:] $ \etaproj(W^k, \xi^k) \leq \tol $.
	\item [{\bf Output}:]$ (W^k, \xi^k) $.
\end{algorithmic}
\end{algorithm}

The convergence of {\sgspgm} is stated as follows.
\begin{theorem}\label{converg-thm-sgs-apg}
	Suppose that $(\Wstar, \xistar)\in \mathbb{S}_+^n\times \mathbb{R}^{m} $ is an optimal solution of problem \eqref{prob:projection-dual} and that $\mathcal{A}\mathcal{A}^*$ is nonsingular. Let $ \{(W^{k}, \xi^{k})\} $ be the sequence generated by Algorithm \ref{alg:sgs-apg}. Then there exists a constant $ c \geq 0 $ such that 
	\[
		0\leq f(W^{k}, \xi^k) - f(\Wstar, \xistar) \leq \frac{c}{(k+1)^2}.
	\]
\end{theorem}

The proof of Theorem \ref{converg-thm-sgs-apg} can be taken directly from \cite[Proposition 2]{Li2019}. Note that the proposed method in this subsection is in fact a direct application of the one proposed in \cite[Section 4]{Li2019} which is designed for a more general convex composite quadratic programming model with multiple blocks. Moreover, \cite[Section 4]{Li2019} also investigates the inexact computation for the corresponding subproblem. Thus, when $\mathcal{A}\mathcal{A}^*$ cannot be factorized efficiently, one may apply an iterative solver
% , such as a \red{preconditioned minimum residual method (MINRES) \cite{Paige-Saunders1975} or a preconditioned symmetric quasi-minimal residual method (PSQMR) \cite{Freund1994}}, 
for solving the linear systems in \eqref{eq:sgs-apg-1} of Algorithm \ref{alg:sgs-apg}. We omit the details here since for the applications studied in this paper, factorizing $\mathcal{A}\mathcal{A}^*$ can always be done efficiently due to the sparsity of $\calA$.
\blue{We further remark that for the computation of $W^k$ in \eqref{eq:sgs-apg-1}, 
we can make use of the expected low-rank property of the projection $\Pi_{\mathbb{S}_+^n}
(\mathcal{A}^*\tilde{\xi}^k + Z)$ corresponding to the primal variable. 
Again, we use the subroutine {\tt dsyevx} in LAPACK to compute
only the positive eigen-pairs of ${\mathcal{A}^*\tilde{\xi}^k + Z}$.
}

%%%%%%%%%%%%%%%%%%%%%%%%%%%%%%%%%%%%%%%%%%%%%%%%%%%%%%%%%%%%%%%%%%%%%%%%%%%%%%%%%%%%%%%%%%%%%%%%%%%%%%%%%%%%%%
%%%%%%%%%%%%%%%%%%%%%%%%%%%%%%%%%%%%%%%%%%%%%%%%%%%%%%%%%%%%%%%%%%%%%%%%%%%%%%%%%%%%%%%%%%%%%%%%%%%%%%%%%%%%%%
%%%% Phase Two: LBFGS
%%%%%%%%%%%%%%%%%%%%%%%%%%%%%%%%%%%%%%%%%%%%%%%%%%%%%%%%%%%%%%%%%%%%%%%%%%%%%%%%%%%%%%%%%%%%%%%%%%%%%%%%%%%%%%
%%%%%%%%%%%%%%%%%%%%%%%%%%%%%%%%%%%%%%%%%%%%%%%%%%%%%%%%%%%%%%%%%%%%%%%%%%%%%%%%%%%%%%%%%%%%%%%%%%%%%%%%%%%%%%

\subsection{Phase Two: A Modified Limited-memory BFGS Method} 
\label{subsec:lbfgs}
% In the last subsection, we proposed an sGS based PGM for computing the dual projection problem onto the feasible set $ \mathcal{F}_P $. We also presented the convergence result for the sGS based PGM by applying the sGS decomposition Theorem \ref{thm:sGS} and Theorem \ref{thm-no-A} directly. However, 
Even though the {\sgspgm} method, as presented in Algorithm~\ref{alg:sgs-apg}, converges as $ O(1/k^2)$ in terms of objective value, the constant $ c$ depends on the distance between the initial point and the optimal solution, \blue{which can be quite large if the initial point
is not well chosen}. Therefore, it may require many iterations to reach a satisfactory solution when the initial point is far way from the optimal solution set. 
\blue{Moreover, the {\sgspgm} method may not have a fast
local linear convergent property since the problem \eqref{prob:projection-dual} is not strongly convex.}
To speed up the computation of the dual projection problem \eqref{prob:projection-dual}, we now propose a modified version of the classical limited-memory BFGS (\lbfgs) method that is warm-started by the solution from Algorithm \ref{alg:sgs-apg}. 

To proceed, we observe that, fixing the unconstrained $\xi$, the dual projection problem~\eqref{prob:projection-dual} becomes finding the closest $W \in \psd{n}$ to the matrix $-(\calAadj \xi + Z)$ and admits a closed-form expression
\bea\label{eq:Sofy}
W 	= \Pi_{\psd{n}} \parentheses{-\calAadj \xi - Z}.
\eea 
As a result, problem~\eqref{prob:projection-dual} can be further simplified, after inserting~\eqref{eq:Sofy}, as 
\begin{equation}
	\label{prob:projetion-dual-y}
	\min_{\xi \in \Real{m}} \quad \phi(\xi):= \frac{1}{2}\left\lVert \Pi_{\mathbb{S}_+^n}(\mathcal{A}^*\xi + Z) \right\rVert^2 - \left\langle b, \xi \right \rangle,
\end{equation}
with the gradient of $ \phi(\xi) $ given as
\begin{equation}
	\label{eq:grad-phi}
	\nabla \phi(\xi) = \mathcal{A}\Pi_{\mathbb{S}_+^n}(\mathcal{A}^*\xi + Z) - b.
\end{equation}
We have used the following equality in getting~\eqref{prob:projetion-dual-y}
\bea\nonumber
\Pi_{\psd{n}} \parentheses{-\calAadj \xi - Z} + \calAadj \xi + Z = \Pi_{\psd{n}}(\calAadj \xi + Z)
\eea
due to {the Moreau identity}~\cite[Theorem 2.2]{Malick06sva-projsdp}.
Thus, if $ \xistar $ is an optimal solution for problem~\eqref{prob:projetion-dual-y}, then we can recover $\Wstar$ from~\eqref{eq:Sofy}. Formulating the dual projection problem as~\eqref{prob:projetion-dual-y} has appeared multiple times; {see, for instances} \cite{Zhao2010,Malick09siopt-sdpregularization}.
\blue{We note that the function $\phi(\cdot)$ is smooth but not twice continuously 
differentiable, and it is convex but not strongly convex.}

Now that~\eqref{prob:projetion-dual-y} is an unconstrained convex minimization problem in $ \xi \in \mathbb{R}^{m} $ with gradient given in \eqref{eq:grad-phi}, 
many efficient algorithms are available \blue{for solving the problem}, such as (accelerated) gradient \blue{descent} methods~\cite{Nesterov18book-convexopt}, nonlinear conjugate gradient methods~\cite{Dai99siopt-ncg}, quasi-Newton methods~\cite{Jorge2006} and the semi-smooth Newton method~\cite{Zhao2010}. For this paper, we decide to propose a modified limited-memory BFGS (\lbfgs) method because \lbfgs~is easy to implement, can handle very large unconstrained optimization problems, and is typically the ``the algorithm of choice'' for large-scale problems.
 Empirically, we observed that our proposed \lbfgs~warm-started by \sgspgm~is efficient and robust for various class of applications (shown in Section~\ref{sec:numexps}). This observation also supports the discussions made in~\cite[Chapter 7]{Jorge2006}. To the best of our knowledge, this is the first work that demonstrates the effectiveness of \lbfgs, or in general quasi-Newton methods, on solving large and challenging SDPs.\footnote{From our extensive numerical trials, we found that accelerated gradient descent, nonlinear conjugate gradient, and semismooth Newton with CG iterative solver fail to give satisfactory performance, especially when $m$ is very large.}

The template for the modified \lbfgs~method is presented in Algorithm \ref{alg:lbfgs}.

\begin{algorithm}
\caption{A modified Limited-memory BFGS method for \eqref{prob:projetion-dual-y}.}
\label{alg:lbfgs}
\begin{algorithmic}
	\item[{\bf Input}:] Starting point $ \xi^0 $, integer $ {\tt mem} > 0 $, a positive constant $K$, termination tolerance $ {\tt tol} > 0 $, $ \mu\in (0, 1/2) $, $ \mu\leq \mu' < 1 $, $ \rho \in (0,1) $ and $ \tau_1,\tau_2\in (0,\infty) $.
	\item[{\bf Iterate} the following steps for $ k = 1,\dots $:]
	\item[\quad {\bf Step 1 (Search direction).}] Choose $ Q_k^0\succ 0 $, $ \beta_k:= \tau_1 \left\lVert \nabla \phi(\xi^{k}) \right\rVert^{\tau_2} $ and compute $ d^k = -\beta_k \nabla \phi(\xi^k) - g^k $ where $ g^k := Q_k\nabla \phi(\xi^k) $ with $Q_k\succeq 0$, is obtained via the two-loop recursion as in \cite[Algorithm 7.4]{Jorge2006}. If $\norm{d^k}\geq K$, then choose $d^k = -\beta_k\nabla \phi(\xi^k)$ (i.e., set $Q_k = 0$).
	\item[\quad {\bf Step 2 (Line search).}] Set $ \alpha_k = \rho^{m_k} $, where $ m_k $ is the smallest nonnegative integer $ m $ such that
		\[
			\phi(\xi^{k}+\rho^{k}d^k) \leq \phi(\xi^k) + \mu \rho^m \left\langle \nabla \phi(\xi^k), d^k \right \rangle. 
		\]
		Compute $ \xi^{k+1} = \xi^{k} + \alpha_k d^{k} $. 
	\item[\quad {\bf Step 3 (Update memory).}] Compute and save $ u^k = \xi^{k+1} - \xi^{k} $ and $ w^{k} = \nabla \phi(\xi^{k+1}) - \nabla \phi(\xi^{k}) $. If $ k > {\tt mem} $, discard the vector $ \{u^{k-{\tt mem}}, w^{k-{\tt mem}}\} $ from storage.
	\item[{\bf Until}:] $ \etaproj(W^{k+1}, \xi^{k+1}) \leq \tol $ with $ W^{k+1} = \Pi_{\mathbb{S}_+^n}(-\mathcal{A}^{*}\xi^{k+1} - Z) $. 
	\item [{\bf Output}:]$ (W^{k+1}, \xi^{k+1}) $.
\end{algorithmic}
	
\end{algorithm}

In Algorithm \ref{alg:lbfgs}, we always choose $ Q_k^{0} = I_m$ for all $ k\geq 0 $ which implies that $ Q_k\succeq 0 $ (see~\eg~\cite[eq. (7.19)]{Jorge2006}). Therefore, the matrix $ Q_k + \beta_k I_m \succ 0 $ when $ \nabla \phi(\xi^k) $ is not zero, and $ d^{k} $ computed in {\bf Step 1} in Algorithm \ref{alg:lbfgs} can be shown to be a descent direction {(\ie the algorithm is well-defined)}. We state the convergence property of Algorithm \ref{alg:lbfgs} in the following theorem, {whose proof is presented in Appendix \ref{sec:proof:theorem-lbfgs}.
\begin{theorem}
	\label{thm:convergence-lbfgs}
	Suppose that $\mathcal{A}\mathcal{A}^*$ is nonsingular and the Slater condition holds for \eqref{eq:primalSDP}. Then Algorithm \ref{alg:lbfgs} is well defined. Let the sequence $ \{\xi^{k}\} $ be generated by Algorithm \ref{alg:lbfgs} such that $ \nabla \phi(\xi^k) \neq 0 $, $ \forall\,k\geq 0 $. Then, the sequence $ \{\xi^{k}\} $ is bounded and any accumulation point $ \hat{\xi} $ of this sequence is an optimal solution of \eqref{prob:projetion-dual-y}. 
\end{theorem}
%The proof of Theorem \ref{thm:convergence-lbfgs} is presented in Appendix \ref{sec:proof:theorem-lbfgs}.
%%%%%%%%%%%%%%%%%%%%%%%%%%%%%%%%%%%%%%%%%%%%%%%%%%%%%%%%%%%%%%%%%%%%%%%%%%%%%%%%%%%%%%%%%%%%%%%%%%%%%%%%%%%%%%
%%%%%%%%%%%%%%%%%%%%%%%%%%%%%%%%%%%%%%%%%%%%%%%%%%%%%%%%%%%%%%%%%%%%%%%%%%%%%%%%%%%%%%%%%%%%%%%%%%%%%%%%%%%%%%
%%%% Back to original problem
%%%%%%%%%%%%%%%%%%%%%%%%%%%%%%%%%%%%%%%%%%%%%%%%%%%%%%%%%%%%%%%%%%%%%%%%%%%%%%%%%%%%%%%%%%%%%%%%%%%%%%%%%%%%%%
%%%%%%%%%%%%%%%%%%%%%%%%%%%%%%%%%%%%%%%%%%%%%%%%%%%%%%%%%%%%%%%%%%%%%%%%%%%%%%%%%%%%%%%%%%%%%%%%%%%%%%%%%%%%%%

\subsection{Connection between the Projection SDP and the Original SDP}
\label{sec:connectprojectsdpandoriginalsdp}
Recall from Algorithm \ref{alg-iPGMnlp} that, $Z = X^{k-1} - \sigma_k C$ at iteration $k$ for $k\geq 1$. Also recall from~\eqref{eq:XfromSandy} that $X(W, \xi) = \calAadj \xi + Z + W $ with $W = \Pi_{\mathbb{S}_+^n}(-\calA^*\xi - Z) \in \mathbb{S}_+^n$. Combining these equations, we have $X(W, \xi) = \Pi_{\mathbb{S}_+^n}(\calA^*\xi + Z) \in \mathbb{S}_+^n$ and hence, $\inprod{X(W, \xi)}{W} = 0$. Furthermore, $X(W, \xi) = \calAadj \xi + W + Z$ implies
\bea
\frac{1}{\sigma_k}\left(X(W, \xi) - X^{k-1}\right) = \calAadj\left(\frac{1}{\sigma_k}\xi \right) + \frac{1}{\sigma_k}W - C,
\eea
which means that when $X(W, \xi)$ is close to $X^{k-1}$ (which is the case when $X^{k-1}$ is approximately optimal), then $\left(\frac{1}{\sigma_k}W, \frac{1}{\sigma_k}\xi\right)$ is an approximate solution for the original dual problem~\eqref{eq:dualSDP}. Therefore, from the solution $\xi$ of {\lbfgs}, we output 
\bea
W = & \Pi_{\psd{n}}(-\calAadj \xi - Z) \\
X(W,\xi) = & \calAadj \xi + Z + W  = \Pi_{\psd{n}}(\calAadj \xi + Z) \\
(\barX^{k}, y^{k}, S^{k}) = & \left(X(W,\xi), \frac{1}{\sigma_k}\xi, \frac{1}{\sigma_k}W \right)
\eea
for {\bf Step 1} of the $k$-th iteration of Algorithm~\ref{alg-iPGMnlp}. In particular, if 
\[
\norm{\nabla \phi(\xi)} = \norm{\calA(\barX^k) - b} \leq \varepsilon_k
\]
for given $\varepsilon_k\geq 0$, then $(\barX^{k}, y^{k}, S^{k})$ satisfies the inexact conditions \eqref{eq-iPGM-inexact-conditions}.

%!TEX root = main.tex
\section{Applications and Numerical Experiments}
\label{sec:numexps}
In this section, we conduct numerical experiments by using~\name~to solve SDP relaxations for {several important classes of} POPs. In Section~\ref{sec:bqp} and~\ref{sec:q4s}, we solve \emph{dense} second-order moment relaxation of random \emph{binary quadratic programming} problems, and random \emph{quartic optimization problems over the unit square}, both with increasing number of variables $d$. We demonstrate, for the first time, that the relaxation is always tight for the instances we have generated up to $d=60$. In Section~\ref{sec:wahba} and~\ref{sec:STLS}, we solve \emph{sparse} second-order moment relaxations coming from two engineering applications, namely an \emph{outlier-robust Wahba problem} that underpins many computer vision applications, and a \emph{nearest structured rank deficient matrix problem} that finds extensive applications in control, statistics, computer algebra, among others. We demonstrate that, with a sparse relaxation scheme, we can globally solve POP problems with $d$ up to $1000$, far beyond the reach of SDP solvers used in the corresponding engineering literature. When solving the outlier-robust Wahba problem in Section~\ref{sec:wahba}, we additionally show that (i) leveraging domain-specific POP heuristics for primal initialization can further speed up {\name} by 2-3 times, and (ii) {\name} can \emph{certifiably optimally} solve two real applications of the outlier-robust Wahba problem, namely image stitching and scan matching. Our experiments show that {\name} is the \emph{only} solver that can consistently solve rank-one tight semidefinite relaxations to high accuracy (\eg KKT residuals below $1\ee{-9}$), in the presence of millions of equality constraints. 

% \red{Delete?In Section~\ref{sec:primalinit}, we revisit the two POPs in Section~\ref{sec:wahba} and~\ref{sec:STLS}, but this time leveraging domain-specific POP heuristics to generate better initialization, as mentioned in Section~\ref{subsec:initialpoints}. We show that~\name~achieves \red{2-10} times speedup when combined together with domain-specific heuristics. }
% Lastly, in Section~\ref{sec:exp:maxcut}, we solve two types of classical first-order MAXCUT relaxations. The first type is those appearing in standard benchmarks such as GSET and admit high-rank solutions. The second type is those that have recently been discovered by Waldspurger and Waters~\cite{Waldspurger20SIOPT-BMRank}, which admit \emph{spurious} rank-one optimal solutions that require high rank for the Burer-Monteiro approach to converge. We demonstrate that (i) \name~compares favorably with existing solvers on solving MAXCUT problems with high-rank solutions; (ii)~\name~is not subject to the high-rank difficulty exhibited by B-M factorization and finds rank-one optimal MAXCUT solutions easily.

Before presenting each application, let us describe the details about implementation and experimental setup.

{\bf Implementation}. We implement the \name Algorithm~\ref{alg-iPGMnlp} in MATLAB R2020a, with core subroutines, such as projection onto the PSD cone, implemented in C for efficiency. Our implementation is available at
\begin{center}
\url{https://github.com/MIT-SPARK/STRIDE}
\end{center}
and also supports SDP problems with multiple PSD blocks (hence, it can also handle relaxations of \eqref{eq:pop} problems with inequality constraints). 

Note that different from generic SDP solvers, we also implement the rounding, local search and lifting procedures required in Step 2 of Algorithm~\ref{alg-iPGMnlp}, for each POP. Since these procedures are problem dependent, we will describe them in the corresponding subsections. 

{\bf Stopping Conditions}.
To measure the feasibility and optimality at a given approximate solution $ (X, y, S)\in \psd{n}\times \mathbb{R}^{m}\times \psd{n} $, we define the following standard relative KKT 
\blue{residues}:
\begin{equation}
	\label{eq:relative-KKT}
	\eta_p = \frac{\left\lVert \mathcal{A}(X) - b \right\rVert}{1 + \left\lVert b \right\rVert},\ \eta_d = \frac{\left\lVert \mathcal{A}^*y + S - C \right\rVert}{1 + \left\lVert C \right\rVert},\ \eta_g = \frac{|\left\langle C, X \right \rangle - \left\langle b, y \right \rangle|}{1 + |\left\langle C, X \right \rangle| + |\left\langle b, y \right \rangle|}.
\end{equation}
For a given tolerance $ \tol > 0 $, we terminate \name~when $\max\{\eta_p,\eta_d,\eta_g\}\leq \tol$, and we choose $\tol = 1\ee{-8}$ for all our experiments. Because our goal is to obtain a solution of the original \eqref{eq:pop} with an optimality or suboptimality certificate, we also compute a relative suboptimality gap from the SDP solution $(X,y,S)$ as
\bea \label{eq:relsubopt}
\subopt = \frac{\abs{ p(\hatx) - (\inprod{b}{y} + M_b \lambda_\min (C- \calAadj y)) } }{1 + \abs{p(\hatx)} + \abs{\inprod{b}{y} + M_b \lambda_\min (C- \calAadj y)} },
\eea
where $\hatx \in \setpop$ is a feasible approximate solution to the \eqref{eq:pop} that is rounded from the leading eigenvector of $X$,\footnote{For a generic POP, even finding a feasible $x \in \setpop$ is NP-hard. However, in all our numerical examples, rounding a feasible point from the approximate SDP solution is easy.} $\lambda_{\min}$ denotes the minimum eigenvalue, and $M_b \geq \trace{X}$ is a bound on the trace of $X$ when $X$ is generated by a rank-one lifting. We easily have $p(\hatx) \geq \fpopstar \geq \inprod{b}{y} + M_b \lambda_\min (C- \calAadj y)$, and hence $\subopt = 0$ (\ie the upper bound and the lower bound coincide) certifies that $p(\hatx) = \fpopstar$ and $\hatx$ is a global minimizer of the nonconvex \eqref{eq:pop}.

{\bf Baseline Solvers.}
We compare~\name~with a diverse set of existing SDP solvers. We choose \sdpt~\cite{Toh99OMS-sdpt3} and~\mosek~\cite{mosek} as representative interior point methods; \cdcs~\cite{Zheng20MP-CDCS} and~\sketchy~\cite{Yurtsever21SIMDS-scalableSDP} as representative first-order methods; and~\sdpnal~\cite{Yang15MPC-sdpnalplus} as a representative method that combines first-order and second-order Newton-type methods. For~\sdpt~and~\mosek, we use default parameters. For~\cdcs, we use the \sos~solver with maximum $20,000$ iterations instead of the default homogeneous self-dual embedding solver because we found~\sos~to typically perform better. 
% We reduce the maximum number of iterations where we find it to take an unreasonable amount of time. 
For \sketchy, we use the default parameters with sketching size $10$. We set the maximum runtime to be $10,000$ seconds and maximum number of iterations to be $20,000$. For~\sdpnal, we use $1\ee{-8}$ as the tolerance, and we run it for maximum $20,000$ iterations and $10,000$ seconds. For very large problems (\eg $m$ above one million), we increase the maximum runtime of {\sketchy} and {\sdpnal}, which will be described in relevant subsections. 

{\bf Hardware.} All experiments are performed on a Linux PC with 12-core Intel i9-7920X CPU@2.90GHz and 128GB RAM.

% \red{Present results on average to save space.}

\subsection{Binary quadratic programming}
\label{sec:bqp}
Consider minimizing a \emph{quadratic} polynomial over the $d$-dimensional binary cube
\begin{equation}\label{eq:bqp}
\min_{x_i \in \{+1,-1\}, i=1,\dots,d} \inprod{c}{[x]_2} \tag{BQP}
\end{equation}
where $[x]_2: \Real{d} \rightarrow \Real{\bard_2}$ is the vector of monomials in $x$ of degree up to 2, and $c \in \Real{\bard_2}$ contains the coefficients of all monomials. Problem~\eqref{eq:bqp} is a classical NP-hard combinatorial problem with examples such as the maximum cut (MAXCUT) problem~\cite{Goemans95JACM-maxcut}, the $0$-$1$ knapsack problem~\cite{Helmberg00JCO-SDPKnapsack}, the number partitioning problem~\cite{Mertens06CCSP-npp,Gamarnik21arxiv-npp}, and the linear quadratic regulator control problem with binary inputs~\cite{Wu18TAC-LQRswitchedsystems}. It is well known that the standard Shor's semidefinite relaxation for~\eqref{eq:bqp} is typically not tight, \eg in MAXCUT problems, and hence a globally optimal solution cannot be obtained with an optimality certificate (albeit a lower bound can be obtained). 

We consider the \emph{second-order} dense moment relaxation for~\eqref{eq:bqp}, which creates a \emph{positive semidefinite} moment matrix $X := [x]_2 [x]_2\tran$, using which the cost function in~\eqref{eq:bqp} can be written as $\inprod{C}{X}$. The moment matrix $X \in \sym{\bard_2}$ contains all the monomials in $x$ of degree up to 4, hence is a linear subspace of $\sym{\bard_2}$ with dimension $\bard_4$. As a result, $X$ must satisfy $\mmom = \trinum(\bard_2) - \bard_4$ linearly independent equality constraints, referred to as the \emph{moment} constraints. In addition, $X$ must satisfy \emph{redundant} equality constraints, obtained by the fact that since $x_i^2 - 1 = 0$, it also holds $[x]_2 (x_i^2 - 1) = 0$ for each $i=1,\dots,d$. This leads to a total of $\mloc = d \times \bard_2$ linear equality constraints.
% \footnote{In~\eqref{eq:bqp}, $[x]_2$ contains $1$ and $x_i^2$, $i=1,\dots,d$, which are all equal to $1$ due to the binary constraint, hence some of the localizing constraints are linearly dependent. However, the linearly dependent constraints can be removed before solving the SDP.} 
Last but not least, the top-left entry of $X$ is equal to 1 due to the leading element of $[x]_2$ is 1 (the zero-order monomial).\footnote{The reader can refer to~\cite{Mai20arXiv-spectralPOP} or our code for details about generating SDP data $\calA,b,C$.} Therefore, the second-order relaxation for~\eqref{eq:bqp} leads to an SDP with size
\bea
n = \bard_2, \quad m = \mmom + \mloc + 1 = \trinum(\bard_2) - \bard_4 + d \times \bard_2 + 1,
\eea
which grows rapidly with $d$. Lasserre~\cite{Lasserre01ICIPCO-SDPbinary} showed that the second-order relaxation is empirically tight on a set of 50 randomly generated MAXCUT problems. However, due to the limitation of interior point methods back then, the experiments were performed on problems with small size $d=10$.
% To the best of our knowledge, the largest $d$ for which the second-order relaxation of~\eqref{eq:bqp} has been successfully solved is $d=20$ (with $n=$, $m=$) due to the limitation of interior-point methods. 

In this paper, we aim to solve \eqref{eq:bqp} instances with much larger $d$.  We generate random instances of~\eqref{eq:bqp} by sampling the coefficients vector $c$ from the standard zero-mean Gaussian distribution,~\ie~$c_i \sim \calN(0,1), i = 1,\dots,\bard_2$. At $d \in \{10,20,30,40,50,60 \}$, we randomly generate three instances each and solve the second-order moment relaxation using \sdpt, \mosek, \sdpnal, \cdcs, \sketchy, and \name. 
% We use the default parameters for all solvers, except that (i) we set the tolerance parameter in \sdpnal~to $1\ee{-8}$ to obtain more accurate solutions; (ii) we use the ``{\tt sos}'' solver in \cdcs~to exploit partial orthogonality of the SDP linear constraints and achieve slightly better convergence; (iii) we set the maximum number of iterations of~\cdcs~and~\sketchy~to $20,000$ to allow more computation time (but unfortunately still not enough to obtain accurate solutions). 
For~\name, we use the standard~\fmincon~interface in Matlab with an interior point solver (supplied with analytical objective and constraint gradients) as the \nlp~method. For \rounding~hypotheses from the moment matrix, we follow
\bea
X  = \sum_{i=1}^{n} \lambda_i v_i v_i\tran, \quad v_i \leftarrow v_i / v_i[1], \quad \barx_i = \mathrm{sgn}(v_i[x]),
\eea
which first performs a spectral decomposition of $X$ with $\lambda_1 \geq \dots \geq \lambda_n$ in non-increasing order, then normalizes the $i$-th eigenvector $v_i$ so that its leading element is $1$, and finally generates $\barx_i$ by taking the sign of the elements of $v_i$ that correspond to the order-one monomials. We round $r=5$ hypotheses using the first $5$ eigenvectors ($i=1,\dots,5$). 
% After \nlp~converges to $\hatx_i$ starting from each $x_i$, we take the best $\hatx_i$ with minimum objective value, denoted as $\hatx$, and lift it to be $\hatX = [\hatx]_2[\hatx]_2\tran$, as the lifting procedure in Algorithm~\ref{alg-iPGMnlp}. 
In order to compute $\subopt$, we set $M_b = n$ because the diagonal entries of $X$ are all equal to 1.

%!TEX root = main.tex
\begin{table}
% \vspace{-7mm}
\caption{Results on solving second-order relaxation of random~\eqref{eq:bqp} instances. ``$**$'' indicates solver out of memory.}
\label{table:pop:binary} 
\vspace{-1mm}
\adjustbox{max width=\textwidth}{%
\centering
\begin{tabular}{|c|c|c|ccccc||c|}
\hline
 Dimension & Run & Metric & \sdpt~\cite{Toh99OMS-sdpt3} & \mosek~\cite{mosek} & \cdcs~\cite{Zheng20MP-CDCS} & \sketchy~\cite{Yurtsever21SIMDS-scalableSDP} & \sdpnal~\cite{Yang15MPC-sdpnalplus} & \name \\
 \hline
 \hline
\multirow{15}{*}{\Large$ \substack{d:\ 10 \\ \\ n:\ 66 \\ \\ m:\ 1871}$ } & \multirow{5}{*}{$\#1$} 
												   & $\pfeas$ & $1.4\mathrm{e}{-10}$ & $2.3\mathrm{e}{-9}$ & $9.7\mathrm{e}{-9}$ & $0.45$ & $2.4\mathrm{e}{-9}$ & ${8.2\mathrm{e}{-16}}$ \\

							& 				       & $\dfeas$ & $4.9\mathrm{e}{-11}$ & $1.0\mathrm{e}{-9}$ & $8.4\mathrm{e}{-9}$ & $4.1\mathrm{e}{-6}$ & $7.8\mathrm{e}{-9}$ & ${7.9\mathrm{e}{-16}}$ \\
							& 					   & $\gap$   & $1.5\mathrm{e}{-8}$ & $2.8\mathrm{e}{-11}$ & $8.0\mathrm{e}{-11}$ & $0.0112$ & $8.8\mathrm{e}{-8}$ & ${5.0\mathrm{e}{-13}}$ \\
							& 					   & $\subopt$& $1.2\mathrm{e}{-12}$ & $5.6\mathrm{e}{-14}$ & $1.6\mathrm{e}{-12}$ & $3.0\mathrm{e}{-9}$ & $1.6\mathrm{e}{-12}$ & ${3.1\mathrm{e}{-16}}$ \\
							&   				   & time     & $2.2$ & $1.3$ & $8.8$ & $50.4$ & ${0.85}$ & $1.1$ \\
\cline{2-9}
							& \multirow{5}{*}{$\#2$} & $\pfeas$ & $8.1\mathrm{e}{-11}$ & $7.4\mathrm{e}{-11}$ & $8.6\mathrm{e}{-9}$ & $0.022$ & $6.7\mathrm{e}{-9}$ & ${8.2\mathrm{e}{-16}}$\\
							& 				       & $\dfeas$ & $1.2\mathrm{e}{-11}$ & $3.0\mathrm{e}{-11}$ & $9.9\mathrm{e}{-9}$ & $1.8\mathrm{e}{-6}$ & $6.0\mathrm{e}{-10}$ & ${8.8\mathrm{e}{-16}}$ \\
							& 					   & $\gap$   & $6.0\mathrm{e}{-10}$ & ${5.7\mathrm{e}{-13}}$ & $1.1\mathrm{e}{-11}$ & $0.0037$ & $6.2\mathrm{e}{-9}$ & $6.1\mathrm{e}{-13}$ \\
							& 					   & $\subopt$& $5.9\mathrm{e}{-14}$ & $2.5\mathrm{e}{-16}$ & $2.4\mathrm{e}{-13}$ & $1.7\mathrm{e}{-9}$ & $2.5\mathrm{e}{-14}$ & ${0.0}$ \\
							&   				   & time     & $1.8$ & $1.3$ & $3.0$ & $49.6$ & ${0.63}$ & $1.03$ \\
\cline{2-9}
							& \multirow{5}{*}{$\#3$} & $\pfeas$ & $2.4\mathrm{e}{-10}$ & $5.5\mathrm{e}{-12}$ & $9.2\mathrm{e}{-9}$ & $1.37$ & $2.89\mathrm{e}{-9}$ & ${8.2\mathrm{e}{-16}}$ \\
							& 				       & $\dfeas$ & $2.8\ee{-11}$ & $5.8\ee{-12}$ & $9.3\ee{-9}$ & ${1.3\ee{-14}}$ & $5.3\ee{-10}$ & $6.3\ee{-14}$ \\
							& 					   & $\gap$   & $3.5\ee{-8}$ & ${5.1\ee{-15}}$ & $3.0\ee{-11}$ & $0.384$ & $5.6\ee{-9}$ & $5.1\ee{-13}$ \\
							& 					   & $\subopt$& $2.3\ee{-10}$ & $7.8\ee{-17}$ & $4.4\ee{-12}$ & $0.0033$ & $2.1\ee{-13}$ & ${0.0}$\\
							&   				   & time     & $1.9$ & $1.4$ & $10.3$ & $49.8$ & ${0.89}$ & $0.96$ \\
\hline
\hline
\multirow{15}{*}{\Large$ \substack{d:\ 20 \\ \\ n:\ 231 \\ \\ m:\ 20,791}$ } & \multirow{5}{*}{$\#1$} 
							                       & $\pfeas$ & $3.0\ee{-9}$ & $1.1\ee{-8}$ & $1.1\ee{-5}$ & $7.34$ & $6.7\ee{-9}$ & ${1.6\ee{-15}}$ \\
							&					   & $\dfeas$ & $3.9\ee{-9}$ & $3.0\ee{-9}$ & $2.0\ee{-4}$ & $0.046$ & $9.6\ee{-11}$ & ${5.1\ee{-14}}$ \\
							& 					   & $\gap$   & $1.8\ee{-5}$ & $1.9\ee{-11}$ & $9.6\ee{-9}$ & $0.124$ & $3.3\ee{-9}$ & ${8.1\ee{-13}}$ \\
							& 					   & $\subopt$& $9.9\ee{-8}$ & $5.7\ee{-13}$ & $2.3\ee{-4}$ & $0.0076$ & $1.8\ee{-16}$ & ${1.2\ee{-16}}$ \\
							&   				   & time     & $351.8$ & $246.9$ & $119.6$ & $405.4$ & $12.55$ & ${10.67}$ \\
\cline{2-9}
							& \multirow{5}{*}{$\#2$} 
												   & $\pfeas$ & $2.4\ee{-10}$ & $1.8\ee{-12}$ & $2.9\ee{-5}$ & $4.64$ & $3.7\ee{-11}$ & ${1.61\ee{-15}}$\\
							& 				       & $\dfeas$ & $7.8\ee{-9}$ & $4.3\ee{-13}$ & $0.0013$ & $0.0438$ & $1.4\ee{-8}$ & ${8.2\ee{-14}}$ \\
							& 					   & $\gap$   & $2.1\ee{-6}$ & $6.9\ee{-15}$ & $3.9\ee{-7}$ & $0.1529$ & $2.6\ee{-7}$ & ${1.3\ee{-12}}$ \\
							& 					   & $\subopt$& $1.5\ee{-12}$ & $5.0\ee{-16}$ & $0.0013$ & $0.0047$ & $5.0\ee{-16}$ & ${9.9\ee{-17}}$ \\
							&   				   & time     & $341.5$ & $250.0$ & $119.5$ & $403.2$ & $19.6$ & ${10.4}$ \\
\cline{2-9}
							& \multirow{5}{*}{$\#3$} 
												   & $\pfeas$ & $1.6\ee{-11}$ & $4.0\ee{-11}$ & $3.9\ee{-5}$ & $2.0963$ & $2.2\ee{-10}$ & ${1.6\ee{-15}}$ \\
							& 				       & $\dfeas$ & $6.7\ee{-10}$ & $1.0\ee{-11}$ & $2.3\ee{-5}$ & $0.0409$ & $2.8\ee{-11}$ & ${3.3\ee{-14}}$ \\
							& 					   & $\gap$   & $1.3\ee{-7}$ & ${5.5\ee{-14}}$ & $3.5\ee{-8}$ & $0.0391$ & $4.2\ee{-10}$ & $5.4\ee{-13}$\\
							& 					   & $\subopt$& $4.6\ee{-13}$ & $3.4\ee{-16}$ & $1.8\ee{-5}$ & $2.9\ee{-4}$ & ${1.1\ee{-16}}$ & ${1.1\ee{-16}}$\\
							&   				   & time     & $342.8$ & $223.5$ & $117.4$ & $404.3$ & $13.0$ & ${10.0}$ \\
\hline
\hline
\multirow{15}{*}{\Large$ \substack{d:\ 30 \\ \\ n:\ 496 \\ \\ m:\ 91,761}$ } & \multirow{5}{*}{$\#1$} 
											       & $\pfeas$ &\multirow{5}{*}{**} & \multirow{5}{*}{**} & $5.0\ee{-5}$ & $7.86$ & $2.4\ee{-10}$ & ${2.4\ee{-15}}$ \\
							& 				       & $\dfeas$ & & & $0.0025$ & $0.1199$ & $5.5\ee{-9}$ & ${1.7\ee{-15}}$ \\
							& 					   & $\gap$   & & & $4.8\ee{-7}$ & $0.1850$ & $8.1\ee{-7}$ & ${2.6\ee{-12}}$ \\
							& 					   & $\subopt$& & & $0.0448$ & $0.0208$ & $5.4\ee{-15}$ & ${4.5\ee-16}$ \\
							&   				   & time     & & & $403.6$ & $2028$ & $61.9$ & ${40.1}$ \\
\cline{2-9}
							& \multirow{5}{*}{$\#2$}& $\pfeas$&\multirow{5}{*}{**} &\multirow{5}{*}{**} & $1.0\ee{-4}$ & $3.93$ & $2.8\ee{-9}$ & ${2.4\ee{-15}}$ \\
							& 				       & $\dfeas$ & & & $0.0023$ & $0.1311$ & $7.0\ee{-10}$ & ${6.8\ee{-14}}$ \\
							& 					   & $\gap$   & & & $2.8\ee{-7}$ & $0.0624$ & $1.4\ee{-8}$ & ${1.6\ee{-12}}$ \\
							& 					   & $\subopt$& & & $0.0082$ & $0.0049$ & ${0}$ & $2.3\ee{-16}$ \\
							&   				   & time     & & & $398.0$ & $2051$ & $94.2$ & ${56.7}$ \\
\cline{2-9}
							& \multirow{5}{*}{$\#3$}& $\pfeas$&\multirow{5}{*}{**} &\multirow{5}{*}{**} & $1.1\ee{-4}$ & $5.62$ & $5.0\ee{-11}$ & ${2.4\ee{-15}}$ \\
							& 				       & $\dfeas$ & & & $0.0032$ & $0.0995$ & $1.6\ee{-8}$ & ${1.0\ee{-13}}$ \\
							& 					   & $\gap$   & & & $1.0\ee{-6}$ & $0.2739$ & $1.4\ee{-6}$ & $2.9\ee{-12}$ \\
							& 					   & $\subopt$& & & $0.0115$ & $0.0084$ & $6.3\ee{-15}$ & ${6.4\ee{-17}}$ \\
							&   				   & time     & & & $405.5$ & $2026$ & $95.5$ & ${56.3}$ \\
\hline
\hline
\multirow{15}{*}{\Large$ \substack{d:\ 40 \\ \\ n:\ 861 \\ \\ m:\ 269,781}$ } & \multirow{5}{*}{$\#1$} 
												   & $\pfeas$ &\multirow{5}{*}{**} &\multirow{5}{*}{**} & $1.2\ee{-4}$ & $8.79$ & $3.4\ee{-9}$ & ${3.2\ee{-15}}$\\
							& 				       & $\dfeas$ & & & $0.0030$ & $0.1386$ & $1.5\ee{-9}$ & ${1.5\ee{-14}}$ \\
							& 					   & $\gap$   & & & $7.8\ee{-7}$ & $0.4433$ & $5.2\ee{-8}$ & ${5.7\ee{-13}}$ \\
							& 					   & $\subopt$& & & $0.0602$ & $0.0828$ & $3.1\ee{-16}$ & ${2.3\ee{-16}}$ \\
							&   				   & time     & & & $1836$ & $8962$ & $1384$ & ${461.6}$ \\
\cline{2-9}
							& \multirow{5}{*}{$\#2$}& $\pfeas$&\multirow{5}{*}{**} &\multirow{5}{*}{**} & $9.5\ee{-5}$ & $8.60$ & $3.66\ee{-10}$ & ${3.2\ee{-15}}$ \\
							& 				       & $\dfeas$ & & & $0.0027$ & $0.1340$ & $2.1\ee{-9}$ & ${5.8\ee{-14}}$ \\
							& 					   & $\gap$   & & & $6.8\ee{-7}$ & $0.3349$ & $8.1\ee{-8}$ & ${2.2\ee{-12}}$ \\
							& 					   & $\subopt$& & & $0.0395$ & $0.1746$ & ${4.0\ee{-16}}$ & ${4.0\ee{-16}}$ \\
							&   				   & time     & & & $1752$ & $8781$ & $1295$ & ${460.6}$ \\
\cline{2-9}
							& \multirow{5}{*}{$\#3$}& $\pfeas$&\multirow{5}{*}{**} &\multirow{5}{*}{**} & $2.1\ee{-4}$ & $7.23$ & $9.1\ee{-9}$ & ${3.2\ee{-15}}$ \\
							& 				       & $\dfeas$ & & & $0.0043$ & $0.1331$ & $2.8\ee{-9}$ & ${1.7\ee{-13}}$ \\
							& 					   & $\gap$   & & & $7.8\ee{-7}$ & $0.2462$ & $1.0\ee{-7}$ & ${5.7\ee{-12}}$ \\
							& 					   & $\subopt$& & & $0.0444$ & $0.0295$ & $2.9\ee{-16}$ & ${7.2\ee{-17}}$ \\
							&   				   & time     & & & $1780$ & $8845$ & $414.6$ & ${385.8}$ \\
\hline
%%%%%%%%%%%%%%%%%%%%%%%%%%%%%%%%%%%%%% Continued %%%%%%%%%%%%%%%%%%%%%%%%%%%%%%%%%%%%%%%%%%%
\hline
\multirow{15}{*}{\Large$ \substack{d:\ 50 \\ \\ n:\ 1,326 \\ \\ m:\ 629,851}$ } & \multirow{5}{*}{$\#1$} & $\pfeas$ & \multirow{5}{*}{**} & \multirow{5}{*}{**} & $5.9\ee{-5}$ & $18.07$ & $6.9\ee{-12}$ & ${4.0\ee{-15}}$ \\
							& 				       & $\dfeas$ & & & $0.0020$ & $0.5285$ & $0.0012$ & ${7.7\ee{-14}}$ \\
							& 					   & $\gap$   & & & $3.6\ee{-7}$ & $0.9033$ & $0.0133$ & ${4.0\ee{-12}}$ \\
							& 					   & $\subopt$& & & $0.0446$ & $0.9990$ & $0.0490$ & ${4.1\ee{-16}}$ \\
							&   				   & time     & & & $5232$ & $7654$ & $10000$ & ${1180}$ \\
\cline{2-9}
							& \multirow{5}{*}{$\#2$} & $\pfeas$ & \multirow{5}{*}{**} & \multirow{5}{*}{**} & $4.2\ee{-5}$ & $18.06$ & $9.2\ee{-10}$ & ${4.0\ee{-15}}$ \\
							& 				       & $\dfeas$ & & & $0.0016$ & $0.5526$ & $3.5\ee{-7}$ & ${3.0\ee{-14}}$ \\
							& 					   & $\gap$   & & & $2.3\ee{-7}$ & $0.7747$ & $0.0013$ & ${1.6\ee{-12}}$ \\
							& 					   & $\subopt$& & & $0.0532$ & $0.9394$ & $6.2\ee{-8}$ & ${5.7\ee{-17}}$ \\
							&   				   & time     & & & $5203$ & $7656$ & $5805$ & ${1440}$ \\
\cline{2-9}
							& \multirow{5}{*}{$\#3$} & $\pfeas$ & \multirow{5}{*}{**} & \multirow{5}{*}{**} & $6.8\ee{-5}$ & $17.98$ & $1.3\ee{-7}$ & ${4.0\ee{-15}}$ \\
							& 				       & $\dfeas$ & & & $0.0024$ & $0.5607$ & $5.3\ee{-8}$ & ${1.6\ee{-15}}$ \\
							& 					   & $\gap$   & & & $3.6\ee{-7}$ & $0.7953$ & $2.5\ee{-6}$ & ${8.2\ee{-13}}$ \\
							& 					   & $\subopt$& & & $0.0749$ & $0.9023$ & $7.4\ee{-14}$ & ${4.4\ee{-16}}$ \\
							&   				   & time     & & & $5312$ & $7562$ & $9884$ & ${1368}$ \\
\hline
\hline
\multirow{15}{*}{\Large$ \substack{d:\ 60 \\ \\ n:\ 1,891 \\ \\ m:\ 1,266,971}$ } & \multirow{5}{*}{$\#1$} & $\pfeas$ & \multirow{5}{*}{**} & \multirow{5}{*}{**} & $5.8\ee{-5}$ & $25.91$ & $8.3\ee{-13}$ & ${4.7\ee{-15}}$ \\
							& 				       & $\dfeas$ & & & $0.0026$ & $2.17$ & $4.7\ee{-4}$ & ${4.5\ee{-13}}$ \\
							& 					   & $\gap$   & & & $2.8\ee{-7}$ & $0.9721$ & $0.0209$ & ${3.0\ee{-11}}$\\
							& 					   & $\subopt$& & & $0.0891$ & $0.9756$ & $0.0015$ & ${1.2\ee{-15}}$\\
							&   				   & time     & & & $9731$ & $8387$ & $10000$ & ${4083}$ \\
\cline{2-9}
							& \multirow{5}{*}{$\#2$} & $\pfeas$ & \multirow{5}{*}{**} & \multirow{5}{*}{**} & $6.8\ee{-5}$ & $25.87$ & $1.4\ee{-11}$ & ${4.7\ee{-15}}$ \\
							& 				       & $\dfeas$ & & & $0.0028$ & $2.16$ & $3.5\ee{-5}$ & ${1.6\ee{-13}}$ \\
							& 					   & $\gap$   & & & $3.2\ee{-7}$ & $0.8612$ & $0.0019$ & ${1.1\ee{-11}}$\\
							& 					   & $\subopt$& & & $0.0411$ & $0.9998$ & $9.1\ee{-5}$ & ${1.0\ee{-15}}$\\
							&   				   & time     & & & $9732$ & $8361$ & $10000$ & ${3868}$ \\
\cline{2-9}
							& \multirow{5}{*}{$\#3$} & $\pfeas$ & \multirow{5}{*}{**} & \multirow{5}{*}{**} & $6.0\ee{-5}$ & $25.87$ & $0.5488$ & ${4.7\ee{-15}}$ \\
							& 				       & $\dfeas$ & & & $0.0026$ & $2.23$ & $2.3\ee{-7}$ & ${2.6\ee{-13}}$ \\
							& 					   & $\gap$   & & & $2.9\ee{-7}$ & $0.8746$ & $0.9969$ & ${1.7\ee{-11}}$ \\
							& 					   & $\subopt$& & & $0.0973$ & $0.9993$ & $0.7549$ & ${1.7\ee{-16}}$\\
							&   				   & time     & & & $10567$ & $8320$ & $10000$ & ${3804}$\\
\hline
\end{tabular}
}%
\end{table}

Table~\ref{table:pop:binary} gives the numerical results of different solvers. We make the following observations. (i) We first look at the performance of interior point methods (IPMs, {\sdpt} and {\mosek}). For small-scale problems ($d=10$), IPMs can solve the SDPs efficiently to high accuracy (\eg around $1$ second). For medium-scale problems ($d=20$), although IPMs can still obtain solutions with high accuracy, the computational time starts to grow significantly (\eg $300$-$400$ seconds). Moreover, for large-scale problems ($d\geq 30$), IPMs cannot be executed on ordinary workstations due to intensive memory consumption. The fundamental challenge of IPMs lies in solving large and dense linear systems (\ie the $m \times m$ Schur complement system) at each iteration. Although it is possible to use iterative solvers to solve the linear system \cite{toh04siopt-iterativesolver}, they are known to suffer from slow convergence as interior-point iterates approach optimality. (ii) First-order solvers ({\cdcs} and {\sketchy}) can solve problems to medium or low accuracy for $d \leq 20$, but their runtime can be worse than IPMs (although {\cdcs} is faster than {\mosek} for $d=20$, its accuracy is orders of magnitude worse than {\mosek}). However, first-order methods are indeed advantageous in terms of memory consumption and they can still be executed for problems with up to $d=60$. Nevertheless, they are not able to compute POP solutions of certified global optimality within reasonable time, as shown by the nonzero $\subopt$ in Table \ref{table:pop:binary} (\ie they are not able to show that the relaxation is indeed tight). This phenomenon further stresses the challenge for solving degenerate rank-one SDP relaxations. We also observe that, between {\cdcs} and {\sketchy}, {\cdcs} seems to perform much better for such problems. This suggests that sketching may not be the best choice for degenerate SDPs with large $m$. (iii) {\sdpnal} has the best performance among existing solvers for \eqref{eq:bqp}. It can solve \eqref{eq:bqp} instances to certified global optimality for up to $d=40$, and it is over $10$ times faster than {\mosek} and {\cdcs} when $d=20$. However, when $d=50$ and $d=60$, {\sdpnal} cannot solve the SDPs to sufficient accuracy (within $10000$ seconds), and hence the POP solution cannot be certified as globally optimal (\cf the nonzero $\subopt$). (iv) Finally, we look at the performance of our solver {\name}. We observe that {\name} solved all the SDPs to high accuracy, certified the global optimality of the POP solutions, and demonstrated the tightness of the SDP relaxations (\cf the numerically zero $\subopt$). For small and medium problems ($d=10$ and $20$), {\name} attains accuracy that is comparable to {\mosek}, while being about $30$ times faster when $d=20$. For large problems ($d\geq 30$), {\name} attains accuracy that is superior to {\sdpnal}, while being $5$-$10$ times faster when $d=50$. At $d=60$ with $m$ over a million, {\name} becomes the only solver that can obtain solutions of high accuracy.

\subsection{Quartic programming on the sphere}
\label{sec:q4s}
Consider minimizing a \emph{quartic} polynomial over the $d$-dimension unit sphere
\begin{equation} \label{eq:Q4S}
\min_{x \in \usphere{d-1}} \inprod{c}{[x]_4} \tag{Q4S}
\end{equation}
where $[x]_4: \Real{d} \rightarrow \Real{\bard_4}$ is the vector of monomials in $x$ of degree up to 4, and $c \in \Real{\bard_4}$ is the vector of known coefficients. Problem~\eqref{eq:Q4S} is known to be NP-hard with important examples such as computing the largest stable set of a graph~\cite[Theorem 3.4]{deKlerk08CEJOR-simplexhypercubesphere}, computing the $2\rightarrow 4$ norm of a matrix~\cite{Barak12ACM-hypercontractivity}, and the best separable state problem in quantum information theory~\cite{Doherty04PRA-quantumseparability}. See~\cite{Fang20MP-SOSsphere,Ling10SIOPT-biquadratic} and references therein for a thorough discussion about problem~\eqref{eq:Q4S}.

Here we consider the dense second-order (also the lowest order) moment relaxation of~\eqref{eq:Q4S} and numerically show that they are indeed tight and admit rank-one solutions. By following the same relaxation scheme as in Section~\ref{sec:bqp} (\ie build the moment matrix $X = [x]_2 [x]_2\tran$ and add equality constraints), we can count the size of the SDP relaxation to be
\bea
n = \bard_2, \quad m = \trinum(\bard_2) - \bard_4 + \bard_2 + 1.
\eea

At each $d \in \{10,20,30,40,50,60\}$, we generate three random instances of \eqref{eq:Q4S} by drawing $c$ from the standard normal distribution. We solve the resulting SDP using \sdpt, \mosek, \cdcs, \sketchy, \sdpnal, and \name. 
% We use the same parameters for all solvers as in Section~\ref{sec:bqp}. 
For~\name, we exploit the manifold structure of the sphere constraint and adopt~\manopt~\cite{manopt} with a trust region solver as the~\nlp~method. One can also treat~\eqref{eq:Q4S} as a standard nonlinear programming and solve it with~\fmincon, but we found that~\manopt~is faster and more robust for this problem. To generate hypotheses for \nlp~from the moment matrix, we follow
\bea
X  = \sum_{i=1}^{n} \lambda_i v_i v_i\tran, \quad v_i \leftarrow v_i / v_i[1], \quad \barx_i = \frac{v_i[x]}{\norm{v_i[x]}},
\eea
which first performs a spectral decomposition of $X$ with $\lambda_1 \geq \dots \geq \lambda_n$ in non-increasing order, then normalizes the $i$-th eigenvector $v_i$ so that its leading element is $1$, and finally generates $\barx_i$ by projecting the elements of $v_i$ that correspond to the order-one monomials onto the unit sphere. We generate $r=5$ hypotheses using the first $5$ eigenvectors. A valid upper bound $M_b$ on $\trace{X}$ can be obtained as 
\bea\nonumber
\begin{split}
\trace{X} & = \trace{[x]_2 [x]_2\tran} = [x]_2\tran [x]_2 \\
& = 1 + \sum_{i=1}^d x_i^2 + \sum_{1\leq i \leq j \leq d} (x_i x_j)^2  \leq 2 + \parentheses{\sum_{i=1}^d x_i^2}^2 = 3 := M_b.
\end{split}
\eea

%!TEX root = main.tex
\begin{table}
% \vspace{-7mm}
\caption{Results on solving second-order relaxation of random~\eqref{eq:Q4S} instances. ``$**$'' indicates solver out of memory.}
\label{table:pop:sphere} 
\vspace{-1mm}
\adjustbox{max width=\textwidth}{%
\centering
\begin{tabular}{|c|c|c|ccccc||c|}
\hline
 Dimension & Run & Metric & \sdpt~\cite{Toh99OMS-sdpt3} & \mosek~\cite{mosek} & \cdcs~\cite{Zheng20MP-CDCS} & \sketchy~\cite{Yurtsever21SIMDS-scalableSDP} & \sdpnal~\cite{Yang15MPC-sdpnalplus} & \name \\
 \hline
 \hline
\multirow{15}{*}{\Large$ \substack{d:\ 10 \\ \\ n:\ 66 \\ \\ m:\ 1277}$ } & \multirow{5}{*}{$\#1$} 
												   & $\pfeas$ & $4.5\ee{-9}$ & $3.3\ee{-8}$ & $1.2\ee{-11}$ & $0.1388$ & $4.6\ee{-12}$ & $1.5\ee{-16}$ \\
							& 				       & $\dfeas$ & $5.1\ee{-11}$ & $5.1\ee{-9}$ & $9.7\ee{-11}$ & $8.7\ee{-4}$ & $1.9\ee{-11}$ & $3.5\ee{-15}$ \\
							& 					   & $\gap$   & $5.7\ee{-9}$ & $3.4\ee{-9}$ & $9.4\ee{-13}$ & $0.1823$ & $6.9\ee{-12}$ & $6.8\ee{-10}$ \\
							& 					   & $\subopt$& $3.6\ee{-12}$ & $3.2\ee{-11}$ & $1.4\ee{-12}$ & $0.0116$ & $2.1\ee{-11}$ & $2.0\ee{-11}$ \\
							&   				   & time     & $1.2$ & $0.6$ & $2.5$ & $58.0$ & $0.6$ & $1.2$ \\
\cline{2-9}
							& \multirow{5}{*}{$\#2$} & $\pfeas$ & $1.7\ee{-9}$ & $1.3\ee{-8}$ & $8.8\ee{-11}$ & $0.0291$ & $1.8\ee{-11}$ & $1.6\ee{-16}$ \\
							& 				       & $\dfeas$ & $1.4\ee{-10}$ & $2.0\ee{-9}$ & $2.4\ee{-11}$ & $6.0\ee{-4}$ & $1.7\ee{-11}$ & $1.0\ee{-11}$\\
							& 					   & $\gap$   & $6.5\ee{-9}$ & $2.1\ee{-10}$ & $8.2\ee{-11}$ & $0.0244$ & $1.7\ee{-11}$ & $5.0\ee{-11}$ \\
							& 					   & $\subopt$& $1.3\ee{-12}$ & $6.5\ee{-10}$ & $4.0\ee{-12}$ & $1.8\ee{-4}$ & $2.4\ee{-13}$ & $2.9\ee{-12}$ \\
							&   				   & time     & $1.0$ & $0.8$ & $0.7$ & $57.7$ & $0.6$ & $1.0$ \\
\cline{2-9}
							& \multirow{5}{*}{$\#3$} & $\pfeas$ & $2.8\ee{-9}$ & $1.6\ee-{12}$ & $9.8\ee{-11}$ & $0.0427$ & $2.5\ee{-11}$ & $1.4\ee{-16}$ \\
							& 				       & $\dfeas$ & $1.5\ee{-11}$ & $9.1\ee{-11}$ & $3.2\ee{-11}$ & $6.4\ee{-4}$ & $3.6\ee{-11}$ & $5.5\ee{-12}$ \\
							& 					   & $\gap$   & $1.8\ee{-10}$ & $1.1\ee{-13}$ & $8.5\ee{-11}$ & $0.0384$ & $3.4\ee{-11}$ & $2.1\ee{-11}$ \\
							& 					   & $\subopt$& $2.2\ee{-12}$ & $3.9\ee{-15}$ & $3.1\ee{-12}$ & $2.3\ee{-4}$ & $1.4\ee{-12}$ & $1.9\ee{-12}$\\
							&   				   & time     & $1.0$ & $0.6$ & $0.6$ & $58.2$ & $0.6$ & $1.1$ \\
\hline
\hline
\multirow{15}{*}{\Large$ \substack{d:\ 20 \\ \\ n:\ 231 \\ \\ m:\ 16,402}$ } & \multirow{5}{*}{$\#1$} 
							                       & $\pfeas$ & $3.0\ee{-11}$ & $2.2\ee{-9}$ & $3.6\ee{-12}$ & $0.0800$ & $6.1\ee{-11}$ & $1.0\ee{-16}$ \\
							&					   & $\dfeas$ & $1.3\ee{-10}$ & $1.3\ee{-10}$ & $9.9\ee{-11}$ & $4.9\ee{-4}$ & $2.9\ee{-11}$ & $4.6\ee{-15}$\\
							& 					   & $\gap$   & $3.4\ee{-6}$ & $4.3\ee{-11}$ & $2.3\ee{-11}$ & $0.0875$ & $6.4\ee{-11}$ & $2.7\ee{-11}$ \\
							& 					   & $\subopt$& $4.1\ee{-9}$ & $4.7\ee{-13}$ & $2.5\ee{-12}$ & $4.3\ee{-4}$ & $1.6\ee{-12}$ & $2.5\ee{-12}$ \\
							&   				   & time     & $173.8$ & $99.6$ & $11.2$ & $566.3$ & $2.9$ & $4.4$ \\
\cline{2-9}
							& \multirow{5}{*}{$\#2$} 
												   & $\pfeas$ & $1.5\ee{-11}$ & $2.1\ee{-8}$ & $2.2\ee{-12}$ & $0.1975$ & $1.7\ee{-11}$ & $1.2\ee{-16}$ \\
							& 				       & $\dfeas$ & $7.2\ee{-11}$ & $1.2\ee{-9}$ & $1.0\ee{-10}$ & $0.0010$ & $7.1\ee{-11}$ & $7.9\ee{-14}$ \\
							& 					   & $\gap$   & $2.9\ee{-6}$ & $1.1\ee{-9}$ & $2.2\ee{-12}$ & $0.2571$ & $1.5\ee{-11}$ & $5.5\ee{-13}$\\
							& 					   & $\subopt$& $1.0\ee{-8}$ & $1.5\ee{-10}$ & $4.6\ee{-13}$ & $0.0207$ & $7.4\ee{-13}$ & $8.6\ee{-14}$ \\
							&   				   & time     & $192.0$ & $98.7$ & $30.7$ & $553.3$ & $3.0$ & $14.7$\\
\cline{2-9}
							& \multirow{5}{*}{$\#3$} 
												   & $\pfeas$ & $2.2\ee{-10}$ & $7.7\ee{-9}$ & $9.4\ee{-11}$ & $0.0373$ & $8.8\ee{-11}$ & $2.7\ee{-16}$ \\
							& 				       & $\dfeas$ & $1.2\ee{-10}$ & $4.5\ee{-10}$ & $4.5\ee{-11}$ & $6.3\ee{-4}$ & $1.9\ee{-11}$ & $4.5\ee{-13}$ \\
							& 					   & $\gap$   & $2.2\ee{-6}$ & $9.3\ee{-11}$ & $9.1\ee{-11}$ & $0.0323$ & $1.3\ee{-10}$ & $4.5\ee{-13}$ \\
							& 					   & $\subopt$& $2.6\ee{-9}$ & $5.1\ee{-11}$ & $1.7\ee{-11}$ & $3.8\ee{-4}$ & $1.5\ee{-12}$ & $2.7\ee{-13}$ \\
							&   				   & time     & $174.8$ & $90.3$ & $9.0$ & $551.7$ & $2.7$ & $4.5$ \\
\hline
\hline
\multirow{15}{*}{\Large$ \substack{d:\ 30 \\ \\ n:\ 496 \\ \\ m:\ 77,377}$ } & \multirow{5}{*}{$\#1$} 
											       & $\pfeas$ &\multirow{5}{*}{**} & \multirow{5}{*}{**} & $8.8\ee{-13}$ & $0.2407$ & $4.4\ee{-13}$ & $2.1\ee{-16}$ \\
							& 				       & $\dfeas$ & & & $1.0\ee{-10}$ & $7.3\ee{-4}$ & $1.0\ee{-12}$ & $5.0\ee{-15}$ \\
							& 					   & $\gap$   & & & $1.2\ee{-11}$ & $0.3913$ & $9.6\ee{-10}$ & $7.9\ee{-12}$ \\
							& 					   & $\subopt$& & & $2.1\ee{-12}$ & $0.0037$ & $9.5\ee{-13}$ & $1.3\ee{-12}$ \\
							&   				   & time     & & & $202.3$ & $1092$ & $15.1$ & $17.3$ \\
\cline{2-9}
							& \multirow{5}{*}{$\#2$}& $\pfeas$&\multirow{5}{*}{**} &\multirow{5}{*}{**} & $8.5\ee{-13}$ & $0.2198$ & $5.8\ee{-11}$ & $1.8\ee{-16}$ \\
							& 				       & $\dfeas$ & & & $1.0\ee{-10}$ & $7.3\ee{-4}$ & $7.4\ee{-11}$ & $4.4\ee{-14}$ \\
							& 					   & $\gap$   & & & $4.7\ee{-13}$ & $0.3468$ & $7.5\ee{-11}$ & $5.7\ee{-13}$ \\
							& 					   & $\subopt$& & & $1.0\ee{-13}$ & $0.0070$ & $2.6\ee{-13}$ & $1.0\ee{-13}$ \\
							&   				   & time     & & & $98.0$ & $1052$ & $10.6$ & $18.5$ \\
\cline{2-9}
							& \multirow{5}{*}{$\#3$}& $\pfeas$&\multirow{5}{*}{**} &\multirow{5}{*}{**} & $8.8\ee{-13}$ & $0.2637$ & $1.5\ee{-12}$ & $1.5\ee{-16}$ \\
							& 				       & $\dfeas$ & & & $1.0\ee{-10}$ & $1.3\ee{-4}$ & $8.3\ee{-11}$ & $1.1\ee{-12}$ \\
							& 					   & $\gap$   & & & $5.2\ee{-12}$ & $0.4770$ & $1.3\ee{-11}$ & $1.5\ee{-11}$ \\
							& 					   & $\subopt$& & & $9.7\ee{-13}$ & $0.0196$ & $2.4\ee{-12}$ & $2.7\ee{-12}$ \\
							&   				   & time     & & & $224.0$ & $1046$ & $12.3$ & $16.1$ \\
\hline
\hline
\multirow{15}{*}{\Large$ \substack{d:\ 40 \\ \\ n:\ 861 \\ \\ m:\ 236,202}$ } & \multirow{5}{*}{$\#1$} 
												   & $\pfeas$ &\multirow{5}{*}{**} &\multirow{5}{*}{**} & $1.3\ee{-12}$ & $0.3876$ & $2.4\ee{-12}$ & $2.9\ee{-16}$ \\
							& 				       & $\dfeas$ & & & $1.0\ee{-10}$ & $1.8\ee{-14}$ & $4.4\ee{-11}$ & $5.2\ee{-15}$ \\
							& 					   & $\gap$   & & & $5.0\ee{-12}$ & $0.7867$ & $7.4\ee{-12}$ & $3.9\ee{-12}$ \\
							& 					   & $\subopt$& & & $1.0\ee{-12}$ & $0.0804$ & $1.1\ee{-12}$ & $6.8\ee{-13}$ \\
							&   				   & time     & & & $413.0$ & $3005$ & $32.7$ & $46.8$ \\
\cline{2-9}
							& \multirow{5}{*}{$\#2$}& $\pfeas$&\multirow{5}{*}{**} &\multirow{5}{*}{**} & $5.1\ee{-7}$ & $0.3401$ & $9.4\ee{-11}$ & $2.1\ee{-16}$ \\
							& 				       & $\dfeas$ & & & $8.3\ee{-6}$ & $2.0\ee{-14}$ & $7.2\ee{-11}$ & $5.5\ee{-13}$ \\
							& 					   & $\gap$   & & & $1.0\ee{-5}$ & $0.7936$ & $5.9\ee{-10}$ & $4.2\ee{-12}$\\
							& 					   & $\subopt$& & & $0.0063$ & $0.1052$ & $7.5\ee{-12}$ & $7.2\ee{-12}$ \\
							&   				   & time     & & & $1690$ & $3034$ & $50.5$ & $76.7$ \\
\cline{2-9}
							& \multirow{5}{*}{$\#3$}& $\pfeas$&\multirow{5}{*}{**} &\multirow{5}{*}{**} & $9.5\ee{-13}$ & $0.3397$ & $1.5\ee{-11}$ & $2.1\ee{-16}$ \\
							& 				       & $\dfeas$ & & & $1.0\ee{-10}$ & $2.1\ee{-14}$ & $5.4\ee{-11}$ & $5.2\ee{-15}$ \\
							& 					   & $\gap$   & & & $6.1\ee{-12}$ & $0.8206$ & $3.1\ee{-11}$ & $6.0\ee{-11}$ \\
							& 					   & $\subopt$& & & $1.1\ee{-12}$ & $0.1558$ & $3.1\ee{-12}$ & $9.8\ee{-12}$ \\
							&   				   & time     & & & $622.8$ & $3045$ & $32.1$ & $73.4$ \\
\hline
%%%%%%%%%%%%%%%%%%%%%%%%%%%%%%%%%%%%%% Continued %%%%%%%%%%%%%%%%%%%%%%%%%%%%%%%%%%%%%%%%%%%
\hline
\multirow{15}{*}{\Large$ \substack{d:\ 50 \\ \\ n:\ 1,326 \\ \\ m:\ 564,877}$ } & \multirow{5}{*}{$\#1$} & $\pfeas$ & \multirow{5}{*}{**} & \multirow{5}{*}{**} & $8.0\ee{-11}$ & $0.3188$ & $2.6\ee{-13}$ & $2.0\ee{-16}$ \\
							& 				       & $\dfeas$ & & & $1.4\ee{-8}$ & $3.2\ee{-14}$ & $3.4\ee{-10}$ & $5.8\ee{-15}$ \\
							& 					   & $\gap$   & & & $8.4\ee{-10}$ & $0.8868$ & $1.6\ee{-9}$ & $2.7\ee{-11}$ \\
							& 					   & $\subopt$& & & $4.0\ee{-8}$ & $0.3756$ & $7.2\ee{-11}$ & $4.8\ee{-12}$ \\
							&   				   & time     & & & $4201$ & $5364$ & $156.5$ & $240.2$ \\
\cline{2-9}
							& \multirow{5}{*}{$\#2$} & $\pfeas$ & \multirow{5}{*}{**} & \multirow{5}{*}{**} & $3.9\ee{-12}$ & $0.3341$ & $1.1\ee{-12}$ & $2.6\ee{-16}$ \\
							& 				       & $\dfeas$ & & & $1.0\ee{-10}$ & $3.3\ee{-14}$ & $6.1\ee{-11}$ & $5.8\ee{-15}$ \\
							& 					   & $\gap$   & & & $3.0\ee{-11}$ & $0.8910$ & $6.9\ee{-12}$ & $1.7\ee{-12}$ \\
							& 					   & $\subopt$& & & $6.4\ee{-12}$ & $0.4208$ & $1.2\ee{-12}$ & $3.1\ee{-13}$ \\
							&   				   & time     & & & $1015$ & $5356$ & $75.9$ & $138.5$ \\
\cline{2-9}
							& \multirow{5}{*}{$\#3$} & $\pfeas$ & \multirow{5}{*}{**} & \multirow{5}{*}{**} & $9.4\ee{-13}$ & $0.3199$ & $5.3\ee{-14}$ & $1.8\ee{-16}$\\
							& 				       & $\dfeas$ & & & $1.0\ee{-10}$ & $3.1\ee{-14}$ & $5.7\ee{-11}$ & $5.9\ee{-15}$ \\
							& 					   & $\gap$   & & & $2.4\ee{-11}$ & $0.9205$ & $6.5\ee{-9}$ & $1.9\ee{-12}$ \\
							& 					   & $\subopt$& & & $4.2\ee{-12}$ & $0.4937$ & $1.4\ee{-10}$ & $3.2\ee{-13}$ \\
							&   				   & time     & & & $1284$ & $5356$ & $113.3$ & $145.9$ \\
\hline
\hline
\multirow{15}{*}{\Large$ \substack{d:\ 60 \\ \\ n:\ 1,891 \\ \\ m:\ 1,155,402}$ } & \multirow{5}{*}{$\#1$} & $\pfeas$ & \multirow{5}{*}{**} & \multirow{5}{*}{**} & $2.1\ee{-13}$ & $0.3276$ & $3.4\ee{-13}$ & $2.1\ee{-16}$ \\
							& 				       & $\dfeas$ & & & $1.0\ee{-10}$ & $3.8\ee{-14}$ & $1.8\ee{-10}$ & $6.2\ee{-15}$  \\
							& 					   & $\gap$   & & & $4.1\ee{-12}$ & $0.9610$ & $5.2\ee{-9}$ & $3.0\ee{-12}$ \\
							& 					   & $\subopt$& & & $7.9\ee{-13}$ & $0.6813$ & $8.2\ee{-10}$ & $5.7\ee{-13}$ \\
							&   				   & time     & & & $5059$ & $6895$ & $321.9$ & $283.4$ \\
\cline{2-9}
							& \multirow{5}{*}{$\#2$} & $\pfeas$ & \multirow{5}{*}{**} & \multirow{5}{*}{**} & $3.1\ee{-13}$ & $0.3301$ & $4.0\ee{-14}$ & $4.2\ee{-16}$\\
							& 				       & $\dfeas$ & & & $1.0\ee{-10}$ & $5.2\ee{-14}$ & $9.9\ee{-11}$ & $6.3\ee{-15}$ \\
							& 					   & $\gap$   & & & $3.1\ee{-12}$ & $0.9281$ & $1.6\ee{-10}$ & $5.5\ee{-14}$ \\
							& 					   & $\subopt$& & & $7.1\ee{-10}$ & $0.5856$ & $3.5\ee{-11}$ & $8.8\ee{-15}$ \\
							&   				   & time     & & & $7684$ & $6905$ & $286.3$ & $303.5$ \\
\cline{2-9}
							& \multirow{5}{*}{$\#3$} & $\pfeas$ & \multirow{5}{*}{**} & \multirow{5}{*}{**} & $6.2\ee{-13}$ & $0.3133$ & $7.9\ee{-15}$ & $4.5\ee{-16}$ \\
							& 				       & $\dfeas$ & & & $1.0\ee{-10}$ & $4.0\ee{-14}$ & $6.1\ee{-11}$ & $6.1\ee{-15}$ \\
							& 					   & $\gap$   & & & $3.5\ee{-12}$ & $0.9572$ & $4.3\ee{-9}$ & $3.6\ee{-11}$ \\
							& 					   & $\subopt$& & & $6.9\ee{-13}$ & $0.5980$ & $2.6\ee{-10}$ & $6.6\ee{-12}$ \\
							&   				   & time     & & & $3669$ & $6909$ & $268.9$ & $311.5$ \\
\hline
\end{tabular}
}%
\end{table}

Table \ref{table:pop:sphere} gives the numerical results of different solvers. We make the following observations. (i) Similar to Table \ref{table:pop:binary} for the \eqref{eq:bqp} problem, IPMs can solve small and medium problems ($d=10$ and $20$) to high accuracy, although the runtime grows significantly from $d=10$ (about $1$ second) to $d=20$ (about $100$ seconds). (ii) Both {\cdcs} and {\sdpnal} are able to solve all SDPs to high accuracy and certify the tightness of the second-order relaxation (despite that {\cdcs} only attained medium accuracy for $\#2$ at $d=40$). However, {\sdpnal} is significantly faster than {\cdcs}, in most cases $10$-$30$ times faster. Compared to Table \ref{table:pop:binary}, these results suggest that the \eqref{eq:Q4S} relaxation is easier to solve than the \eqref{eq:bqp} relaxation, perhaps because \eqref{eq:Q4S} only has a single unit-norm constraint. (iii) {\sketchy}, however, failed to solve most of the SDPs to high accuracy, despite taking more time than {\cdcs} and {\sdpnal}. (iv) Our solver {\name} achieved similar performance compared to {\sdpnal}. Although {\name} can be slightly slower than {\sdpnal}, it generally attained higher accuracy than {\sdpnal}. 

% \red{The results presented in the Table~\ref{table:pop:sphere} shows that for small-scale problems (say, $m$ is less than 20 thousands), except for \cgal, all other solvers are able to solve the problems successfully. However, \sdpnal and \name turn out to have comparable and more promising performance. For large-scale problems (say, $m$ is larger 70 thousands), \cgal is scalable but not able to compute satisfactory optimal solutions. On the other hand, \cdcs, \sdpnal and \name all able to solve the problem to very high accuracy. Again, \sdpnal and \name have similar performance and outperform \cdcs.}

%%%%%%%%%%%%%%%%%%%%%%%%%%%%%%%%%%%%%%%%%%%%%%%%%%%%%%%%%%%%%%%%%%%%%%%%%%%%%%%%%%%%%%%%%%%%%%%%%%%%%%%%%%%%%%%%%%%%%%%
%%%%%%%%%%%%%%%%%%%%%%%%%%%%%%%%%%%%%%%%%%%%%%%%%%%%%%%%%%%%%%%%%%%%%%%%%%%%%%%%%%%%%%%%%%%%%%%%%%%%%%%%%%%%%%%%%%%%%%%

\subsection{Outlier-robust Wahba problem}
\label{sec:wahba}
Consider the problem of finding the best 3D rotation to align two sets of 3D points while explicitly tolerating \emph{outliers}
\bea \label{eq:tlswahba}
\min_{q \in \usphere{3}} \sum_{i=1}^N \min \cbrace{ \frac{\norm{\tz_i - q \circ \tw_i \circ q\inv }^2 }{\beta_i^2}, 1 }
\eea
where $q \in \usphere{3}$ is the \emph{unit quaternion} parametrization of a 3D rotation, $(z_i \in \Real{3},w_i\in\Real{3})_{i=1}^N$ are given $N$ pairs of 3D points (often normalized to have unit norm), $\tilde{z} \triangleq [z\tran,0]\tran \in \Real{4}$ denotes the zero-homogenization of a 3D vector $z$, $q \inv \triangleq [-q_1,-q_2,-q_3,q_4]\tran$ is the inverse quaternion, ``$\circ$'' denotes the quaternion product defined as
\bea
q \circ p \triangleq \bmat{cccc}
q_4 & -q_3 & q_2 & q_1 \\
q_3 & q_4 & -q_1 & q_2 \\
-q_2 & q_1 & q_4 & q_3 \\
-q_1 & -q_2 & -q_3 & q_4
\emat p, \quad \forall q,p \in \Real{4},
\eea
$\beta_i > 0$ is a given threshold that determines the maximum \emph{inlier} residual, and $\min \cbrace{\cdot, \cdot}$ realizes the so-called \emph{truncated least squares} (TLS) cost function in robust estimation~\cite{Antonante20arxiv-outlier}. Intuitively, the term $q \circ \tw_i \circ q\inv$ is the rotated copy of $w_i$, and the $\ell_2$ norm in~\eqref{eq:tlswahba} measures the Euclidean distance between $z_i$ and $w_i$ after rotation (a metric for the goodness of fit). Problem~\eqref{eq:tlswahba} therefore seeks to find the best 3D rotation that minimizes the sum of (normalized) squared Euclidean distances between $z_i$ and $w_i$ while preventing outliers from damaging the estimation via the usage of the TLS cost function, which assigns a constant value to those pairs of points that cannot be aligned well (\ie~outliers). A pictorial description of the outlier-robust Wahba problem is presented in Fig.~\ref{fig:wahba}. Problem~\eqref{eq:tlswahba} is nonsmooth, but can be equivalently reformulated as
\begin{equation}\label{eq:wahba}
\min_{ \substack{ q \in \usphere{3}, \\ \theta_i \in \{+1,-1\}, i=1,\dots,N} } \sum_{i=1}^N \frac{1+\theta_i}{2} \frac{\norm{\tz_i - q \circ \tw_i \circ q\inv }^2 }{\beta_i^2} + \frac{1-\theta_i}{2} \tag{Wahba}
\end{equation}
by introducing $N$ binary variables $\{ \theta_i \}$ that expose the combinatorial nature. Each $\theta_i$ acts as the selection variable for determining \blue{whether} the $i$-th pair of 3D points $(z_i,w_i)$ is an inlier or an outlier. Problem \eqref{eq:wahba} is a fundamental problem in aerospace, robotics and computer vision, and is the rotation subproblem in point cloud registration~\cite{Yang20TRO-teaser,Yang19rss-teaser}.

%!TEX root = main.tex

\begin{figure}[ht!]
\centering
\includegraphics[width=0.7\linewidth]{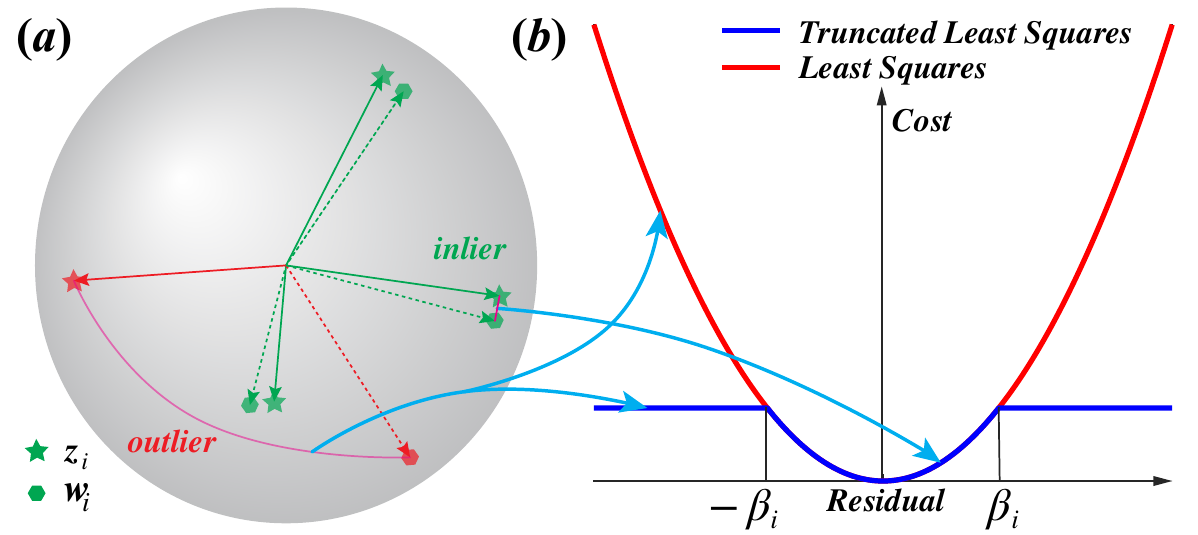}
\caption{An example of the outlier-robust~\eqref{eq:wahba} problem. (a) Four pairs of 3D points $(z_i,w_i)$ lying around a unit sphere, with one of the pairs being an outlier that cannot be aligned well by a 3D rotation. (b) A truncated least squares (TLS) cost function compared with a least squares cost function. The TLS cost function prevents the outlier from contaminating the estimation problem by assigning a constant cost to the outlier. Adapted from~\cite{Yang19ICCV-wahba}.}
\label{fig:wahba}
\end{figure}

To solve~\eqref{eq:wahba} to global optimality, Yang and Carlone~\cite{Yang19ICCV-wahba} proposed the following semidefinite relaxation that was empirically shown to be always tight. Let $x = [q\tran,\theta_1,\dots,\theta_N]\tran \in \Real{d}, d = 4+N$, be the variable of the nonlinear programming \blue{problem}~\eqref{eq:wahba}, construct 
\bea \label{eq:wahbaclone}
[x]_s = [q\tran,\theta_1 q\tran, \dots, \theta_N q\tran]\tran \in \Real{n}, \quad n = 4N+4
\eea
as the sparse set of monomials in $x$ of degree up to 2 (a technique that was dubbed \emph{binary cloning}), and then build $X = [x]_s [x]_s\tran$ as the \emph{sparse} moment matrix. Because of the binary constraint $\theta_i^2 = 1$, it can be easily seen that: (i) the diagonal $4\times 4$ blocks of $X$ are all identical ($\theta_i^2 q q\tran = qq\tran$), and (ii) the off-diagonal $4\times 4$ blocks are symmetric ($\theta_i \theta_j q q\tran \in \sym{4}$). Because of the unit quaternion constraint, $X$ satisfies $\trace{X} = N+1$. Therefore, this leads to a semidefinite relaxation of size
\bea
n = 4N+4, \quad m = 10N + 1 + 3N(N+1).
\eea
Compared to the dense second-order moment relaxation in Section~\ref{sec:bqp} and~\ref{sec:q4s}, this sparse second-order relaxation is much more manageable, and the largest $N$ whose relaxation was successfully solved by interior point method was $N=100$~\cite{Yang19ICCV-wahba}. This sparse second-order relaxation scheme has been shown as a general framework for certifiable outlier-robust machine perception~\cite{Yang20neurips-onering}.
%,Yang21arxiv-certifiable}.

Here we show the scalability of our solver by solving instances of~\eqref{eq:wahba} up to $N=1000$ and obtaining the globally optimal solution. At each $N = 50,100,200,500$, $1000$, we generate three random instances of the~\eqref{eq:wahba} problem as follows. (i) We draw a random 3D rotation $R \in \SOthree$ (a rotation matrix can be converted from and to a unit quaternion easily); (ii) we simulate $N$ 3D unit vectors $w_i,i=1,\dots,N$ uniformly on the unit sphere; (iii) we generate 
\bea\label{eq:wahbasimmodel}
z_i = R w_i + \epsilon_i, \quad \epsilon_i \sim \calN(0,0.01^2),i=1,\dots,N
\eea
by rotating $w_i$ and adding Gaussian noise; (iv) we replace $50\%$ of the $z_i$'s by random unit vectors on the sphere so that they do not follow the generative model \eqref{eq:wahbasimmodel} and are considered as outliers. We then use \sdpt, \mosek, \cdcs, \sketchy, \sdpnal, and \name~to solve the SDP relaxations. For \name, we use \manopt~with a trust region solver as the \nlp~method for solving the nonlinear programming~\eqref{eq:wahba}. Specifically, $q\in \usphere{3}$ is modeled as a sphere manifold, and $\theta \in \{+1,-1\}^N$ is modeled as an oblique manifold of size $1 \times N$ (an oblique manifold of size $n\times m$ is the set of matrices of size $n \times m$  with unit-norm columns), and the problem is treated as an unconstrained problem on the product of two manifolds. To round hypotheses from a moment matrix $X$, we follow
\bea
X = \sum_{i=1}^n \lambda_i v_i v_i\tran, \quad q_i = \frac{v_i[q]}{\norm{v_i[q]}}, \quad \theta_j^i = \sign(q_i\tran v_i[\theta_j q]),\;
j=1,\dots,N,
\eea
where we first perform spectral decomposition of $X$ with eigenvalues in nonincreasing order, then round $q_i$ by normalizing the corresponding entries of $v_i$ to have unit norm, and finally generate $\theta_j^i$ by taking the sign of the dot product between the rounded $q_i$ and the entries of $v_i$ corresponding to each $\theta_j q$ block (the rationale for using this \rounding~method is easily seen from~\eqref{eq:wahbaclone} \blue{where 
$\theta_j^i$ is identified with $q_i\tran (\theta_j q_i)$ for the rounded $q_i$ with unit norm}). We generate $r=5$ hypotheses by rounding $5$ eigenvectors from the moment matrix. We set $M_b=N+1 = \trace{X}$ to compute $\subopt$ as in \eqref{eq:relsubopt}. 

Table \ref{table:wahba} gives the numerical results for different solvers. Notice that at $N=1000$, we increased the maximum runtime of {\sketchy} and {\sdpnal} to be $50000$ seconds for a fair comparison with {\name}. We make the following observations. (i) IPMs can solve small and medium problems ($N=50$ and $100$) to high accuracy and certify global optimality of the POP solutions. However, their runtime grows quickly and they cannot scale to problems with $N \geq 200$. (ii) {\cdcs}, {\sketchy} and {\sdpnal} perform poorly on this problem. Notably, {\cdcs} and {\sketchy} failed on all instances and they cannot certify global optimality and tightness of the relaxation. {\sdpnal} succeeded on problems with $N=50$ but failed to attain high accuracy for all other problems. Comparing Table \ref{table:wahba} with Tables \ref{table:pop:binary}-\ref{table:pop:sphere}, the degraded performance of {\cdcs} and {\sdpnal} seems to suggest that sparse relaxations are more challenging to solve than dense relaxations. (iii) {\name} was able to solve all SDP instances to high accuracy. Particularly, for $N=50$ and $N=100$, {\name} achieved similar accuracy compared to {\mosek}, while being $3$ times faster at $N=50$ and $20$ times faster at $N=100$. For $N \geq 200$, {\name} is the only solver than can attain high accuracy and certify global optimality and tightness (despite taking much less time than the other solvers). 

{\bf {\name} with Domain-Specific Primal Initialization}. In Section \ref{subsec:initialpoints}, we mentioned that {\name} can benefit from domain-specific primal initialization. We now use the \eqref{eq:wahba} problem to support our claim. Although the \eqref{eq:wahba} problem is a combinatorial problem with binary variables, Yang \etal \cite{Yang20ral-GNC} have designed a heuristic method called \emph{graduated non-convexity} (\gnc) that can solve \eqref{eq:wahba} to global optimality with high probability of success (note that {\gnc} only outputs a solution without optimality certificate). Therefore, we first use {\gnc} to solve the combinatorial \eqref{eq:wahba} problem and then use its solution as a primal initialization for {\name}. Particularly, let $(\hatq,\hattheta)$ be the output of {\gnc}, we input a rank-one point $X^0 = [\hatx]_s [\hatx]_s\tran$ to {\name}, where $[\hatx]_s$ is computed from \eqref{eq:wahbaclone} with $\hatx = (\hatq,\hattheta)$. The last column of Table \ref{table:wahba} shows the numerical results for {\name} with primal initialization supplied by {\gnc}. We can see that the {\gnc} primal initialization gives {\name} an additional $2$-$3$ times speedup. 

%!TEX root = main.tex
\begin{table}
% \vspace{-7mm}
\caption{Results on solving sparse second-order relaxation of random~\eqref{eq:wahba} instances. ``$**$'' indicates solver out of memory. The last column shows results for {\name} with primal initialization using graduated non-convexity (\gnc) \cite{Yang20ral-GNC}.}
\label{table:wahba} 
\vspace{-1mm}
\adjustbox{max width=\textwidth}{%
\centering
\begin{tabular}{|c|c|c|ccccc||cc|}
\hline
 Dimension & Run & Metric & \sdpt~\cite{Toh99OMS-sdpt3} & \mosek~\cite{mosek} & \cdcs~\cite{Zheng20MP-CDCS} & \sketchy~\cite{Yurtsever21SIMDS-scalableSDP} & \sdpnal~\cite{Yang15MPC-sdpnalplus} & \name & w/ \gnc\\
 \hline
 \hline
\multirow{15}{*}{\Large$ \substack{N:\ 50 \\ \\ n:\ 204 \\ \\ m:\ 8,151}$ } & \multirow{5}{*}{$\#1$} 
												   & $\pfeas$ & $4.9\ee{-11}$ & $1.6\ee{-13}$ & $2.4\ee{-6}$ & $0.8356$ & $4.6\ee{-14}$ & $1.4\ee{-16}$ & $1.2\ee{-15}$\\
							& 				       & $\dfeas$ & $1.0\ee{-12}$ & $3.7\ee{-15}$ & $3.3\ee{-5}$ & $0.0354$ & $5.0\ee{-10}$ & $2.4\ee{-15}$ & $4.6\ee{-15}$\\
							& 					   & $\gap$   & $1.2\ee{-8}$ & $1.3\ee{-10}$ & $0.0020$ & $0.9996$ & $5.2\ee{-6}$ & $1.0\ee{-9}$ & $4.0\ee{-10}$ \\
							& 					   & $\subopt$& $9.5\ee{-13}$ & $8.9\ee{-15}$ & $0.1974$ & $0.9990$ & $2.1\ee{-12}$ & $5.4\ee{-14}$ & $1.1\ee{-13}$ \\
							&   				   & time     & $65.8$ & $32.7$ & $90.2$ & $302.3$ & $48.8$ & $13.7$ & $4.0$ \\
\cline{2-10}
							& \multirow{5}{*}{$\#2$} & $\pfeas$ & $5.0\ee{-9}$ & $2.2\ee{-12}$ & $2.4\ee{-6}$ & $0.8372$ & $3.4\ee{-11}$ & $2.6\ee{-15}$ & $2.6\ee{-15}$ \\
							& 				       & $\dfeas$ & $1.3\ee{-12}$ & $8.0\ee{-11}$ & $3.4\ee{-5}$ & $0.0346$ & $3.5\ee{-11}$ & $3.6\ee{-15}$ & $5.6\ee{-13}$ \\
							& 					   & $\gap$   & $1.8\ee{-7}$ & $8.6\ee{-10}$ & $0.0021$ & $0.9996$ & $3.8\ee{-7}$ & $8.7\ee{-12}$ & $5.9\ee{-9}$ \\
							& 					   & $\subopt$& $3.6\ee{-12}$ & $1.9\ee{-13}$ & $0.0887$ & $0.9996$ & $9.6\ee{-14}$ & $9.6\ee{-14}$ & $1.3\ee{-13}$  \\
							&   				   & time     & $67.4$ & $33.2$ & $87.2$ & $303.8$ & $39.8$ & $13.3$ & $4.2$ \\
\cline{2-10}
							& \multirow{5}{*}{$\#3$} & $\pfeas$ & $1.3\ee{-9}$ & $3.0\ee{-12}$ & $2.2\ee{-6}$ & $0.8352$ & $4.5\ee{-13}$ & $1.4\ee{-16}$ & $0.0$ \\
							& 				       & $\dfeas$ & $2.2\ee{-12}$ & $3.5\ee{-15}$ & $3.2\ee{-5}$ & $0.0348$ & $1.7\ee{-10}$ & $4.0\ee{-15}$ & $3.3\ee{-13}$ \\
							& 					   & $\gap$   & $3.4\ee{-8}$ & $1.3\ee{-9}$ & $0.0019$ & $0.9996$ & $1.7\ee{-6}$ & $3.9\ee{-9}$ & $3.3\ee{-9}$ \\
							& 					   & $\subopt$& $4.0\ee{-15}$ & $5.8\ee{-13}$ & $0.0717$ & $0.9996$ & $3.9\ee{-14}$ & $3.8\ee{-12}$ & $3.4\ee{-14}$ \\
							&   				   & time     & $64.6$ & $34.7$ & $86.7$ & $301.2$ & $44.7$ & $13.5$ & $3.9$ \\
\hline
\hline
\multirow{15}{*}{\Large$ \substack{N:\ 100 \\ \\ n:\ 404 \\ \\ m:\ 31,301}$ } & \multirow{5}{*}{$\#1$} 
												   & $\pfeas$ & \multirow{5}{*}{**} & $1.3\ee{-11}$ & $2.4\ee{-7}$ & $0.8761$ & $7.1\ee{-15}$ & $5.3\ee{-15}$ & $5.3\ee{-15}$ \\
							& 				       & $\dfeas$ & & $4.5\ee-{15}$ & $2.6\ee{-5}$ & $0.0513$ & $9.6\ee{-6}$ & $1.6\ee{-13}$ & $5.4\ee{-15}$ \\
							& 					   & $\gap$   & & $4.1\ee{-11}$ & $1.5\ee{-4}$ & $0.9999$ & $0.0750$ & $2.3\ee{-9}$ & $4.7\ee{-9}$ \\
							& 					   & $\subopt$& & $7.0\ee{-14}$ & $0.1338$ & $0.9996$ &$0.0219$ & $1.9\ee{-11}$ & $1.9\ee{-15}$ \\
							&   				   & time     & & $974.9$ & $326.8$ & $830.9$ & $526.0$ & $50.7$ & $32.7$ \\
\cline{2-10}
							& \multirow{5}{*}{$\#2$} & $\pfeas$ & \multirow{5}{*}{**} & $3.8\ee{-12}$ & $3.6\ee{-7}$ & $0.8622$ & $1.1\ee{-14}$ & $0.0$ & $4.5\ee{-13}$ \\
							& 				       & $\dfeas$ & & $5.6\ee{-15}$  & $2.5\ee{-5}$ & $0.0518$ & $1.0\ee{-5}$ & $5.5\ee{-15}$ & $4.9\ee{-15}$ \\
							& 					   & $\gap$   & & $5.7\ee{-11}$ & $2.3\ee{-4}$ & $0.9999$ & $0.0797$ & $8.4\ee{-9}$ & $8.3\ee{-10}$ \\
							& 					   & $\subopt$& & $7.8\ee{-14}$ & $0.3539$ & $0.9999$ & $0.0212$ & $5.5\ee{-14}$ & $1.3\ee{-13}$ \\
							&   				   & time     & & $849.6$ & $321.4$ & $831.7$ & $488.2$ & $46.4$ & $17.7$ \\
\cline{2-10}
							& \multirow{5}{*}{$\#3$} & $\pfeas$ & \multirow{5}{*}{**} & $6.7\ee{-12}$ & $2.2\ee{-7}$ & $0.8570$ & $2.2\ee{-14}$ & $1.5\ee{-15}$ & $1.7\ee{-15}$ \\
							& 				       & $\dfeas$ & & $4.9\ee{-15}$ & $2.5\ee{-5}$ & $0.0514$ & $1.1\ee{-5}$ & $2.5\ee{-13}$ & $5.9\ee{-15}$ \\
							& 					   & $\gap$   & & $3.9\ee{-9}$ & $1.4\ee{-4}$ & $0.9999$ & $0.0882$ & $3.7\ee{-9}$ & $2.4\ee{-9}$ \\
							& 					   & $\subopt$& & $1.3\ee{-12}$ & $0.2188$ & $0.9997$ & $0.0186$ & $1.4\ee{-13}$ & $1.4\ee{-13}$ \\
							&   				   & time     & & $880.7$ & $323.9$ & $830.5$ & $454.1$ & $44.0$ & $16.1$ \\
\hline
\hline
\multirow{15}{*}{\Large$ \substack{N:\ 200 \\ \\ n:\ 804 \\ \\ m:\ 122,601}$ } & \multirow{5}{*}{$\#1$} 
							                       & $\pfeas$ & \multirow{5}{*}{**} & \multirow{5}{*}{**} & $1.0\ee{-6}$ & $0.8931$ & $8.9\ee{-15}$ & $0.0$ & $9.1\ee{-15}$  \\
							&					   & $\dfeas$ & & & $1.6\ee{-5}$ & $0.0904$ & $1.2\ee{-5}$ & $2.6\ee{-13}$ & $2.9\ee{-13}$ \\
							& 					   & $\gap$   & & & $4.8\ee{-4}$ & $1.0$ & $0.1086$ & $5.3\ee{-9}$ & $1.8\ee{-9}$ \\
							& 					   & $\subopt$& & & $0.3140$ & $0.9999$ & $0.0029$ & $3.9\ee{-11}$ & $4.3\ee{-9}$ \\
							&   				   & time     & & & $1208$& $1723$ & $1941$ & $291.9$ & $174.6$ \\
\cline{2-10}
							& \multirow{5}{*}{$\#2$} 
												   & $\pfeas$ & \multirow{5}{*}{**} & \multirow{5}{*}{**} & $1.0\ee{-6}$ & $0.8913$ & $4.5\ee{-8}$ & $5.9\ee{-15}$ & $4.5\ee{-15}$ \\
							& 				       & $\dfeas$ & & & $1.8\ee{-5}$ & $0.0919$ & $1.3\ee{-5}$ & $2.1\ee{-13}$ & $2.6\ee{-13}$ \\
							& 					   & $\gap$   & & & $4.8\ee{-4}$ & $1.0$ & $0.0671$ & $1.2\ee{-9}$ & $5.2\ee{-9}$ \\
							& 					   & $\subopt$& & & $0.2921$ & $0.9999$ & $0.0136$ & $5.5\ee{-9}$ & $4.7\ee{-13}$ \\
							&   				   & time     & & & $1206$ & $1646$ & $2349$ & $311.5$ & $185.2$ \\
\cline{2-10}
							& \multirow{5}{*}{$\#3$} 
												   & $\pfeas$ & \multirow{5}{*}{**} & \multirow{5}{*}{**} & $9.7\ee{-7}$ & $0.8892$ & $8.8\ee{-6}$ & $2.6\ee{-13}$ & $2.9\ee{-13}$ \\
							& 				       & $\dfeas$ & & & $1.8\ee{-5}$ & $0.0821$ & $1.3\ee{-5}$ & $8.2\ee{-15}$ & $1.2\ee{-14}$ \\
							& 					   & $\gap$   & & & $4.6\ee{-4}$ & $1.0$ & $0.0693$ & $1.7\ee{-9}$ & $1.5\ee{-9}$ \\
							& 					   & $\subopt$& & & $0.3027$ & $0.9999$ & $0.0398$ & $4.5\ee{-13}$ & $1.8\ee{-11}$ \\
							&   				   & time     & & & $1204$ & $1712$ & $2392$ & $272.6$ & $171.7$ \\
\hline
\hline
\multirow{15}{*}{\Large$ \substack{N:\ 500 \\ \\ n:\ 2004 \\ \\ m:\ 756,501}$ } & \multirow{5}{*}{$\#1$} 
											       & $\pfeas$ &\multirow{5}{*}{**} & \multirow{5}{*}{**} & $4.1\ee{-7}$ & $0.8985$ & $3.1\ee{-14}$ & $4.3\ee{-14}$ & $3.2\ee{-13}$ \\
							& 				       & $\dfeas$ & & & $1.3\ee{-5}$ & $0.1318$ & $1.1\ee{-5}$ & $1.5\ee{-14}$ & $1.4\ee{-13}$ \\
							& 					   & $\gap$   & & & $1.7\ee{-4}$ & $1.0$ & $0.0662$ & $3.8\ee{-11}$ & $4.3\ee{-9}$ \\
							& 					   & $\subopt$& & & $0.3385$ & $1.0$ & $0.2480$ & $2.3\ee{-13}$ & $1.7\ee{-12}$ \\
							&   				   & time     & & & $7659$ & $6035$ & $10001$ & $4389$ & $1819$ \\
\cline{2-10}
							& \multirow{5}{*}{$\#2$}& $\pfeas$&\multirow{5}{*}{**} &\multirow{5}{*}{**} & $1.8\ee{-7}$ & $0.8984$ & $2.6\ee{-4}$ & $7.7\ee{-13}$ & $1.1\ee{-12}$ \\
							& 				       & $\dfeas$ & & & $1.5\ee{-5}$ & $0.1295$ & $1.1\ee{-5}$ & $2.7\ee{-14}$ & $5.0\ee{-14}$ \\
							& 					   & $\gap$   & & & $7.4\ee{-5}$ & $1.0$ & $0.0780$ & $7.8\ee{-10}$ & $1.5\ee{-9}$ \\
							& 					   & $\subopt$& & & $0.3813$ & $1.0$ & $0.1154$ & $2.6\ee{-14}$ & $2.1\ee{-13}$ \\
							&   				   & time     & & & $7718$ & $6023$ & $10001$ & $2848$ & $1335$ \\
\cline{2-10}
							& \multirow{5}{*}{$\#3$}& $\pfeas$&\multirow{5}{*}{**} &\multirow{5}{*}{**} & $5.7\ee{-7}$ & $0.8984$ & $2.5\ee{-4}$ & $2.6\ee{-14}$ & $5.5\ee{-14}$ \\
							& 				       & $\dfeas$ & & & $1.3\ee{-5}$ & $0.1142$ & $1.1\ee{-5}$ & $1.4\ee{-14}$ & $1.7\ee{-14}$ \\
							& 					   & $\gap$   & & & $2.3\ee{-4}$ & $1.0$ & $0.0745$ & $1.1\ee{-10}$ & $1.2\ee{-11}$ \\
							& 					   & $\subopt$& & & $0.3326$ & $1.0$ & $0.0644$ & $4.4\ee{-15}$ & $1.0\ee{-13}$ \\
							&   				   & time     & & & $7849$ & $5994$ & $10219$ & $3743$ & $2316$ \\
\hline
%%%%%%%%%%%%%%%%%%%%%%%%%%%%%%%%%%%%%% Continued %%%%%%%%%%%%%%%%%%%%%%%%%%%%%%%%%%%%%%%%%%%
\hline
\multirow{15}{*}{\Large$ \substack{N:\ 1000 \\ \\ n:\ 4004 \\ \\ m:\ 3,013,001}$ } & \multirow{5}{*}{$\#1$} & $\pfeas$ & & & $1.5\ee{-6}$ & $0.8993$ & $4.0\ee{-14}$ & $1.0\ee{-13}$ & $1.4\ee{-14}$ \\
							& 				       & $\dfeas$ & \multirow{5}{*}{**} & \multirow{5}{*}{**} & $1.7\ee{-5}$ & $0.1747$ & $8.4\ee{-6}$ & $3.1\ee{-14}$ & $3.9\ee{-14}$ \\
							& 					   & $\gap$   & & & $4.0\ee{-3}$ & $1.0$ & $6.5\ee{-3}$ & $1.3\ee{-10}$ & $3.3\ee{-13}$ \\
							& 					   & $\subopt$& & & $0.5787$ & $1.0$ & $0.3102$ & $1.2\ee{-13}$ & $1.2\ee{-9}$ \\
							&   				   & time     & & & $62900$ & $42546$ & $50004$ & $30269$ & $26894$ \\
\cline{2-10}
							& \multirow{5}{*}{$\#2$} & $\pfeas$ & \multirow{5}{*}{**} & \multirow{5}{*}{**} & $1.5\ee{-6}$ & $0.8994$ & $3.7\ee{-8}$ & $8.0\ee{-13}$ & $1.1\ee{-16}$ \\
							& 				       & $\dfeas$ & & & $1.8\ee{-5}$ & $0.1316$ & $8.4\ee{-6}$ & $3.1\ee{-14}$ & $3.3\ee{-14}$ \\
							& 					   & $\gap$   & & & $4.0\ee{-3}$ & $1.0$ & $1.2\ee{-3}$ & $1.4\ee{-10}$ & $1.2\ee{-9}$ \\
							& 					   & $\subopt$& & & $0.5907$ & $1.0$ & $0.3118$ & $2.2\ee{-13}$ & $9.8\ee{-13}$ \\
							&   				   & time     & & & $60554$ & $42504$ & $44929$ & $50423$ & $41647$ \\
\cline{2-10}
							& \multirow{5}{*}{$\#3$} & $\pfeas$ & \multirow{5}{*}{**} & \multirow{5}{*}{**} & $1.5\ee{-6}$ & $0.8994$ & $1.3\ee{-13}$ & $2.4\ee{-13}$ & $3.4\ee{-15}$ \\
							& 				       & $\dfeas$ & & & $1.8\ee{-5}$ & $0.1212$ & $8.7\ee{-6}$ & $2.6\ee{-14}$ & $3.6\ee{-14}$ \\
							& 					   & $\gap$   & & & $4.0\ee{-3}$ & $1.0$ & $0.0174$ & $2.0\ee{-11}$ & $6.2\ee{-10}$ \\
							& 					   & $\subopt$& & & $0.5846$ & $1.0$ & $0.3265$ & $7.3\ee{-14}$ & $2.5\ee{-13}$ \\
							&   				   & time     & & & $64392$ & $42498$ & $50003$ & $48846$ & $45100$ \\
\hline
\end{tabular}
}%
\end{table}

{\bf Outlier-Robust Wahba Problem on Real Data}: 
To show the practical usefulness of {\name}, we test it on two applications of the Wahba problem on real data. The first application is image stitching on PASSTA~\cite{Meneghetti15SCIA-stitching} shown in Fig.~\ref{fig:wahba-result}(a). Given two images taken by the same camera with an unknown relative rotation, we first use SURF~\cite{Bay06eccv-surf} to establish $N=70$ putative keypoint matches, and then use {\name} to solve the SDP relaxation of the Wahba problem ($n=284,m=15,611$) to estimate the relative rotation and stitch the two images. {\name} obtains the globally optimal solution ($\max\{ \pfeas,\dfeas,\gap \}=6.2\ee{-13}$, $\subopt=7.8\ee{-14}$) in {$79$ seconds}. The second application is point cloud registration on 3DMatch~\cite{Zeng17cvpr-3dmatch} shown in Fig.~\ref{fig:wahba-result}(b). Given two point clouds with an unknown relative rotation, we first use FPFH~\cite{Rusu09icra-fast3Dkeypoints} and ROBIN~\cite{Shi21icra-robin} to establish $N=108$ keypoint matches, and then use {\name} to solve the SDP relaxation ($n=436,m=36,397$) to estimate the rotation and register the point clouds. {\name} obtains the globally optimal solution ($\max\{ \pfeas,\dfeas,\gap \}=1.6\ee{-10}$, $\subopt=1.8\ee{-13}$) in {$53$ seconds}.

%!TEX root = main.tex

\newcommand{\mpwtwo}{0.5\textwidth}
\begin{figure}[h]
\vspace{-2mm}
	\begin{center}
	\begin{minipage}{\textwidth}
	\begin{tabular}{cc}%
		%%%%%%%%%%%%%%%%%%%%%%%%%%%%%%%%%%%%%%%%%%%%%%%%%%%%%%%%%%%%%%%%%%%%%%%%%%%%%%%%%%%%%%%%%%%%%%%%%%%%%%%%%
			\hspace{-3mm}
			\begin{minipage}{\mpwtwo}%
			\centering%
			\includegraphics[width=\columnwidth]{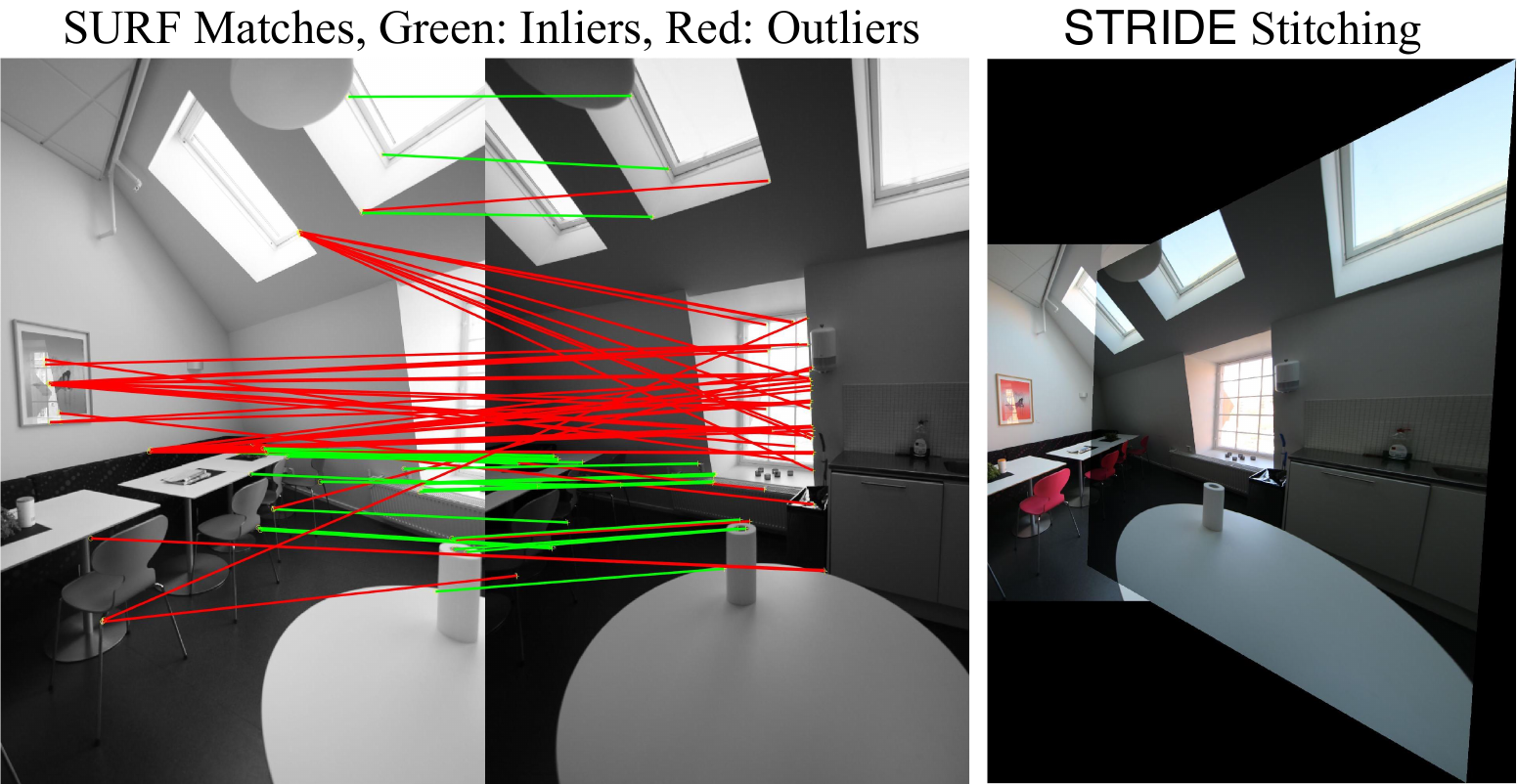} \\
			{\small (a) Image stitching on PASSTA \cite{Meneghetti15SCIA-stitching}.}
			\end{minipage}
		&  \hspace{-4mm}
			\begin{minipage}{\mpwtwo}%
			\centering%
			\includegraphics[width=\columnwidth]{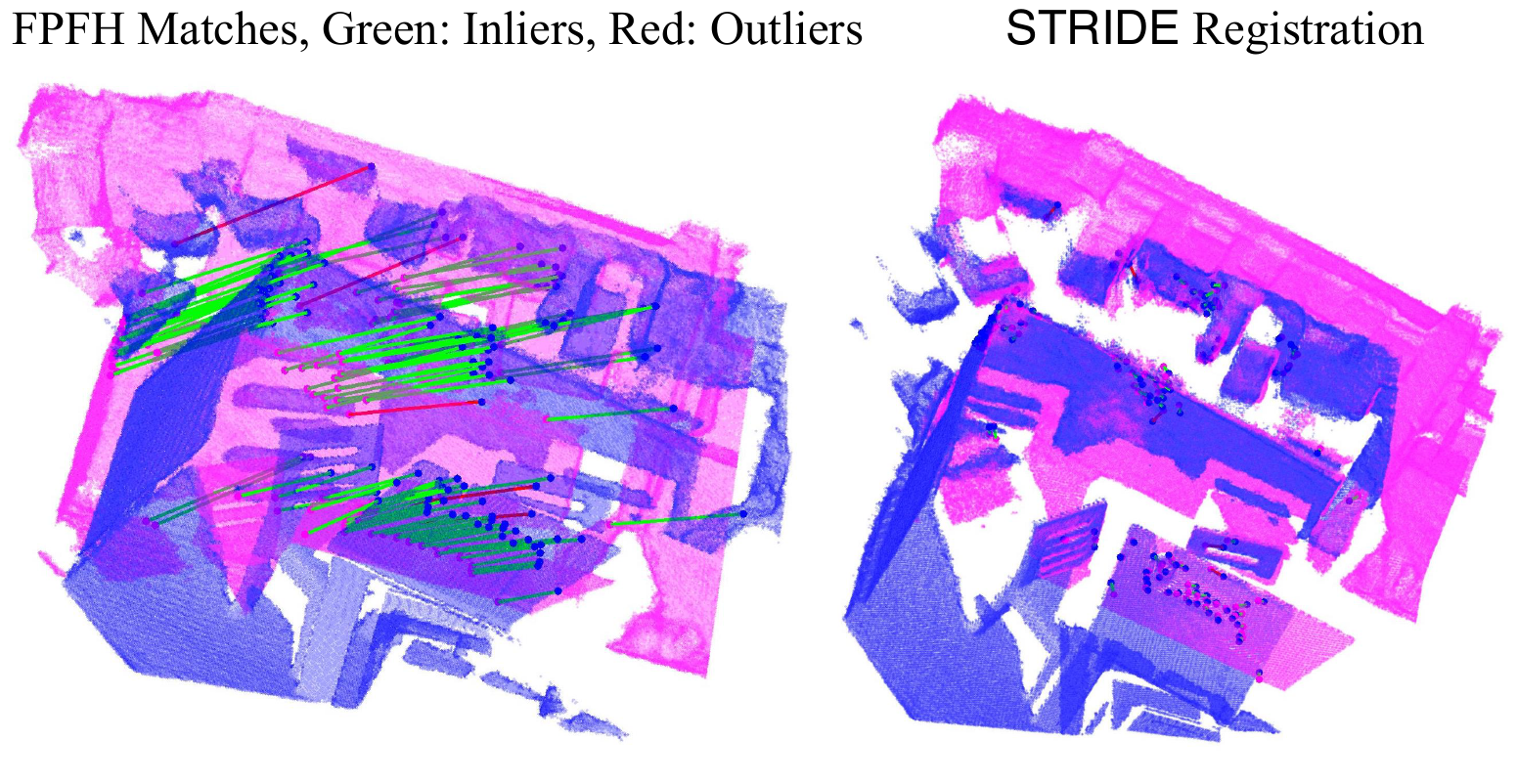} \\
			{\small (b) Point cloud registration on 3DMatch \cite{Zeng17cvpr-3dmatch}.}
			\end{minipage}
	\end{tabular}
	\vspace{-1mm}
	\end{minipage} 
	\caption{{\name} solves outlier-robust Wahba problems on real datasets to global optimality. 
	\label{fig:wahba-result}} 
	\end{center}
	\vspace{-4mm}
\end{figure}

%%%%%%%%%%%%%%%%%%%%%%%%%%%%%%%%%%%%%%%%%%%%%%%%%%%%%%%%%%%%%%%%%%%%%%%%%%%%%%%%%%%%%%%%%%%%%%%%%%%%%%%%%%%%%%%%%%%%%%%
%%%%%%%%%%%%%%%%%%%%%%%%%%%%%%%%%%%%%%%%%%%%%%%%%%%%%%%%%%%%%%%%%%%%%%%%%%%%%%%%%%%%%%%%%%%%%%%%%%%%%%%%%%%%%%%%%%%%%%%

\subsection{Nearest structured rank deficient matrices}
\label{sec:STLS}

Let $N,N_1,N_2$ be positive integers with $N_1 \leq N_2$, and let $\calL: \Real{N} \rightarrow \Real{N_1 \times N_2}$ be an affine map. Consider finding the nearest structured rank deficient matrix problem
\bea\label{eq:stlsraw}
\min_{u \in \Real{N}} \cbrace{ \norm{u - \theta}^2 \mid  \calL(u) \text{ is rank deficient} }, 
\eea
where $\theta \in \Real{N}$ is a given point. Problem~\eqref{eq:stlsraw} is commonly known as the \emph{structured total least squares} (STLS) problem~\cite{Rosen96SIMAA-stls,Markovsky07SP-stls}, and has numerous applications in control, systems theory, statistics~\cite{Markovsky08automatica-stls}, approximate greatest common divisor~\cite{Kaltofen06-approximateGCD}, camera triangulation~\cite{Aholt12eccv-qcqptriangulation}, among others~\cite{Markovsky14JCAM-slra}. Problem~\eqref{eq:stlsraw} can be reformulated as the following polynomial optimization problem
\begin{equation}\label{eq:stls}
\min_{z \in \usphere{N_1 - 1}, u \in \Real{N}} \cbrace{ \norm{u - \theta}^2 \ \middle\vert\ 
\blue{z\tran } \parentheses{L_0 + \sum_{i=1}^N u_i L_i} = 0 }, \tag{STLS}
\end{equation}
where $L_i \in \Real{N_1 \times N_2},i=0,\dots,N$, are the set of independent bases of the affine map $\calL$, $z \in \usphere{N_1 - 1}$ is a unit vector in the left kernel of $\calL(u)$ and acts as a witness of rank deficiency. Problem~\eqref{eq:stls} is easily seen to be nonconvex and the best algorithm in practice is based on local nonlinear programming~\cite{Markovsky14JCAM-slra}.

Recently, Cifuentes~\cite{Cifuentes19arXiv-stls} devised a semidefinite relaxation for~\eqref{eq:stls} and proved that the relaxation is guaranteed to be tight under a low-noise assumption~\cite{Cifuentes17arxiv-localstability}. This semidefinite relaxation is very similar to the relaxation for~\eqref{eq:wahba} presented in Section~\ref{sec:wahba} and is also a sparse second-order moment relaxation. Let $x = [z\tran,u\tran]\tran \in \Real{d},d=N_1 + N$ be the vector of unknowns in~\eqref{eq:stls}, construct the sparse monomial vector of degree up to 2 in $x$ as
\bea\label{eq:stlslift}
[x]_s = [z\tran,u_1 z\tran,\dots, u_N z\tran]\tran \in \Real{n}, \quad n= (N+1)N_1
\eea 
and build the moment matrix $X = [x]_s [x]_s\tran$. It can be easily checked that all the off-diagonal $N_1 \times N_1$ blocks, $u_i u_j z z\tran$, are symmetric by construction. Using $[x]_s$, the $N_2$ equality constraints in~\eqref{eq:stls} can be conveniently written as $[x]_s\tran a_i = 0,i=1,\dots,N_2$, for constant vectors $a_i \in \Real{n}$. In addition, each of the equality constraint also gives rise to $n$ redundant constraints of the form $([x]_s\tran a_i)[x]_s = 0$. Finally, the unit sphere constraint $z \in \usphere{N_1 -1}$ implies the trace of the leading $N_1 \times N_1$ block of $X$ is equal to 1. The construction above leads to a semidefinite relaxation with size
\bea
n = (N+1)N_1,\quad m = 1 + n N_2 + \trinum(N_1 - 1) \times \trinum(N).
\eea
Due to the limitation of interior point methods, Cifuentes~\cite{Cifuentes19arXiv-stls} was only able to numerically verify the tightness of the relaxation for very small problems (\eg$~N_1 \leq N_2 \leq 10, N < 20$).

We aim to compute globally optimal solutions of~\eqref{eq:stls} with much larger dimensions. We perform experiments on random instances of \eqref{eq:stls} where the affine map $\calL$ is structured to be a square Hankel matrix such that $N = N_1 + N_2 -1, N_1 = N_2$. We set $N_1 \in \{10,20,30,40\}$, and at each level we randomly generate three problem instances by drawing $\theta \sim \calN(0,I_N)$ from the standard Gaussian distribution. We then solve the sparse semidefinite relaxation using \sdpt, \mosek, \cdcs, \sketchy, \sdpnal, and \name. For \name, since the nonlinear programing~\eqref{eq:stls} does not admit any manifold structure, we use \fmincon~with an interior point method as the \nlp~method. To generate hypotheses for \nlp~from the moment matrix $X$, we follow
\bea
X = \sum_{i=1}^n \lambda_i v_i v_i\tran, \quad z_i = \frac{v_i[z]}{\norm{v_i[z]}}, \quad u^i_j = z_i\tran v_i[u_j z],\; 
j=1,\dots,N,
\eea
where we first perform a spectral decomposition, then round $z_i$ by projecting the entries of $v_i$ corresponding to block $z$ onto the unit sphere, and  round $u^i_j$ by computing the inner product between $z_i$ and the entries of $v_i$ corresponding to block $u_j z$ (again, the rationale for this rounding comes from the lifting \eqref{eq:stlslift}). We generate $r=5$ hypotheses by rounding $5$ eigenvectors. In order to set $M_b$ for computing $\subopt$, we make the assumption that the search variable $u$ of \eqref{eq:stls} contains random vectors that follow $\calN(0,I_N)$, and its squared norm follows a chi-square distribution of degree $N$. As a result, we choose the quantile corresponding to a $99.9\%$ probability, denoted as $\overline{M}$, as the bound on $\norm{u}^2$, such that $M_b = \overline{M}+1$ can upper bound the trace of $X$. 

Table \ref{table:stls} gives the numerical results of different solvers. Notice that at $N_1 = 40$, we have increased the maximum runtime of {\cdcs} and {\sdpnal} to $20000$ seconds to make a fair comparison with {\name}. Generally, Table \ref{table:stls} has similar results as Table \ref{table:wahba}. (i) IPMs can solve small problems to high accuracy, but cannot handle large problems due to memory issues. (ii) First-order solvers ({\cdcs} and {\sketchy}) are memory efficient but exhibit slow convergence and cannot attain high accuracy to certify global optimality and tightness of the relaxation. {\sdpnal} can solve problems with $N_1=10$ and $20$ to high accuracy, but failed for $N_1=30$ and $40$. The degraded performance of {\cdcs}, {\sketchy}, and {\sdpnal} compared to Tables \ref{table:pop:binary}-\ref{table:pop:sphere} again suggests the difficulty in solving sparse relaxations compared to dense relaxations. (iii) {\name} computed solutions of high accuracy for all test instances and it is about $5$ times faster than {\sdpnal}. 

%!TEX root = main.tex
\begin{table}
% \vspace{-7mm}
\caption{Results on solving sparse second-order relaxation of random~\eqref{eq:stls} instances. ``$**$'' indicates solver out of memory.}
\label{table:stls} 
\vspace{-1mm}
\adjustbox{max width=\textwidth}{%
\centering
\begin{tabular}{|c|c|c|ccccc||c|}
\hline
 Dimension & Run & Metric & \sdpt~\cite{Toh99OMS-sdpt3} & \mosek~\cite{mosek} & \cdcs~\cite{Zheng20MP-CDCS} & \sketchy~\cite{Yurtsever21SIMDS-scalableSDP} & \sdpnal~\cite{Yang15MPC-sdpnalplus} & \name \\
 \hline
 \hline
\multirow{15}{*}{\Large$ \substack{N_1:\ 10 \\ \\ n:\ 200 \\ \\ m:\ 10,551}$ } & \multirow{5}{*}{$\#1$} 
												   & $\pfeas$ & $4.8\ee{-12}$  & $1.6\ee{-11}$ & $3.7\ee{-9}$ & $0.3602$ &$1.0\ee{-10}$ & $6.2\ee{-16}$ \\
							& 				       & $\dfeas$ & $1.2\ee{-10}$  & $2.2\ee{-15}$ & $1.5\ee{-5}$ & $4.3\ee{-14}$ & $1.2\ee{-11}$ & $6.3
							\ee{-15}$ \\
							& 					   & $\gap$   & $5.5\ee{-6}$  & $1.1\ee{-12}$ & $9.4\ee{-9}$ & $9.1\ee{-4}$ & $5.3\ee{-9}$ & $4.8\ee{-12}$ \\
							& 					   & $\subopt$& $5.7\ee{-8}$  & $8.7\ee{-14}$ & $0.0078$ & $0.8808$ & $1.2\ee{-8}$ & $1.1\ee{-11}$ \\
							&   				   & time     & $94.0$ & $48.2$ & $138.6$ & $339.7$ & $13.3$ & $14.1$ \\
\cline{2-9}
							& \multirow{5}{*}{$\#2$} & $\pfeas$ & $1.2\ee{-10}$  & $1.8\ee{-10}$ & $4.5\ee{-6}$ & $0.3787$ & $1.7\ee{-11}$ & $1.5\ee{-15}$ \\
							& 				       & $\dfeas$ & $1.7\ee{-10}$ & $4.1\ee{-11}$ & $3.8\ee{-4}$ & $2.6\ee{-5}$ & $1.0\ee{-10}$ & $2.5\ee{-14}$ \\
							& 					   & $\gap$   & $5.5\ee{-5}$  & $5.2\ee{-9}$ & $1.9\ee{-4}$ & $0.2894$ & $6.3\ee{-10}$ & $2.1\ee{-11}$ \\
							& 					   & $\subopt$& $6.4\ee{-6}$  & $2.4\ee{-9}$ & $0.1179$ & $0.8544$ & $1.5\ee{-8}$ & $2.5\ee{-11}$ \\
							&   				   & time     & $92.2$ & $43.9$& $140.9$& $344.2$& $48.6$ & $18.4$ \\
\cline{2-9}
							& \multirow{5}{*}{$\#3$} & $\pfeas$ & $1.6\ee{-10}$ & $9.0\ee{-12}$ & $8.3\ee{-7}$ & $0.3943$ & $7.5\ee{-11}$ & $8.2\ee{-16}$\\
							& 				       & $\dfeas$ & $2.3\ee{-11}$  & $9.7\ee{-9}$ & $1.4\ee{-4}$ & $1.6\ee{-5}$ & $8.9\ee{-11}$ & $1.0\ee{-13}$ \\
							& 					   & $\gap$   & $1.9\ee{-6}$ & $1.7\ee{-10}$ & $5.3\ee{-5}$ & $0.0164$ & $7.3\ee{-11}$ & $8.0\ee{-11}$ \\
							& 					   & $\subopt$& $2.4\ee{-8}$ & $1.5\ee{-11}$ & $0.1406$ & $0.8993$ & $3.9\ee{-13}$ & $1.9\ee{-10}$\\
							&   				   & time     & $99.2$  & $44.4$ & $139.8$ & $338.1$ & $38.9$ & $17.4$ \\
\hline
\hline
\multirow{15}{*}{\Large$ \substack{N_1:\ 20 \\ \\ n:\ 800 \\ \\ m:\ 164,201}$ } & \multirow{5}{*}{$\#1$} 
							                       & $\pfeas$ & \multirow{5}{*}{**} & \multirow{5}{*}{**} & $2.7\ee{-6}$ & $0.9608$ & $3.2\ee{-13}$ & $1.9\ee{-15}$ \\
							&					   & $\dfeas$ & & & $3.4\ee{-4}$ & $0.0325$ & $7.8\ee{-10}$ & $3.8\ee{-14}$ \\
							& 					   & $\gap$   & & & $2.2\ee{-4}$ & $0.1207$ & $2.8\ee{-7}$ & $1.8\ee{-10}$ \\
							& 					   & $\subopt$& & & $0.3945$ & $0.9863$ & $1.2\ee{-8}$ & $2.8\ee{-10}$ \\
							&   				   & time     & & & $4206$ & $2171$ & $1189$ & $292.5$ \\
\cline{2-9}
							& \multirow{5}{*}{$\#2$} 
												   & $\pfeas$ & \multirow{5}{*}{**} & \multirow{5}{*}{**} & $2.5\ee{-6}$ & $0.9326$ & $5.6\ee{-13}$ & $1.5\ee{-15}$ \\
							& 				       & $\dfeas$ & & & $2.7\ee{-4}$ & $0.0109$ & $1.7\ee{-10}$ & $2.4\ee{-15}$ \\
							& 					   & $\gap$   & & & $1.6\ee{-4}$ & $0.4961$ & $2.0\ee{-10}$ & $2.8\ee{-10}$ \\
							& 					   & $\subopt$& & & $0.4805$ & $0.9714$ & $9.1\ee{-8}$ & $7.2\ee{-10}$ \\
							&   				   & time     & & & $4187$ & $2163$ & $870.6$ & $257.7$ \\
\cline{2-9}
							& \multirow{5}{*}{$\#3$} 
												   & $\pfeas$ & \multirow{5}{*}{**} & \multirow{5}{*}{**} & $4.3\ee{-6}$ & $0.9450$ & $5.3\ee{-11}$ & $1.6\ee{-15}$\\
							& 				       & $\dfeas$ & & & $3.4\ee{-4}$ & $0.0200$ & $7.9\ee{-11}$ & $2.2\ee{-15}$ \\
							& 					   & $\gap$   & & & $3.1\ee{-4}$ & $0.1536$ & $1.3\ee{-8}$ & $1.2\ee{-10}$ \\
							& 					   & $\subopt$& & & $0.4493$ & $0.9798$ & $4.1\ee{-9}$ & $2.5\ee{-10}$ \\
							&   				   & time     & & & $4224$ & $2178$ & $1169$ & $324.4$\\
\hline
\hline
\multirow{15}{*}{\Large$ \substack{N_1:\ 30 \\ \\ n:\ 1800 \\ \\ m:\ 823,951}$ } & \multirow{5}{*}{$\#1$} 
											       & $\pfeas$ &\multirow{5}{*}{**} & \multirow{5}{*}{**} & $1.4\ee{-5}$ & $1.3862$ & $1.7\ee{-12}$ & $5.4\ee{-15}$\\
							& 				       & $\dfeas$ & & & $1.0\ee{-3}$ & $1.3752$ & $2.4\ee{-5}$ & $4.9\ee{-14}$ \\
							& 					   & $\gap$   & & & $1.7\ee{-4}$ & $0.5175$ & $0.1060$ & $1.6\ee{-9}$ \\
							& 					   & $\subopt$& & & $0.87$ & $0.9997$ & $0.2083$ & $6.8\ee{-10}$ \\
							&   				   & time     & & & $12234$ & $6047$ & $10000$ & $1679$ \\
\cline{2-9}
							& \multirow{5}{*}{$\#2$}& $\pfeas$&\multirow{5}{*}{**} &\multirow{5}{*}{**} & $1.6\ee{-5}$ & $1.3826$ & $2.0\ee{-14}$ & $1.6\ee{-8}$\\
							& 				       & $\dfeas$ & & & $0.0011$ &$1.3544$ & $1.2\ee{-5}$ & $2.3\ee{-11}$ \\
							& 					   & $\gap$   & & & $3.3\ee{-4}$ & $0.4171$ & $0.1077$ & $2.8\ee{-7}$ \\
							& 					   & $\subopt$& & & $0.8248$ & $0.9996$ & $0.0454$ & $2.7\ee{-7}$ \\
							&   				   & time     & & & $12324$ & $5994$ & $10001$ & $1813$ \\
\cline{2-9}
							& \multirow{5}{*}{$\#3$}& $\pfeas$&\multirow{5}{*}{**} &\multirow{5}{*}{**} & $1.4\ee{-5}$ & $1.3819$ & $2.2\ee{-12}$ & $4.2\ee{-15}$ \\
							& 				       & $\dfeas$ & & & $0.0010$ & $1.3376$ & $5.8\ee{-6}$ & $2.3\ee{-15}$ \\
							& 					   & $\gap$   & & & $3.7\ee{-4}$ & $0.3730$ & $0.0569$ & $6.9\ee{-10}$ \\
							& 					   & $\subopt$& & & $0.8191$ & $0.9996$ & $0.0640$ & $4.6\ee{-10}$ \\
							&   				   & time     & & & $12360$ & $5976$ & $10000$ & $2050$ \\
\hline
\hline
\multirow{15}{*}{\Large$ \substack{N_1:\ 40 \\ \\ n:\ 3200 \\ \\ m:\ 2,592,801}$ } & \multirow{5}{*}{$\#1$} 
												   & $\pfeas$ &\multirow{5}{*}{**} &\multirow{5}{*}{**} & \multirow{5}{*}{**} & $1.6887$ & $9.7\ee{-14}$  & $5.6\ee{-15}$ \\
							& 				       & $\dfeas$ & & & & $2.8283$ & $2.3\ee{-5}$ & $2.2\ee{-12}$ \\
							& 					   & $\gap$   & & & & $0.4648$ & $0.0439$ & $4.3\ee{-10}$ \\
							& 					   & $\subopt$& & & & $1.0$ & $0.2345$ & $1.3\ee{-7}$ \\
							&   				   & time     & & & & $17222$ & $20003$ & $14541$ \\
\cline{2-9}
							& \multirow{5}{*}{$\#2$}& $\pfeas$&\multirow{5}{*}{**} &\multirow{5}{*}{**} & \multirow{5}{*}{**} & $1.6865$ & $2.3\ee{-5}$ & $4.9\ee{-15}$ \\
							& 				       & $\dfeas$ & & & & $2.8528$ & $2.6\ee{-5}$ & $1.0\ee{-13}$ \\
							& 					   & $\gap$   & & & & $0.4391$ & $0.0750$ & $2.1\ee{-9}$ \\
							& 					   & $\subopt$& & & & $1.0$ & $0.3296$ & $3.5\ee{-9}$ \\
							&   				   & time     & & & & $17270$ & $20008$ & $9994$\\
\cline{2-9}
							& \multirow{5}{*}{$\#3$}& $\pfeas$&\multirow{5}{*}{**} &\multirow{5}{*}{**} & \multirow{5}{*}{**} & $1.6875$ & $3.5\ee{-4}$  & $7.9\ee{-15}$ \\
							& 				       & $\dfeas$ & & & & $2.9918$ & $1.4\ee{-4}$ & $5.3\ee{-12}$ \\
							& 					   & $\gap$   & & & & $0.5329$ & $0.1743$ & $2.0\ee{-8}$\\
							& 					   & $\subopt$& & & & $1.0$ & $0.6361$ & $2.8\ee{-7}$ \\
							&   				   & time     & & & & $17227$ & $20345$ & $10120$\\
\hline
\end{tabular}
}%
\end{table}

\section{Concluding Remarks}
\label{sec:conclusions}

In this paper, we have designed an efficient algorithmic framework, \name, to solve large-scale tight semidefinite programming (SDP) relaxations of polynomial optimization problems (POPs). {\name} employs an inexact projected gradient method (iPGM) as its backbone, but leverages fast nonlinear programming (NLP) methods to seek rapid primal acceleration. By conducting a novel convergence analysis for the iPGM and taking safeguarded steps, we have proved the global convergence of \name. For solving the projection subproblem in every iPGM step, we have designed a modified limited-memory BFGS method and proved its global convergence. In addition to the contributions in algorithmic design and convergence analysis, we have studied several important classes of POPs and their dense or sparse (tight) relaxations. In our extensive numerical experiments, {\name} globally solved all the POPs and the corresponding SDP relaxations to very high accuracy, and it offered state-of-the-art efficiency and robustness. We hope the practical performance of {\name} can encourage more theoretical study towards investigating when and why the NLP iterates can produce effective acceleration to the {\ipgm} backbone.

\begin{acknowledgements}
The authors would like to thank Diego Cifuentes for useful discussions about the nearest structured rank deficient matrix problem, and Jie Wang for clarifications about sparse semidefinite relaxations. H. Yang and L. Carlone were partially funded by ARL DCIST CRA W911NF-17-2-0181 and NSF CAREER award ``Certifiable Perception for Autonomous Cyber-Physical Systems''. 
% L. Liang and K.-C. Toh were partially funded by ......
\end{acknowledgements}

\appendix
%!TEX root = main.tex

\section{Proof of Theorem \ref{thm-iPGM}}
\label{appendix-proof-iPGM}

Let us first prove the following useful lemma for later convenience. 

\begin{lemma}
\label{lemma-fk}
	Under the inexactness conditions in \eqref{eq-iPGM-inexact-conditions}, for any $X\in \mathbb{S}_+^n$, it holds that
	\[
		\inprod{C}{X - X^k} \geq \frac{1}{\sigma_k}\norm{X^{k} - X^{k-1}}^2 + \frac{1}{\sigma_k}\inprod{X^{k-1} - X}{X^{k} - X^{k-1}} + \inprod{\mathcal{A}^*y^k}{X - X^k}.
	\]
\end{lemma}
\begin{proof}
	From the last two conditions in \eqref{eq-iPGM-inexact-conditions}, $\mathcal{A}^*y^k + S^k - C = \frac{1}{\sigma_k}\left(X^k - X^{k-1}\right)$ and $\inprod{X^k}{S^k} = 0$, we have that
	\[
		\begin{aligned}
			\inprod{C}{X - X^k} = &\; \inprod{\mathcal{A}^*y^k + S^k - \frac{1}{\sigma_k}\left(X^k - X^{k-1}\right)}{X - X^k} \\
			= &\; \frac{1}{\sigma_k} \inprod{X^k - X^{k-1}}{X^k - X} + \inprod{\mathcal{A}^*y^k}{X - X^k} + \underbrace{\inprod{S^k}{X}}_{\geq 0} - \underbrace{\inprod{S^k}{X^k}}_{ = 0} \\
			\geq &\; \frac{1}{\sigma_k}\norm{X^{k} - X^{k-1}}^2 + \frac{1}{\sigma_k}\inprod{X^{k-1} - X}{X^{k} - X^{k-1}} + \inprod{\mathcal{A}^*y^k}{X - X^k},
		\end{aligned}
	\]
	where the last inequality is due to the following fact:
	\bea
		\frac{1}{\sigma_k} \inprod{X^k - X^{k-1}}{X^k - X} = &  \displaystyle \frac{1}{\sigma_k} \inprod{X^k- X^{k-1}}{X^k -X^{k-1} + X^{k-1} - X} \nonumber \\
		 = & \displaystyle \frac{1}{\sigma_k}\norm{X^{k} - X^{k-1}}^2 + \frac{1}{\sigma_k}\inprod{X^{k-1} - X}{X^{k} - X^{k-1}}. \nonumber 
	\eea
	Thus, the proof is completed. \qed
\end{proof}

Then, we shall conduct our proof of Theorem \ref{thm-iPGM} as follows.

\begin{proof}
	For notational simplicity, denote $r^k = \mathcal{A}(X^k) - b$ and $e^k = X^k - \Xstar$ for all $k\geq 0$. Then, $\mathcal{A}(e^k) = r^k$ since $\mathcal{A}(\Xstar) = b$. Moreover, we have $\norm{r^k}\leq \varepsilon_k $ from the first condition in \eqref{eq-iPGM-inexact-conditions}. Therefore, the convergence of $\left\{\norm{r^k}\right\}$ is obvious, since $\norm{r^k}\leq \varepsilon_k$ and $\{k\varepsilon_k\}$ is summable. Recall that $\Xstar$ is an optimal solution to problem \eqref{eq:primalSDP}. Let $X = X^{k-1}$ in Lemma \ref{lemma-fk}, we get 
	\begin{equation}
	\label{eq-thm-iPGM-1}
		\inprod{C}{X^{k-1} - X^k} \geq \frac{1}{\sigma_k}\norm{e^k - e^{k-1}}^2 + \inprod{\mathcal{A}^*y^k}{e^{k-1} - e^{k}}	
	\end{equation}
	On the other hand, let $X = \Xstar$ in Lemma \ref{lemma-fk}, we obtain
	\begin{equation}
	\label{eq-thm-iPGM-2}
		\inprod{C}{\Xstar - X^k} \geq \frac{1}{\sigma_k}\norm{e^k - e^{k-1}}^2 + \frac{1}{\sigma_k}\inprod{e^{k-1}}{e^k - e^{k-1}} - \inprod{A^*y^k}{e^k}.	
	\end{equation}
	Multiplying $k-1$ to \eqref{eq-thm-iPGM-1} and add the resulting inequality to \eqref{eq-thm-iPGM-2} yields
	\begin{equation}
		\label{eq-thm-iPGM-3}
		\begin{aligned}
			&\; \inprod{C}{(k-1)\left(X^{k-1} - \Xstar \right) - k\left(X^k - \Xstar\right)  } \\ 
			\geq &\; \frac{k}{\sigma_k}\norm{e^k - e^{k-1}}^2 + \frac{1}{\sigma_k}\inprod{e^{k-1}}{e^k - e^{k-1}} + \inprod{\mathcal{A}^*y^k}{(k-1)e^{k-1} - ke^k} \\
			= &\; \frac{1}{\sigma_k}\left(k\norm{e^k}^2 + (k-1)\norm{e^{k-1}}^2 - (2k -1)\inprod{e^k}{e^{k-1}}\right) + \inprod{y^k}{(k-1)r^{k-1} - k r^k} \\
			\geq &\; \frac{1}{\sigma_k}\left(k\norm{e^k}^2 + (k-1)\norm{e^{k-1}}^2 - \frac{2k-1}{2}\left(\norm{e^k}^2 + \norm{e^{k-1}}^2\right)\right) + \inprod{y^k}{(k-1)r^{k-1} - k r^k} \\
			 = &\; \frac{1}{2\sigma_k}\norm{e^k}^2 - \frac{1}{2\sigma_k}\norm{e^{k-1}}^2 + \inprod{y^k}{(k-1)r^{k-1} - k r^k} \\
			 \geq &\; \frac{1}{2\sigma_k}\norm{e^k}^2 - \frac{1}{2\sigma_{k-1}}\norm{e^{k-1}}^2 + \inprod{y^k}{(k-1)r^{k-1} - k r^k},
		\end{aligned} 
	\end{equation}
	where in the last inequality, we use the fact that $\sigma_k \geq \sigma_{k-1} > 0$ for all $k\geq 1$. By summing the inequality in \eqref{eq-thm-iPGM-3} for $k=1,\dots,k$, we get (recall that $\norm{y^k}\leq M$ for all $k\geq 0$)
	\begin{equation*}
	\begin{aligned}
			-k \inprod{C}{X^k - \Xstar} \geq &\; \frac{1}{2\sigma_k}\norm{e^k}^2 - \frac{1}{2\sigma_0}\norm{e_0}^2 + \sum_{i=1}^k\inprod{y^i}{(i-1)r^{i-1} - i r^i} \\
			\geq &\; \frac{1}{2\sigma_k}\norm{e^k}^2 - \frac{1}{2\sigma_0}\norm{e_0}^2 - 2M\sum_{i=1}^k i\varepsilon_i.
		\end{aligned}	
	\end{equation*}
	The above inequality leads to the following upper bound on $\inprod{C}{X^k - \Xstar}$:
	\begin{equation}
	\label{eq-thm-iPGM-4}
	\inprod{C}{X^k - \Xstar} \leq \frac{1}{k}\left(\frac{1}{2\sigma_0}\norm{e^0}^2 + 2M\sum_{i=1}^k i\varepsilon_i\right).	
	\end{equation}
	
	Next, we shall provide a lower bound for $\inprod{C}{X^k - \Xstar}$ for $k\geq 1$. To this end, let us define the set $\mathcal{F}_k$, which is an enlargement of the primal feasible set $\mathcal{F}_{\mathrm{P}}$ with respect to $\varepsilon_k \geq 0$, as follows:
	\[
		\mathcal{F}_k:= \left\{X\in \mathbb{S}^n\; \middle | \; \norm{\mathcal{A}(X) - b}\leq \varepsilon_k,\; X\in \mathbb{S}_+^n\right\},\quad \forall\, k\geq 1. 
	\]
	Moreover, denote
	\[
		\Xstar_k:=\argmin\left\{\inprod{C}{X}\; \middle | \; X\in \mathcal{F}_k\right\},\quad \forall\, k\geq 1.
	\]
	Since $(\ystar, \Sstar)\in \mathbb{R}^m\times \mathbb{S}_+^n$ is a dual optimal solution, by \cite[Lemma 3.4]{Jiang2012}, we have that 
	\[
		0\leq \inprod{C}{\Xstar - \Xstar_k} \leq \norm{\ystar}\varepsilon_k,\quad \forall\, k\geq 1.
	\]
	As a consequence, because $X^k \in \calF_k$ and $\inprod{C}{X^k} \geq \inprod{C}{\Xstar_k}$, it holds that 
	\begin{equation}
	\label{eq-thm-iPGM-5}
		\inprod{C}{X^k - \Xstar} \geq \inprod{C}{\Xstar_k - \Xstar} \geq -\norm{\ystar}\varepsilon_k = -\frac{\norm{\ystar}k\varepsilon_k}{k}.	
	\end{equation}
	Combining \eqref{eq-thm-iPGM-4} and \eqref{eq-thm-iPGM-5}, we see that 
	\[
		-\frac{\norm{\ystar}k\varepsilon_k}{k} \leq \inprod{C}{X^k - \Xstar} \leq  \frac{1}{k}\left(\frac{1}{2\sigma_0}\norm{e^0}^2 + 2M\sum_{i=1}^k i\varepsilon_i\right).
	\]
	
	Finally, let us estimate the dual infeasibility. To this end, we get from \eqref{eq-thm-iPGM-1} that
	\[
		\frac{1}{\sigma_k}\norm{X^k - X^{k-1}}^2 \leq  \inprod{C}{X^{k-1} - \Xstar} - \inprod{C}{X^{k} - \Xstar} +  \norm{y^k}\left(\norm{r^{k-1}} + \norm{r^k}\right)
	\]
	which implies that 
	\[
		\begin{aligned}
			&\; \norm{\mathcal{A}^*y^k + S^k - C} =  \frac{1}{\sigma_k}\norm{X^{k} - X^{k-1}} \\
			\leq &\; \frac{1}{\sqrt{\sigma_k}}\sqrt{ \frac{1}{k-1}\left(\frac{1}{2\sigma_0}\norm{e^0}^2 + 2M\sum_{i=1}^{k-1} i\varepsilon_i\right) + \norm{\ystar}\varepsilon_k + M\left(\varepsilon_{k-1} + \varepsilon_k\right)} \\ 
			\leq &\; O\left(\frac{1}{\sqrt{k\sigma_k}}\right),
		\end{aligned} 
	\]
	and establishes the convergence of the dual infeasibility. Therefore, the proof is completed. \qed
\end{proof}
%!TEX root = main.tex

\section{Proof for Theorem~\ref{thm:convergence-lbfgs}}
\label{sec:proof:theorem-lbfgs}

From Step 1 of the modified L-BFGS algorithm, we have $d^k = (-\beta_k I+Q_k)\nabla \phi(\xi^k)$, with $Q_k = 0$ if $\norm{d^k} \geq K$. This leads to $\inprod{\nabla \phi(\xi^k)}{d^k}  = \inprod{\nabla \phi(\xi^k)}{-(\beta_k I+Q_k)\nabla \phi(\xi^k)}$, which implies that 
\begin{equation}
\label{eq-pf-lbfgs-1}
\lambda_{\min}(\beta_kI + Q_k) \leq \frac{\inprod{-\nabla \phi(\xi^k)}{d^k}}{\norm{\nabla \phi(\xi^k)}^2},	
\end{equation}
where $\lambda_{\min}(\cdot)$ denotes the minimum eigenvalue. This implies that $d^k$ is 
a descent direction and the line search scheme is well-defined.

Let the sequence $\{\xi^k\}$ be generated by the algorithm. We next show the second part of the theorem, \ie $\{\xi^k\}$ is a bounded sequence. Since $\cA\cA^*$ is nonsigular and the Slater condition holds, by \cite[Theorem 17' \& 18']{Rockafellar1974}, the level set $\cL_\phi:=\left\{\xi\in \R^m\,:\, \phi(\xi)\leq \phi(\xi^0)\right\}$ is compact. Thus, the boundedness of the sequence is proven. % Moreover, as a consequence of the boundedness, there also exists a positive constant $M_0$ such that $\norm{\nabla \phi (\xi^k)}\leq M_0$, for all $k$.

Let us now prove the third part of the theorem, which states that every accumulation point of $\{\xi^k\}$ is an optimal solution solution of problem \eqref{prob:projetion-dual-y}. To this end, without loss of generality, we assume that $\{\xi^k\}$ is an infinite sequence (the case for a finite sequence is trivial and we omit it here) and $\xi^k\rightarrow \hat{\xi}$ as $k\rightarrow \infty$, and we shall prove \blue{that} 
$\nabla \phi(\hat{\xi}) = 0$. We plan to prove the latter result by contraction. Let us assume that $\nabla \phi(\hat{\xi}) \neq 0$. Under this assumption, we will deduce a contradiction in three steps.
\begin{description}
	\item [Step A.] We first claim that when $\nabla \phi(\hat{\xi}) \neq 0$, $\norm{d^k}$ is bounded from below and above, \ie there exist two positive constants $M_1$ and $M_2$ such that $0 < M_1 \leq \norm{d^k} \leq M_2$ for all $k$. 
	
	On the one hand, by Step 1 of the L-BFGS Algorithm \ref{alg:lbfgs}, we know that either $d^k = -(\beta_k I+Q_k)\nabla\phi(\xi^k)$, in which case $\norm{d^k}\leq K$ (otherwise the algorithm rejects it), or $d^k = -\beta_k\nabla\phi(\xi^k)$, in which case $\norm{d^k}\leq \tau_1M_0^{1+\tau_2}$ (due to $\beta_k = \tau_1\norm{\nabla \phi(\xi^k)}^{\tau_2}$). Therefore, we have $\norm{d^k}\leq \max\left\{K, \tau_1M_0^{1+\tau_2}\right\}:=M_2 < \infty$, and $\norm{d^k}$ is upper bounded. 
	
	On the other hand, if for some index set $J$, $\norm{d^{k_j}} \rightarrow 0$ for $j\in J$, then we can deduce that 
	\[
		\norm{\nabla \phi(\xi^{k_j})} = \norm{(\beta_{k_j}I + Q_{k_j})^{-1}d^{k_j}} \leq \frac{1}{\beta_{k_j}}\norm{d^{k_j}} = \frac{1}{\tau_1\norm{\nabla\phi(\xi^{k_j}}}\norm{d^{k_j}},\quad j\in J.
	\]
	The above inequality implies that $\norm{\nabla \phi(\hat{\xi})} = \lim_{j\in J} \norm{\nabla \phi(\xi^{k_j})} = 0$, a contradiction. Thus, $\norm{d^k}$ must also be lower bounded, by some constant $M_1 > 0$. 
	\item [Step B.] We next claim that $\{\rho^{m_k}\}$ is bounded away from zero. 
	
	Suppose that, for some index set $J$, $\rho^{m_{k_j} }\rightarrow 0$, $j\in J$ (\ie $m_{k_j}\rightarrow \infty$ due to $\rho\in (0,1)$). Recall that $m_{k_j}$ is the smallest nonnegative integer to satisfy the inequality in the line search scheme (so that $m_{k_j}-1$ does not satisfy the inequality). \blue{Thus,} it holds that 
	\begin{equation}
		\label{eq-pf-lbfgs-2}
		\phi(\xi^{k_j}+ \rho^{m_{k_j}-1} d^{k_j}) > \phi(\xi^{k_j}) + \mu \rho^{m_{k_j}-1} \inprod{\nabla \phi(\xi^{k_j})}{d^{k_j}}. 
	\end{equation}
	Meanwhile, we can expand $\phi(\xi^{k_j} + \rho^{m_{k_j}-1}d^{k_j})$ as 
	\begin{equation}
	\label{eq-pf-lbfgs-3}
		\phi(\xi^{k_j} + \rho^{m_{k_j}-1}d^{k_j}) = \phi(\xi^{k_j}) + \int_0^1 
		\inprod{\nabla \phi(\xi^{k_j} + t\rho^{m_{k_j}-1}d^{k_j})}{\rho^{m_{k_j}-1}d^{k_j}} dt.	
	\end{equation}
	By substituting \eqref{eq-pf-lbfgs-3} into \eqref{eq-pf-lbfgs-2}, we obtain the following inequality:
	\begin{equation}
	\label{eq-pf-lbfgs-4}
		\int_0^1\inprod{\nabla\phi(\xi^{k_j} + t\rho^{m_{k_j}-1}d^{k_j})}{d^{k_j}}dt > \mu \inprod{\nabla \phi(\xi^{k_j})}{d^{k_j}},\quad j\in J.
	\end{equation}
	Since the sequence $\{d^{k_j}\}$ is bounded from both below and above, by taking a subset of $J$ if necessary, we assume that there exists $\hat{d}$ such that $d^{k_j}\rightarrow \hat{d}$, $j\in J$. Taking limit with respect to $j\in J$ in \eqref{eq-pf-lbfgs-4} yields that $\inprod{\nabla\phi(\hat{\xi})}{\hat{d}} \geq \mu \inprod{\nabla \phi(\hat{\xi})}{\hat{d}}$, which further implies that $\inprod{\nabla \phi(\hat{\xi})}{\hat{d}} \geq 0$ by the fact that $\mu\in (0, 1/2)$. Also, since $\inprod{\nabla \phi(\hat{\xi})}{\hat{d}}\leq 0$, we derive $\inprod{\nabla \phi(\hat{\xi})}{\hat{d}} = 0$. However, denoting $\hat{\beta}$ as the limit of $\beta_{k_j}$ as $k_j\rightarrow\infty$, $j\in J$, and \blue{noting that} $Q_k \succeq 0$ for all $k$, the inequality in \eqref{eq-pf-lbfgs-1} shows that 
		\[
			0 = -\inprod{\nabla \phi(\hat{\xi})}{\hat{d}} \geq \hat{\beta} \norm{\nabla \phi(\hat{\xi})}^2 \geq \tau_1 \norm{\nabla \phi(\hat{\xi})}^{2+\tau_2},
		\]
		which implies that $\nabla \phi(\hat{\xi}) = 0$. This is also a contradiction. Therefore, $\{\rho^{m_k}\}$ is bounded away from zero.
	\item[Step C.] Now, using the fact that $\phi$ is bounded from below on the bounded sequence $\{\xi^k\}$, we have that $\{\phi(\xi^{k+1}) - \phi(\xi^k)\}\rightarrow 0$ (since the sequence $\{\phi(\xi^k)\}$ is nonincreasing and $\phi(\cdot)$ is continuous). This implies, by the line search rule again, that $\rho^{m_k}\inprod{\nabla \phi(\xi^k)}{d^k}\rightarrow 0$, as $k\rightarrow \infty$. Therefore, it holds that $\inprod{\nabla \phi(\hat{\xi})}{\hat{d}} = 0$, since $\{\rho^{m_k}\}$ is bounded away from zero, which again implies that $\nabla \phi(\hat{\xi}) = 0$, and this contradicts our initial hypothesis that $\nabla \phi(\hat{\xi}) \neq 0$. 
\end{description}

At this point, we have shown that $\nabla \phi(\hat{\xi}) = 0$, as desired. At last, by the convexity of $\phi(\cdot)$, vanishing of the gradient at $\hat{\xi}$ implies that $\hat{\xi}$ is an optimal solution of problem \eqref{prob:projetion-dual-y}. Therefore, the proof is completed. \qed


%\clearpage
\bibliographystyle{spmpsci}      % mathematics and physical sciences

% \bibliographystyle{siam}
% \bibliographystyle{plain}
% \bibliographystyle{plain}
\bibliography{references}

% % Non-BibTeX users please use
% \begin{thebibliography}{}
% %
% % and use \bibitem to create references. Consult the Instructions
% % for authors for reference list style.
% % Format for Journal Reference
% Author, Article title, Journal, Volume, page numbers (year)
% Author, Book title, page numbers. Publisher, place (year)
% \end{thebibliography}

\end{document}